\documentclass[11pt]{amsart} 

\usepackage{amssymb,amsmath,amscd,amsxtra}
\usepackage{diagrams}
\usepackage[margin=4cm]{geometry}
\usepackage[all]{xy}

\newtheorem{theorem}{Theorem}[section]
\newtheorem{lemma}[theorem]{Lemma}

\newtheorem{proposition}[theorem]{Proposition}  
\newtheorem{corollary}[theorem]{Corollary} 

\theoremstyle{definition}
\newtheorem{remark}[theorem]{Remark} 
\newtheorem{example}[theorem]{Example} 

\newtheorem{definition}[theorem]{Definition}

\newcommand{\beq}{\begin{equation}}
\newcommand{\eeq}{\end{equation}}
\newcommand{\beqa}{\begin{eqnarray}}
\newcommand{\eeqa}{\end{eqnarray}}
\newcommand{\beaa}{\begin{eqnarray*}}
\newcommand{\ben}{\begin{eqnarray*}}
\newcommand{\eaa}{\end{eqnarray*}}
\newcommand{\een}{\end{eqnarray*}}

\newcommand{\ZZ}{\mathbb{Z}}
\newcommand{\CC}{\mathbb{C}}
\newcommand{\PP}{\mathbb{P}}

\renewcommand{\H}{\mathcal{H}}
\newcommand{\T}{\mathcal{T}}
\renewcommand{\O}{\mathcal{O}}
\newcommand{\F}{\mathcal{F}}

\newcommand{\cL}{\mathcal{L}}

\newcommand{\bX}{\overline{X}}
\newcommand{\bpi}{\overline{\pi}}

\begin{document} 

\title[Primitive forms and Frobenius structures on the Hurwitz spaces]
{Primitive forms and Frobenius structures on the Hurwitz spaces}
\author{Todor Milanov}

\address{Kavli IPMU (WPI), UTIAS, The University of Tokyo, Kashiwa, Chiba 277-8583, Japan}
\email{todor.milanov@ipmu.jp}

\begin{abstract}
The main goal of this paper is to introduce the notion of a primitive
form for a generic family of Hurwitz covers of $\mathbb{P}^1$
with a fixed ramification profile over infinity. We prove that primitive forms
are in one-to-one correspondence with semi-simple Frobenius
structures on the base of the family. Furthermore, we introduce the
notion of a polynomial primitive form and show that the corresponding
class of Frobenius manifolds contains the Hurwitz Frobenius manifolds
of Dubrovin. Finally, we apply our theory to investigate the relation
between the Eynard--Orantin recursion and Frobenius manifolds.
\end{abstract}
\maketitle

\setcounter{tocdepth}{2}
\tableofcontents

\section{Introduction}

\subsection{Motivation}
The main motivation behind this paper is the higher-genus
reconstruction of Givental and its applications to the classical
Riemann--Hilbert problem. Recall that if $X$ is a smooth projective
variety with semi-simple quantum cohomology, then Givental's higher
genus reconstruction expresses all Gromov--Witten (GW)
invariants of $X$ in terms of the  semi-simple Frobenius structure
underlying quantum cohomology. The reconstruction was proposed in
\cite{G1,G2} and proved in full generality by C. Teleman in
\cite{Te}. Using Givental's reconstruction as a definition we can
define the analogues of GW-invariants for any semi-simple
Frobenius manifold. The generating function of all genus invariants is
called the {\em total descendent potential} of the corresponding
semi-simple Frobenius manifold. On the other hand, a semi-simple Frobenius manifold can be
defined as a solution to a Riemann--Hilbert problem (see \cite{Du}). Usually, we have
a good knowledge of the corresponding monodromy data, while the
semi-simple Frobenius manifold depends on it in a highly
transcendental way. We are interested in the question 
whether we can express the invariants in terms of the
monodromy data in an algebraic way, e.g., via some explicit
recursions. 

One possible way to answer this question was proposed 
in our joint work \cite{BM}. We have defined a
$\mathcal{W}$-algebra (depending only on the monodromy data) and
proved that each state in the $\mathcal{W}$-algebra provides differential constraints for the 
total descendent potential of a simple singularity. The constraints
can be interpreted as recursion that uniquely determines the coefficients of
the total descendent potential (see \cite{LZ}), so we get a positive answer of the
question raised above for the Frobenius manifolds corresponding to
simple singularities. The ideas from
\cite{BM} are straightforward to generalize to any semi-simple
Frobenius manifold, but the problem is that it is very hard to find
states in the $\mathcal{W}$-algebra.

In the paper \cite{ML}, we
have tested, in the settings of $A_N$-singularity, a new idea  to construct states in the
$\mathcal{W}$-algebra based on the {\em topological recursion} (see \cite{EO,BE}).
Let us recall some basic settings for the recursion. The
starting pont is a triple $(\Sigma,x,y)$, where $x:\Sigma \to \PP^1$
is a branched covering and $y$ is a meromorphic function satisfying some
genericity assumptions. The Riemann surface $\Sigma$ is also known as
the {\em spectral curve}.  Using this data Eynard and Orantin have
proposed a recursion that produces a set of symmetric forms
\ben
\omega_{g,n}(q_1,\dots,q_n)\in T^*_{q_1}\Sigma\otimes \cdots\otimes
T^*_{q_n}\Sigma,\quad 2g-2+n>0,
\een
defined for all $(q_1,\dots,q_n)\in \Sigma^n$ such that $q_i$ is not a
ramification point of $x:\Sigma\to \PP^1$ and having at most finite order poles if some $q_i$ is a
ramification point (see \cite{EO} for precise definitions). We will
refer to the recursion as the topological recursion or the Eynard--Orantin (EO) recursion.
We are interested to find other
examples of semi-simple Frobenius manifolds for which the topological
recursion can be used to construct states in the
$\mathcal{W}$-algebra. This problem can be splited into two parts.
The first part is to classify all semi-simple Frobenius manifolds that correspond to
topological recursion. 
The second part is to determine whether the
topological recursion has a global contour formulation in a sense to
be clarified below. 

The problem of describing the correspondence between semi-simple
Frobenius manifolds and topological recursion was essentially solved in the recent
paper \cite{DNOPS2}. More precisely,
the authors proved that the Hurwitz Frobenius manifolds (see
\cite{Du}) correspond to an EO recursion, provided we relax the
condition that $y$ is a meromorphic function to the condition that the 1-form $dy$ is
holomorphic in a neighborhood of the finite ramification points. 
Although it is not stated explicitly in \cite{DNOPS2}, using the
results of this paper it is not very hard to prove 
that the Hurwitz Frobenius manifolds are the only
Frobenius manifolds that correspond to an Eynard--Orantin recursion
(with the relaxed condition on $y$). In particular, the set of Frobenius
manifolds that correspond to the original EO recursion (as defined in
\cite{EO}) is contained
in the set of Hurwitz Frobenius manifolds. More precisely, recall that
a Hurwitz Frobenius manifold (see \cite{Du}) depends on the choice of a {\em primary
differential} $\phi$ on $\Sigma$ (see also Section
\ref{sec:p-dif}). In order, to have the original EO recursion we have to require
that $\phi=dy$ 
for some meromorphic function $y$. Therefore, from the point of view
of the original EO recursion, it is still an open problem to
classify primary differential that are exact as meromorphic forms.

The starting point of the current paper is the
observation that many of the constructions in \cite{DNOPS2} have a
very natural interpretation in terms of K. Saito's theory of primitive
forms \cite{Sa}. Our main goal is to develop in a systematic way
the notion of primitive forms for families of Hurwitz covers. The main
outcome of our approach is that we were able to find an interesting
generalization of the topological recursion which allows us to extend
the correspondence described in \cite{DNOPS2} to include a wider class of
semi-simple Frobenius manifolds. We also found a very interesting identity
expressing the so-called {\em descendant correlators} in terms of
oscillatory integrals whose integrands are defined by the topological
recursion.

\subsection{Organization of the paper}
The main part of the paper is devoted to constructing Frobenius
structures on the base of a generic family of Hurwitz covers. Our
approach is based on K. Saito's theory of primitive forms (see
\cite{He, SaT}).  The general theory  is developed in Sections \ref{sec:results},
\ref{sec:pf}, \ref{sec:gb}, and \ref{sec:pf-classif}.  Many of the
results here can be obtained from the work of A. Douai and C. Sabbah
(see \cite{DSa1,DSa2}), which generalizes the Hodge theoretic approach
of M. Saito (see \cite{MSa}). Our approach has the advantage of
being elementary. We were able to prove the existence of primitive
forms without relying on Hodge theory thanks to a special set of
holomorphic forms that were introduced essentially in the work of
V. Shramchenko \cite{Shr}. 

Let us point out that if we allow an arbitrary primitive form, then we
can obtain any semi-simple Frobenius manifold. However, the analytic
properties of the underlying Frobenius manifold would not be captured
by the geometry of the underlying spectral curve. On the other hand,
our original goal is to understand $\mathcal{W}$-constraints (and
integrable systems)  using the geometry of the spectral curve. That is
why in Section \ref{sec:pf-poly}, we have proposed the notion of {\em polynomial primitive
  forms} and established several basic properties. It will be
interesting to classify all Frobenius manifolds that correspond to
polynomial primitive forms. We carry out this classification for
Frobenius manifolds of dimension 2.  We obtained a discrete 
set of Frobenius manifolds, which includes the Frobenius manifolds
corresponding to $A_2$-singularity and quantum cohomology of
$\mathbb{P}^1$. We also prove that the primary differentials (see
\cite{Du}) are polynomial primitive forms. 

In Section \ref{sec:tr-pf} we define the problem of comparing
semi-simple Frobenius manifolds and the topological recursion. We prove that
a semi-simple Frobenius manifold corresponds to a topological
recursion if and only if it is one of the Dubrovin's Hurwitz Frobenius manifolds.  

Finally, in Section \ref{sec:desc}  we express the descendant
correlators in terms of oscillatory integrals, whose integrands are
precisely the forms defined by the EO recursion. Similar formulas were
derived in the settings of equivariant mirror symmetry for
$\mathbb{P}^1$ in \cite{FLZ}. Motivated by these
results we propose to think of the EO recursion as defining twisted de
Rham cohomology classes. We found a generalization of the EO
recursion that deserves a further investigation from the point of view
of mirror symmetry.

\subsection{Future directions}

Let us
point out that the EO recursion is defined in terms of sums of
residues over the ramification points of $x:\Sigma\to \PP^1$ and it is
local in a sense that the forms whose residues are computed are defined
only locally near each ramification point. Nevertheless, if we require
that $y$ satisfies some extra properties, then Bouchard and Eynard
proved in \cite{BE} that the sum of the local residues of the EO recursion can be
replaced with a contour integral of a global meromorphic form on
$\Sigma$. This special class of EO-recursions that admits a global
contour integral presentation will be called Bouchard--Eynard (BE)
recursions. It is the BE recursion that was
identified with $\mathcal{W}$-constraints in \cite{ML}. Although the EO
recursion can be generalized in various ways, for our
purposes, the interesting generalizations are the ones that admit a
global contour formulation. It is very interesting to investigate the
existence of a BE-type recursion for the generalization of the
topological recursion proposed in this paper. In fact the existence of a BE recursion
is an open problem even for the generalization of the EO recursion
(proposed in \cite{DNOPS2}) corresponding to the Hurwitz Frobenius
manifolds.  

Our classification in the rank 2 case, shows that most of the
Frobenius manifolds corresponding to finite reflection groups do not
correspond to polynomial primitive forms. Nevertheless, based on our
recent work in \cite{Mi2}, it is clear how to generalize the
methods of this paper in order to include the case of Frobenius
manifolds corresponding to finite reflection groups. Namely, we have
to allow for the spectral curve to be
an orbifold Riemann surface. It will be interesting to search for a
corresponding generalization of the EO and BE recursions.

{\bf Acknowledgements.}
I would like to thank Ilya Karzhemanov and especially Tomoyuki Abe for
very useful discussions on sheaf cohomology. Also, I would like to
thank Sergey Lando for e-mailing me a draft of his work in progress on
a related subject. This work is partially supported by JSPS Grant-In-Aid 26800003 
and by the World Premier International Research Center Initiative (WPI
Initiative),  MEXT, Japan. I would like also to thank the mathematical research
institute MATRIX in Australia where part of this research was
performed.

\section{Statement of the main results}\label{sec:results}

Let us begin by introducing the main ingredients of our construction
and stating the results that we would like to prove.

\subsection{Branched coverings with a fixed ramification
  profile over infinity}\label{sec:br_cov}
Let $f:\Sigma\to \PP^1$ be a branched covering. Recall that a point $p\in
\Sigma$ is called a {\em ramification point} if $df(p)=0$ and that $u\in
\PP^1$ is called a {\em branch point} if $f^{-1}(u)$ contains a
ramification point. We will be interested only in branched
coverings $f$ that are generic in the following sense: if $u\in \PP^1\setminus{\{\infty\}}$ is a branch
point, then the fiber $f^{-1}(u)$ contains only 1 ramification point $p$
and the local degree of $f$ near $p$ is 2, i.e., there is a local
coordinate $t$ on $\Sigma$ near $p$, s.t., $f(q)=u+\frac{1}{2}t(q)^2$ for
all $q\in \Sigma$ sufficiently close to $p$.  

Let $u_i^\circ$ $( 1\leq i\leq N)$ be the set of finite branch
points and  $\infty_i$ $(1\leq i\leq d)$ be the points of the fiber
$f^{-1}(\infty)$. The Riemann--Hurwitz formula yields
\ben
N = 2g-2+d +\sum_{i=1}^d m_i,
\een 
where $g$ is the genus of $\Sigma$ and $m_i$ $(1\leq i\leq d)$ is the local
degree of $f$ near the point $\infty_i$ $(1\leq i\leq d)$. Let us
choose a reference point $u^\circ_0$ on
$\CC\setminus{\{u^\circ_1,\dots,u^\circ_N\}}$. Then we have a
monodromy representation
\ben
\rho:\pi_1(\CC\setminus{\{u^\circ_1,\dots,u^\circ_N\}},u^\circ_0)\to
S_m,
\quad m:=m_1+\cdots +m_d,
\een
where $S_m$ is the group of permutations of the points in the fiber
$f^{-1}(u^\circ_0)$.  

Let us denote by $B$ the universal cover of the configuration space 
\ben
\{u\in (\PP^1)^{N+1}\ |\ u_i\neq u_j \mbox{ for } i\neq j,\mbox{ and }
u_{N+1}=\infty\}.
\een
The points of $B$ consists of pairs $(u,[\gamma])$ of a point $u$ in
the configuration space and the homotopy class of a path (in the
configuration space) from
$u^\circ:=(u_1^\circ,\dots,u_N^\circ)$ to $u$.
It is known that $B$ is a contractible Stein manifold. Put
\ben
D_i=\{(\widetilde{u},\lambda)\in B\times \PP^1\ |\ \lambda=u_i\},\quad 1\leq
i\leq N+1,
\een
where $\widetilde{u}=(u,[\gamma])\in B$ and $u_i$ is the $i$-th
component of $u$. Note that the projection on the first factor
\ben
\pi:B\times \PP^1\setminus{\{D_1\cup \cdots \cup D_{N+1}\}} \to B
\een
is a smooth fibration with fibers diffeomorphic to $\CC
\setminus{\{u_1^\circ,\dots,u_N^\circ\}}$. The latter is identified
with the fiber of $\pi$ over the point $(u^\circ,[1])$. The long exact sequence of
homotopy groups and the contractibility of $B$ imply that the natural
inclusion of the fiber induces an isomorphism 
\ben
\pi_1(\CC\setminus{\{u^\circ_1,\dots,u^\circ_N\}})\cong 
\pi_1(B\times \PP^1\setminus{\{D_1\cup \cdots \cup
  D_{N+1}\}}),
\een
where we suppressed the base points in the above notation, but they are
uniquely determined from $u_0^\circ$ and the embedding of
the fiber in the total space of the fibration. In particular, the
monodromy representation $\rho$ from above
defines a representation 
\ben
\rho: \pi_1(B\times \PP^1\setminus{\{D_1\cup \cdots \cup
  D_{N+1}\}}) \to S_m.
\een
Let 
\ben
\widetilde{\varphi}:
\widetilde{X}\to \PP^1\setminus{\{D_1\cup \cdots \cup
  D_{N+1}\}}
\een
be the degree-$m$ covering whose monodromy representation is
$\rho$. Using the Riemann's extension theorem we can extend
$\widetilde{\varphi}$ to a branched covering (for more details see
\cite{Fu}, Proposition 1.2)
\ben
\overline{\varphi}: \overline{X}\to B\times \PP^1. 
\een

\subsection{Saito structure}\label{sec:Saito_str}
Let us denote by $\Omega^{p}_{\overline{X}/B}$ ($p=0,1$) the sheaf of relative holomorphic
$p$-forms on $\overline{X}$ relative to 
\ben
\overline{\pi}:=\operatorname{pr}_B\circ\,  \overline{\varphi}:\overline{X}\to B.
\een 
Put $D_\infty:=\overline{\varphi}^{-1}(B\times \{\infty\})$  and let 
\beq\label{dfs}
\Omega^{p}_{\overline{X}/B}(m):=\Omega^p_{\bX/B}\otimes
\O_{\bX}(m\,D_\infty),\quad m\geq 0
\eeq
be the sheaf of relative $p$-forms
that are holomorphic except for a possible pole of order at most $m$ along the divisor 
$D_\infty$.  The set \eqref{dfs} is naturally a filtered directed
system. Let us define 
\ben
\Omega^{p,\infty}_{\overline{X}/B} :=\varinjlim_m\ 
\Omega^{p}_{\overline{X}/B}(m),\quad p=0,1.
\een
Put $X:=\overline{X}\setminus{D_\infty}$ and let
\ben
\varphi:=\overline{\varphi}|_X\ :X\to B\times \CC.
\een 
Following K. Saito \cite{Sa}, we define
the {\em twisted de Rham cohomology} $\mathcal{H}$ of the holomorphic function 
\ben
F:= \operatorname{pr}_\CC\circ \varphi\ :X\to \CC
\een
by the following formula 
\ben
\mathcal{H}:=
\overline{\pi}_* \Omega^{1,\infty}_{\overline{X}/B}[z]/
(zd_{\overline{X}/B} +d_{\overline{X}/B}F\wedge) 
\overline{\pi}_* \Omega^{0,\infty}_{\overline{X}/B}[z]),
\een
where $d_{\overline{X}/B}$ is the relative de Rham differential. 
This is a sheaf on $B$ and the completion $\widehat{\H}:=\H\otimes\CC[\![z]\!]$ turns out to be isomorphic to
$\mathcal{T}_B[\![z]\!]$, where $\mathcal{T}_B$ is the sheaf of
holomorphic vector fields on $B$. 
Every section $\omega$ of $\mathcal{H}\otimes \CC[z,z^{-1}]$ defines a family of oscillatory integrals
\ben
\Gamma \mapsto (-2\pi z)^{-1/2} \int_\Gamma e^{F/z} \omega,
\een
where the integration cycle $\Gamma=\Gamma_{u,z}$ depends on the choice of
parameters $(u,z)\in B\times \CC^*$ and it is an element of the
relative homology group   
\ben
H_1(X_u,\operatorname{Re}(F/z)\ll 0;\CC):=
\varprojlim_{m\in \mathbb{Z}_{>0}} \
H_1(X_u,\operatorname{Re}(F/z)<-m;\CC)\cong \CC^N.
\een
The above homology groups form a vector bundle on $B\times
\mathbb{C}^*$ equipped with a flat Gauss--Manin connection. The sheaf
$\mathcal{H}\otimes \CC[z,z^{-1}]$ has an induced connection $\nabla$, which is also called Gauss--Manin
connection, s.t., 
\ben
\int_\Gamma e^{F/z}\nabla_v\omega = v \int_\Gamma e^{F/z}\omega,
\een
where $v\in \mathcal{T}_B$ is a vector field and $\Gamma$ is any {\em flat}
family of semi-infinite cycles. Using the above formula with
$v=\partial_z$ we can also extend $\nabla$ in the $z$-direction.

The sheaf $\mathcal{H}$ is equipped also with a non-degenerate pairing $K$,
which as we will see later on coincides with K. Saito's higher residue
pairing. Namely,
\ben
K:\mathcal{H}\otimes_{\mathcal{O}_B} \mathcal{H}\to \mathcal{O}_B[\![z]\!]z
\een
is defined by
\beq\label{hrp}
K(\omega_1,\omega_2) = \frac{1}{2\pi\sqrt{-1}}\, \sum_{i=1}^N \int_{\Gamma_i}
e^{F/z} \omega_1\  \int_{\Gamma^\vee_i}
e^{-F/z} \omega^*_2,
\eeq
where ${ }^*$ denotes the involution $z\mapsto -z$,
$\{\Gamma_i\}_{i=1}^N$ is any basis of flat semi-infinite cycles,
$\{\Gamma^\vee_j\}_{j=1}^N$ is a dual basis, and we are identifying
$\Gamma_j^\vee$ with homology cycles via the intersection pairing
\ben
H_1(X_u,\operatorname{Re}(F/z)\ll 0;\mathbb{Z})\times 
H_1(X_u,\operatorname{Re}(F/z)\gg 0;\mathbb{Z}) \to \mathbb{Z}.
\een
Note that this is a perfect pairing, i.e., the intersection matrix in
an appropriate integral basis is non-degenerate with determinant $\pm 1$. 
\begin{remark}
The above homology groups and the intersection pairing can be computed
via Morse theory. Namely, for fixed $(u,z)\in B\times \CC^*$ the
function $g:=\operatorname{Re}(F/z):X_u\to \mathbb{R}$ is a real Morse
function. Note that the critical points of $g$ are by definition the
ramification points of the covering $X_u\to \CC$. In particular, the
cycles $\Gamma_i$ (resp. $\Gamma_i^\vee$) can be constructed as the
gradient trajectories of $-g$ (resp. $g$) that flow out of the $i$-th
critical point.
\end{remark}
\begin{remark} If $\F$ is a sheaf of Abelian groups on a manifold $M$, then we denote by
  $\F[z]$ the sheaf on $M$ obtained
  by sheafification of the presheaf  
$
V\mapsto \F(V)[z].
$
Note that if $V$ has finitely many connected components, then
$\F[z](V)=\F(V)[z]$. Similarly we define $\F[\![z]\!]$ and $\F[z,z^{-1}]$.
\end{remark}
\begin{remark}
The notation $\H\otimes \CC[\![z]\!]$ (resp. $\H\otimes
\CC[z,z^{-1}]$) is for the sheaf defined in the same way as $\H$
except the we replace the sheaves $\bpi_*\Omega^{p,\infty}_{\bX/B}[z]$
with $\bpi_*\Omega^{p,\infty}_{\bX/B}[\![z]\!]$
(resp. $\bpi_*\Omega^{p,\infty}_{\bX/B}[z,z^{-1}]$). 
\end{remark}
 
Both the connection $\nabla$ and the higher-residue pairing $K$ extend
uniquely to the completion $\widehat{\H}$. The data
$(\widehat{\mathcal{H}},K,\nabla)$ will be called {\em Saito
  structure}. It allows us to introduce the notion of a {\em primitive form} in
$\widehat{\mathcal{H}}$ (see Section \ref{sec:pf}). 

\subsection{Primitive forms and Frobenius structures}

Let us outline how to construct the Frobenius
structure using a primitive form $\omega\in \widehat{\mathcal{H}}(V)$,
where $V\subseteq B$ is an open subset (see \cite{He,SaT}). The critical values of $F$, i.e., the branch points
$u_i$, $1\leq i\leq N$ will turn out to be canonical coordinates so
the multiplication is 
\ben
\partial_{u_i}\bullet \partial_{u_j} := \delta_{i,j}\, \partial_{u_j},
\een 
the Frobenius pairing is
\ben
(\partial_{u_i},\partial_{u_j}) := 
K^{(0)}(z\nabla_{\partial_{u_i}}[\omega], z\nabla_{\partial_{u_j}}[\omega]),
\een
and the Euler vector field is 
\ben
E=u_1\partial_{u_1}+\cdots+u_N\partial_{u_N}.
\een
Using this data we define a connection $\widetilde{\nabla}$ on the vector bundle
$\operatorname{pr}_V^*TV$, where $\operatorname{pr}_V:V\times \CC^*\to
V$ is the projection. Namely 
\ben
\widetilde{\nabla}_{\partial_{u_i}}:= \nabla^{\rm L.C.}_{\partial_{u_i}} +
z^{-1} \partial_{u_i}\bullet,\quad 1\leq i\leq N,
\een
and
\ben
\widetilde{\nabla}_{\partial_z}:= \partial_z -z^{-1} \theta - z^{-2} E\bullet,
\een
where $\nabla^{\rm L.C.}$ is the Levi--Civita connection of the
Frobenius pairing, for a vector field $v\in \mathcal{T}_B(V)$ we denoted
by $v\bullet\in \operatorname{End}(\mathcal{T}_B)(V)$ the linear operator
of Frobenius multiplication by $v$, and $\theta\in
\operatorname{End}(\mathcal{T}_B)(V)$ is defined by
\beq\label{def:theta}
\theta(v) = \nabla^{\rm L.C.}_{v}E - \Big(1-\frac{D}{2})\, v,
\eeq
where the constant $D$, known as {\em conformal dimension}, is such that $\theta$ is skew-symmetric with
respect to the Frobenius pairing. 
\begin{remark}
We do not assume that $\theta$ is diagonalizable. 
\end{remark}
Given a primitive form $\omega\in \widehat{\mathcal{H}}(V)$, we
can construct a period isomorphism
\ben
\mathcal{T}_B(V)[\![z]\!]\cong \widehat{\mathcal{H}}(V),\quad v\mapsto z\nabla_{v} \omega.
\een
The axioms of a primitive form imply that the period isomorphism
intertwines the connection $\widetilde{\nabla}$ and the Gauss--Manin
connection, i.e., 
\ben
z\nabla_v \nabla_w (z^{-1/2}\omega) = 
z\nabla_{\widetilde{\nabla}_v  w} ( z^{-1/2}\omega) 
\een
for all $v,w\in \mathcal{T}_{B\times \CC^*}$. In particular, the
flatness of the Gauss--Manin connection implies that
$\widetilde{\nabla}$ is flat. It remain only to recall that 
the axioms of a Frobenius structure
are equivalent to the flatness of the connection $\widetilde{\nabla}$ (see
\cite{Du}). 
\begin{theorem}\label{t1}
Let $V\subset B$ be an open subset, then every semi-simple Frobenius structure on
$V$, such that, the canonical vector fields coincide with
$\partial_{u_i}$, $1\leq i\leq N$, is the Frobenius structure
associated to a primitive form in $\widehat{\mathcal{H}}(V)$.
\end{theorem}

Theorem \ref{t1} might look a bit surprising at first, because it
essentially says that the branch covering $f:\Sigma\to \PP^1$ that we
used to set up the entire theory is irrelevant. The reason for this is
that if we work with the completion $\widehat{\H}$, then there is an
excision principle (see Proposition \ref{per-iso}, b)) which allows us to
replace the complex manifold $\bX$ with a tubular neighborhood of the
relative critical variety 
\ben
C=\{p\in X\ |\ d_{X/B}F(p)=0\} .
\een
The variety $C\cong B^N$, which explains why the choice of the branch
covering is irrelevant. In order to obtain Frobenius manifolds that
depend on the branch covering $f:\Sigma\to \PP^1$ in an essential way,
we propose to work with primitive forms $\omega\in \H$. Such forms can
be represented by holomorphic forms that depend polynomially on
$z$. We refer to them as {\em polynomial primitive forms}. 
\begin{theorem}\label{t2}
The primary differentials of Dubrovin are polynomial primitive forms.
\end{theorem}
The primary differentials are divided into 5 types. Except for type
IV, they are elements of $\Omega^{1,\infty}_{\bX/B}(\bX)$, so each
primary differential determines naturally an element in $\H(B)$. The
primary differentials of type IV are multivalued, but they have
analytic branches in a neighborhood of $C$. Using the excision
principle from Proposition \ref{per-iso}, b), we can construct
cohomology classes in $\widehat{\H}$, which turn out to be in $\H$. 
As a byproduct of our argument proving Theorems \ref{t1} and \ref{t2}
we got the following important result. 
\begin{corollary}\label{c1} The sheaf $\H$ is a free $\O_B[z]$-module of
  rank $N$. Moreover, there exists an $\O_B[z]$-basis
  $\{\omega_i\}_{i=1}^N\subset \H(B)$ such that $K(\omega_i,\omega_j)=z\delta_{ij}$ .
\end{corollary}

The formalism of primitive forms allows us to answer the question when
does an EO recursion define a Frobenius manifold. More precisely,
following \cite{Mi1} we introduce the notion of a local EO recursion. It
is proved in \cite{Mi1} that every semi-simple Frobenius manifold
provides a solution to a local EO recursion. We will prove that every
EO recursion also provides a solution to a local EO
recursion. Therefore, we can define the notion of a semi-simple
Frobenius structure that is a solution to an EO recursion. 
\begin{theorem}\label{t3}
A semi-simple Frobenius manifold is a solution to an EO recursion if
and only if it is a Hurwitz Frobenius manifold in the sense of Dubrovin,
i.e., the corresponding primitive form is a sum of homogeneous primary
differentials of the same degree. 
\end{theorem}

\section{The period map}\label{sec:pf}

The goal of this section is to introduce the notion of a primitive
form in $\widehat{\mathcal{H}}$. Our settings are slightly different from the
original ones (see \cite{Sa}), but the necessary modifications are
straightforward. 
We are going to use the following notation. Let
$\pi:=\overline{\pi}|_X:X\to B$ and $\Omega^1_{X/B}$ be the sheaf of sections of the relative
cotangent bundle $T^*_{X/B}:=T^*X/\pi^* T^*B$. The fiber of $T^*_{X/B}$
at a point $p\in X$ is $T_p^*X_{\pi(p)}$. If $\omega\in \Omega^1_{X/B}(U)$
and $p\in U$, then we denote by $\omega(p)\in T_p^*X_{\pi(p)}$ the value of the
section $\omega$ at $p$.

\subsection{The Kodaira--Spencer isomorphism}\label{sec:KS-iso}
Let us define the following sheaves of vector fields on $\bX$:
\ben
\mathcal{T}_{\overline{X}}^\infty=\varinjlim_m \ \T_{\bX}(m),
\quad
\mathcal{T}_{\overline{X}/B}^\infty=\varinjlim_m \ \T_{\bX/B}(m),
\een 
where  $\mathcal{T}_{\overline{X}}$ (resp. 
$\mathcal{T}_{\overline{X}/B}$)  is the sheaf of 
holomorphic (resp. holomorphic relative) vector fields and
for every sheaf $\F$ of $\O_{\bX}$-modules we put
$\F(m):=\F\otimes \O_{\bX}(m D_\infty)$. 
\begin{lemma}\label{vf-lift}
Let $U\subset B$ be an open Stein subset. Then every holomorphic
vector field $v\in \T_B(U)$ admits a lift $\widetilde{v}\in
\T_{\bX}^\infty(\bpi^{-1}(U))$, such that $v=\bpi_* \widetilde{v}$.
\end{lemma}
\proof
Recall that we have the following exact sequence 
\ben
0\to 
\mathcal{T}_{\overline {X}/B}(m) \to 
\mathcal{T}_{\overline{X}}(m)\to
\overline{\pi}^*\mathcal{T}_{B}\otimes
\mathcal{O}_{\overline{X}}(m)\to 
0,
\een
for every $m\in \ZZ$. 
Let us choose $m>4g-4$. Note that the
cohomology groups $H^1(\bX_u,\T_{\bX_u}(m))=0$, because $\T_{\bX_u}(m)$ is a line
bundle on $\bX_u$ of degree $2-2g+m$ and all line bundles of degree
more than the degree of the canonical bundle have vanishing higher
cohomologies. Furthermore, the map $\bpi:\bX\to B$ is a proper regular
map between complex manifolds, where regular means that 
$d_p\bpi:T_p\bX\to  T_{\bpi(p)}B$ 
is surjective for all $p\in \bX$. We get that the higher direct image sheaf
$R^1\bpi_*(\T_{\bX}(m))=0$, because it is a coherent sheaf on $B$ with vanishing
fibers (see \cite{GR}, Theorem 10.5.5). We get the following exact
sequence of $\O_B$-modules:  
\ben
0\to 
\overline{\pi}_*\mathcal{T}_{\overline {X}/B}(m) \to 
\overline{\pi}_*\mathcal{T}_{\overline{X}}(m)\to
\mathcal{T}_{B}\otimes
\overline{\pi}_*\mathcal{O}_{\overline{X}}(m)\to 
0.
\een
On the other hand, we have an injective map
\ben
0\to \mathcal{T}_B(U)\to (\mathcal{T}_B\otimes
\overline{\pi}_*\mathcal{O}_{\overline{X}}(m))(U) \cong
\mathcal{T}_{\overline{X}}(m)(\overline{\pi}^{-1}(U)) /\mathcal{T}_{\overline{X}/B} (m)(\overline{\pi}^{-1}(U)),
\een
where the isomorphism holds, because $U$ is Stein and
$\overline{\pi}_*\mathcal{T}_{\overline {X}/B}(m)$ is a coherent
sheaf. 
In particular, every holomorphic vector field $v\in \mathcal{T}_B(U)$
admits a lift $\widetilde{v}\in \T_{\bX}(m)(\bpi^{-1}(U)) \subset
\mathcal{T}^\infty_{\overline{X}}(\bpi^{-1}(U)) $.
\qed

\smallskip

If $v\in \T_B(U)$, then we cover $U=\cup_i U_i$ with open Stein
subsets, choose a lift $\widetilde{v}_i$ of $v|_{U_i}$ for every $i$ (see Lemma
\ref{vf-lift}), and define the section $\widetilde{v}(F)|_C\in
\pi_*\O_C(U)$ by gluing the 
sections $\widetilde{v}_i(F)|_C\in \pi_*\O_C(U_i)$. It is straightforward to
check that the gluing is possible and that the construction is
independent of the choices of lifts. Therefore we get a map (of
$\O_B$-modules)  
\ben
\mathcal{T}_B\to \pi_*\mathcal{O}_C,\quad 
v\mapsto \left. \widetilde{v}(F)\right|_C,
\een
which will be called the {\em Kodaira--Spencer map}. 
Note that $C$ is a disjoint union of $N$ analytic varieties, each
isomorphic to $B$ via the map $\pi:=\overline{\pi}|_X:X\to B$. Using this fact, it is
straightforward to check that the Kodaira--Spencer map is an
isomorphism and that the induced multiplication $\bullet$ on $\mathcal{T}_B$
takes the form 
\ben
\partial_{u_i}\bullet \partial_{u_j}= \delta_{ij}\partial_{u_j}.
\een
Indeed, the connected components of $C$ are
\ben
C_i=\{p\in C\ |\ F(p)=u_i\},\quad 1\leq i\leq N.
\een
In a tubular neighborhood of $C_i$ we can choose a holomorphic
function $t_i$, s.t., 
\ben
F(p)=u_i+\frac{1}{2} t_i(p)^2.
\een
The image of $\partial_{u_i}$ under the Kodaira--Spencer map is the
direct sum of the $N$ functions
$\widetilde{\partial}_{u_i}F|_{C_j} = \delta_{ij}$, $1\leq j\leq N$. The above
statements are obvious. 

\subsection{The Gauss--Manin connection} 

Let $U\subset B$ be an open connected Stein subset. Note that 
$\bpi^{-1}(U)\cap D_\infty$ has finitely many (in fact $d$) connected
components.
\begin{lemma}\label{h0-colimit}
If $V\subset \bX$ is an open subset such that $V\cap D_\infty$ has
finitely many connected components, then 
\ben
H^0\Big(V\, ,\, \varinjlim_m\ \Omega^p_{\bX/B}(m)\Big) = 
\varinjlim_m\  
H^0\left(V\,,\,\Omega^p_{\bX/B}(m)\right) .
\een
\end{lemma}
\proof
Note that the RHS is canonically embedded into the LHS. We will prove
that the LHS is a subset of the RHS. 
By definition, a section of $\varinjlim \Omega^p_{\bX/B}$ over
the open subset $V$ corresponds to a collection of pairs
$\{(V_\alpha,s_\alpha)\}_{\alpha \in \mathcal{A}}$, such that $\{V_\alpha\}_{\alpha\in \mathcal{A}}$ is
an open covering of $V$, 
\ben
s_\alpha\in \varinjlim
H^0(V_\alpha,\Omega^p_{\bX/B}(m))=\bigcup_{m=0}^\infty H^0(V_\alpha,\Omega^p_{\bX/B}(m)),
\een 
and $s_\alpha|_{V_\alpha\cap V_\beta}= s_\beta|_{V_\alpha\cap V_\beta}
$. Let us denote by $m_\alpha$ the order of the pole of $s_\alpha$
along $V_\alpha\cap D_\infty$. If $V_\alpha$ and $V_\beta$ intersect
the same connected component of $V\cap D_\infty$, then
$m_\alpha=m_\beta$. Since the connected components of $V\cap D$ are
finitely many, there exists an integer $m>0$, such that $m_\alpha<m$
for all $\alpha$. Therefore $s_\alpha\in
H^0(V_\alpha,\Omega^p_{\bX/B}(m) )$ and since $\Omega^p_{\bX/B}(m)$ is
a sheaf, we get that there exists a section $s\in
H^0(V,\Omega^p_{\bX/B}(m) )$ such that
$s_\alpha=s|_{V_\alpha}$. This is exactly what we had to prove.
\qed
\begin{lemma}\label{form-lift}
Let $U\subset B$ be an open  connected Stein subset, then the quotient
map
\ben
\operatorname{rel}: 
H^0( \bpi^{-1}(U),\Omega^{1,\infty}_{\bX}) \to 
H^0( \bpi^{-1}(U),\Omega^{1,\infty}_{\bX/B}) 
\een
is surjective.
\end{lemma}
\proof
We have the following exact sequence
\ben
0\to 
\bpi^*\Omega^1_B\otimes \O_{\bX}(m) \to 
\Omega^1_{\bX}(m)\to 
\Omega^1_{\bX/B}(m) \to 0,
\een
for every $m\in \ZZ$. Note that $\Omega^1_B =\O_B^{\oplus N}$ is a
free sheaf of rank $N$, because the cotangent bundle $T^*B$ is
trivial.  Let us choose $m>2-2g$. Since $\bpi^*\Omega^1_B
=\O_{\bX}^{\oplus N}$, we get $R^1\bpi_*(\bpi^*\Omega^1_B\otimes
\O_{\bX}(m) )=0$ and the following exact sequence
\ben
0\to 
\Omega^1_B\otimes \bpi_* \O_{\bX}(m) \to 
\bpi_*\Omega^1_{\bX}(m)\to 
\bpi_*\Omega^1_{\bX/B}(m) \to 0.
\een 
Let us apply to the above sequence the functor $\varinjlim H^0(U,-
)$. The exactness of the sequence will be preserved, because
direct limits preserve exact sequences (of Abelian groups), the open subset $U$ is
Stein, and the sheaf $\Omega^1_B\otimes 
\bpi_* \O_{\bX}(m)$ is coherent. In particular, we get that the map 
\ben
\varinjlim_m\ H^0(\bpi^{-1}(U),\Omega^1_{\bX}(m))\to
\varinjlim_m\ H^0(\bpi^{-1}(U),\Omega^1_{\bX/B}(m)) \to 0
\een
is surjective. Using that $U$ is connected, we get that
$\pi^{-1}(U)\cap D_\infty$ has finitely many connected components. It
remains only to recall Lemma \ref{h0-colimit}.
\qed

\medskip

If $\omega\in
\bpi_*\Omega^{1,\infty}_{\bX/B}(U)$ is a relative holomorphic form,
then we denote by  
\ben
[\omega]:=\int e^{F/z} \omega
\een
the equivalence class of $\omega$ in $\mathcal{H}(U)$. We have the
following formula for the Gauss--Manin connection (see \cite{ACG})
\beq\label{GM-conn}
z\nabla_v[\omega] = [\operatorname{rel}\circ 
\iota_{\widetilde{ v}} ((zd_{\bX}+d_{\bX}F\wedge)\widetilde{\omega})],
\eeq
where $\widetilde{v}\in \T^\infty_{\bX}(\bpi^{-1}(U))$ and $\widetilde{\omega}\in
\Omega^{1,\infty}_{\bX}(\bpi^{-1}(U))$ are lifts of $v\in \T_B(U)$ and $\omega$, $\iota_{\widetilde{ v}}$ is contraction
by the vector field $\widetilde{ v}$, and $d_{\bX}$ is the de Rham
differential on $\bX$.  

Let us derive more explicit formula for the Gauss--Manin
connection. To begin with, note that \eqref{GM-conn} is still true if we choose the lift
$\widetilde{v}$, s.t., 
\ben
\widetilde{v}(F)=\iota_{\widetilde{v}}\, dF=0.
\een
The above condition, together with $\widetilde{v}(u_i)=v(u_i)$, $1\leq
i\leq N$, uniquely determines $\widetilde{v}$. Note that $\widetilde{v}$ might have a
pole of order 1 along $C$. Nevertheless such a lift works too. Let us also fix a lift of
$\omega$ and write it in the form
\ben
\widetilde{\omega}(p)= h(p)dF(p) + \sum_{j=1}^N h_j(p) du_j,\quad p\in \pi^{-1}(U)\setminus{C}.
\een
Since $\widetilde{\omega}$ extends to a holomorphic 1-form on
$\pi^{-1}(U)$, using  local coordinates $(u_1,\dots,u_N,t_i)$ near
$C_i$, s.t., $F=u_i+\frac{1}{2} t_i^2$, we get that the functions
$h(p)$ and $h_i(p)$ ($1\leq i\leq N$) satisfy the following conditions: $h(p)$ has a pole
of order at most 1 at $C_i$ for all $i$, $h_i(p)$ has a pole of order
at most 1 along $C_i$ and it is holomorphic along $C_j$ for $j\neq i$,
and $h(p)+h_i(p)$ is holomorphic at $p=p_i$. Note that $\omega=h
d_{X/B}F$, so only the functions $h_i$ $(1\leq i\leq N)$ depend on the choice of the
lift and the ambiguity of each $h_i$ is up to a holomorphic function on
$\pi^{-1}(U)$.  The formula for the Gauss--Manin connection takes the
form
\beq\label{ui-derivative}
z\nabla_{\partial_{u_i}} [\omega] = (-h_i(p) + z h_i^{(1)}(p)) d_{X/B}F(p),
\eeq
where $h_i^{(1)} = \widetilde{\partial}_{u_i} h
-(d_{X/B}h_i/d_{X/B}F)$. 

\subsection{The period isomorphism}\label{sec:per-iso}
The key result of this section can be stated as follows.
\begin{proposition}\label{per-iso}
Let $U\subset B$ be an open connected Stein subset and 
$$
\omega=\sum_{n=0}^\infty \omega_n z^n \in
\bpi_*\Omega^{1,\infty}_{\bX/B}(U)[\![z]\!]
$$ 
be such that $\omega_0(p)\neq 0$  for all $p\in C\cap \bpi^{-1}(U)$. 

a) The period map
\ben
\mathcal{T}_B(U)[\![z]\!]\cong \widehat{\mathcal{H}}(U),\quad v\mapsto z\nabla_{v} \omega
\een
is an isomorphism.

b) Suppose that $U=B$ and that $X'\subset X$ is an open subset that
contains $C$.  Then the natural restriction map
$\Omega^{i,\infty}_{\bX/B}(\bX)\to \Omega^i_{X/B}(X')$, $i=0,1$, 
induces an isomorphism 
\ben
\widehat{\mathcal{H}}(B) \cong 
\Omega^1_{X/B}(X')[\![z]\!]/(z
d_{X/B}+d_{X/B}F\wedge) \Omega^0_{X/B}(X') [\![z]\!] .
\een
\end{proposition}
\proof
a)
It is enough to prove that 
\ben
\mathcal{T}_B(U)[\![z]\!]\cong \bpi_*\Omega^{1,\infty}_{\bX/B}(U)[\![z]\!]/(z
d_{\bX/B}+d_{\bX/B}F\wedge) \bpi_*\Omega^{0,\infty}_{\bX/B}(U) [\![z]\!] .
\een
Since $\T_B [\![z]\!]$ is a sheaf, the above identity shows that the
RHS coincides with the space of holomorphic sections over $U$ of the
quotient sheaf $\widehat{\H}$. Let us consider only the case when
$\omega_k=0$ for $k>0$. The argument in the general case is the same,
but the notation is a bit more cumbersome.

First we prove that the period map is surjective. Let
\ben
\psi=\sum_{k=0}^\infty \psi_k z^k\in
\bpi_*\Omega^{1,\infty}_{\bX/B}(U)[\![z]\!]. 
\een 
We want to prove that there are sequences of holomorphic functions
$c_{k,i}\in \O_B(U)$ and $\eta_k\in \O^\infty_{\bX}(\bpi^{-1}(U))$,
$1\leq i\leq N$,  $k\geq 0$, such that  
\beq\label{surj-eqn}
\psi = 
\sum_{k=0}^\infty 
\Big(
(zd_{\bX/B}+d_{\bX/B}F\wedge) \eta_k z^k+ 
\sum_{i=1}^N c_{k,i} z^k\,
z\nabla_{\partial_{u_i}} [\omega] \Big).
\eeq
Comparing the coefficients in front of $z^k$ for $k>0$, we get
\beq\label{equation:zk}
\psi_k(p) = 
d_{\bX/B} \eta_{k-1}(p) + \Big( \eta_k(p) +
\sum_{i=1}^N
\Big( -h_i(p) c_{k,i}(u) + c_{k-1,i}(u) h_i^{(1)}(p) \Big) 
\Big) d_{\bX/B}F(p).
\eeq
Similarly, comparing the coefficients in front of $z^0$ we get 
\beq\label{equation:z0}
\psi_0 (p)= \eta_0(p) d_{\bX/B}F(p) - \sum_{i=1}^N
h_i(p) c_{0,i}(u) d_{\bX/B}F(p).
\eeq
First, we find $c_{0,i}$. Let us fix $u\in U$ and denote by $p_i\in
X_u$ the critical points of $F$. Let $(u_1,\dots,u_N,t_i)$ be the
local coordinates near $C_i$, s.t., $F=u_i+\frac{1}{2}
t_i^2$. Dividing both sides of \eqref{equation:z0} by $t_i(p)$ and
computing the residue at $p=p_i$ we get
\ben
\operatorname{res}_{p=p_i} \frac{\psi_0(p)}{t_i(p)} =
-c_{0,i}(u)\operatorname{res}_{p=p_i} h_i(p)\, dt_i = c_{0,i}(u) \operatorname{res}_{p=p_i} h(p)\,dt_i.
\een
Since $\omega(p_i) = dt_i\, \operatorname{res}_{p=p_i} h(p)\,dt_i\neq
0$, we can solve uniquely for $c_{0,i}(u)$. This choice of $c_{0,i}$
makes the relative 1-form
\ben
\psi_0 (p)+ \sum_{i=1}^N
h_i(p) c_{0,i}(u) d_{\bX/B}F(p)
\een
vanish for $p\in C$, which implies that we can write it uniquely in
the form $\eta_0(p)d_{\bX/B}F(p)$ for some $\eta_0\in
\O_{\bX}^\infty(\bpi^{-1}(U))$. Assuming that we have found $c_{k',i}$ and
$\eta_{k'}$ for all $k'=0,1,\dots,k-1$, then applying the above
argument to equation \eqref{equation:zk}, we find first $c_{k,i}$ and
then $\eta_k$. This completes the proof of the surjectivity. 

To prove the injectivity, it is enough to verify that if $\psi=0$,
then the above algorithm gives $c_{k,i}=0$. This is straightforward to
do, which completes the proof of part a). 

b) Note that if we assume that $\psi\in \Omega^1_{\bX/B}(X')[\![z]\!]$
is any form, then the equation \eqref{surj-eqn} 
still has a unique solution $(\eta_k,c_{k,i})$, $k\geq 0,1\leq i\leq
N$, with $\eta_k\in \Omega^0_{\bX/B}(X')$ and $c_{k,i}\in \O_B(B)$. The
existence and uniqueness of the solution of \eqref{surj-eqn} is
equivalent to the statement in part b).  
\qed

Let us assume that $V\subset B$ is an open subset and that $\omega\in
\widehat{\H}(V)$ is a cohomology class. We would like to formulate a
condition on $\omega$ under which the period map 
\beq\label{V:per-iso}
\T_V[\![z]\!]\to \widehat{\H}|_V,\quad v\mapsto z\nabla_v \omega
\eeq
is an isomorphism. If $U\subset V$ is an open connected Stein subset,
then according to Proposition \ref{per-iso} $\omega|_U$ can be
represented by a holomorphic form
\ben
\omega_U=\sum_{n=0}^\infty \omega_U^{(n)} z^n,\quad \omega_U^{(n)}\in
H^0(\bpi^{-1}(U), \Omega^{1,\infty}_{\bX/B})
\een
We require that $\omega^{(0)}_U(p)\neq 0$ for all $p\in
\bpi^{-1}(U)\cap C$. Note that for each $p\in \bpi^{-1}(V)\cap C$ the
value of $\omega_U^{(0)}(p)$ is independent of the choices of a Stein
neighborhood $U$ of $p$ and a form $\omega_U$ representing the
cohomology class $\omega|_U$. If this condition on $\omega$ is satisfied, then we
will say that the leading order term of $\omega$ 
is a {\em volume form along $\bpi^{-1}(V)\cap C$}. Since open
connected Stein neighborhoods form a basis for the topology of $B$, we
get from Proposition \ref{per-iso}, Part a) that the period map
\eqref{V:per-iso} is an isomorphism for every 
$\omega\in \widehat{\H}(V)$ whose leading order term is a volume form
along $C$. 

\begin{lemma}\label{hv-form}
There exists  $\omega\in
\Omega^{1,\infty}_{\bX/B}(\bX)$, s.t., $\omega(p)\neq 0$ for all $p\in C$. 
\end{lemma} 
\proof
Note that the divisor $D_\infty=\sum_{i=1}^d m_i \infty_i$. We have
the following exact sequence
\ben
0 \rTo  \O_{\bX}(k\, D_\infty)  \rTo^{dF} 
\Omega^1_{\bX/B}(kD_\infty +\sum_{i=1}^d(m_i+1)\infty_i) \rTo  \mathcal{Q}
\rTo  0,  
\een
where $dF$ denotes the map of multiplication by the meromorphic relative
1-form $d_{\bX/B}F$, which is holomorphic on $X$ and has a pole of
order $m_i+1$ at $\infty_i$. The sheaf $\mathcal{Q}$ is defined to be the
quotient of the preceding two sheaves, so that the sequence is
exact. Finally, the number $k$ is a sufficiently large integer, so
that $H^1(\bX_u, \O_{\bX_u}(kD_\infty)) = 0$ (e.g. any $k>2g-2$
works). 
Using the Riemann--Roch formula we get that 
\ben
\operatorname{dim}_\CC\, H^0(\bX_u, \O_{\bX_u}(kD_\infty))
=1-g+k \sum_{i=1}^d m_i
\een 
and 
\ben
\operatorname{dim}_\CC\, H^0(\bX_u,\Omega^1_{\bX_u}(kD_\infty
+\sum_{i=1}^d (m_i+1)\infty_i ) = 
g-1 +d+(k+1)\sum_{i=1}^d m_i
\een
are independent of $u$. 
Furthermore, using that $\bpi:\bX\to B$ is a regular map, we get that 
$R^1\bpi_*(\O_{\bX_u}(kD_\infty))=0$ and that 
\ben
0 \rTo  \bpi_*\O_{\bX}(k\, D_\infty)  \rTo^{dF} 
\bpi_*\Omega^1_{\bX/B}(kD_\infty +\sum_{i=1}^d (m_i+1)\infty_i) \rTo
 \bpi_*\mathcal{Q}  \rTo   0
\een
is an exact sequence of vector bundles. Note that the rank of
$\bpi_*\mathcal{Q}$ is $2g-2+d+\sum_{i=1}^d m_i = N$. Since $B$ is contractible
and Stein the vector bundles must be trivial.  Note that 
\ben
\bpi_*\mathcal{Q}(B)= 
H^0(\bX, \Omega^1_{\bX/B}(kD_\infty +\sum_{i=1}^d(m_i+1)\infty_i) )/
H^0(\bX,\O_{\bX}(k\, D_\infty) dF )
\een
Let us fix holomorphic forms 
\ben
\omega_i\in
H^0(\bX,\Omega^1_{\bX/B}(kD_\infty+\sum_{i=1}^d
(m_i+1)\infty_i)),\quad 1\leq i\leq N,
\een 
inducing a trivialization of the vector bundle $\bpi_*\mathcal{Q}$. Let us define a
square holomorphic matrix $\Phi(u)=(\Phi_{ij}(u))_{i,j=1}^N$ of size
$N$  by $\omega_i(p_j) = \Phi_{ij}(u) dt_j.$ We claim that
$\Phi(u)$ is invertible for every $u\in B$. Otherwise, there is
$u\in B$ and constants $c_i$, such that $\sum_i c_i \omega_i(p)$
vanishes for all $p=p_j$. However, this would imply that $\sum_i c_i
\omega_i(p)$  is proportional to $dF$, which would imply that the
projections of $\omega_i$ in the fiber over $u$ of the vector bundle
$\bpi_*\mathcal{Q}$ are linearly dependent. This contradicts the fact that the
forms induce a trivializing frame in every fiber of $\bpi_*\mathcal{Q}$.  Let
us define 
\ben
(c_1(u),\dots,c_N(u))^T:= \Phi(u)^{-1} \, (1,\dots,1)^T,
\een
then the form $\omega(p):=\sum_{i=1}^N c_i(u)\omega_i(p)$ would
satisfy $\omega(p_j) = 1\neq 0$ for all $j=1,2,\dots, N$. 
\qed
 
The existence of the form $\omega$ in Lemma \ref{hv-form} implies that
$\widehat{\H}$ is a free $\O_B[\![z]\!]$-module of rank $N$ and that the excision
principle in Proposition \ref{per-iso}, Part b) holds.

\section{Good basis}\label{sec:gb}

Following K. Saito, we introduce the so-called {\em higher residue pairing}
\ben
K:\widehat{\H}\otimes \widehat{\H} \to \O_B[\![z]\!]z
\een
and prove the existence of a {\em good basis} on $B$.
Recall that if $V\subset B$ is an open subset, then a set 
of cohomology classes $\{\omega_i\}_{i=1}^N\subset \widehat{\H}(V)$ is
called a good basis on $V$, if 
\begin{enumerate}
\item[(GB1)] $\{\omega_i\}_{i=1}^N$ is a $\O_B[\![z]\!]$-basis of
  $\widehat{\H}|_V$.
\item[(GB2)] The higher residues vanish, i.e.,
  $K(\omega_i,\omega_j)\in z\, \O_B(V)$.
\end{enumerate}

\subsection{Higher-residue pairing}
Let us define $K$ locally on every open connected Stein subset 
$U\subset B$. It is straightforward to check that the local
definitions can be glued. 
Let $\psi^{(i)}\in \Omega^{1,\infty}_{\bX/B}(\bpi^{-1}(U))[\![z]\!]$, $i=1,2$ be
arbitrary. Following K. Saito \cite{Sa} we define the higher residue pairing 
\ben
K(\psi^{(1)},\psi^{(2)})\in \O_B(U) [\![z]\!]z
\een 
as the sum of the residues 
\beq\label{hrp-2}
\sum_{i=1}^N \operatorname{res}_{p=p_i(u)}
\phi^{(1)}(p,z)\psi^{(2)}(p,-z)\, z
\eeq
where we fix $u\in U$, denote by $p_i(u)\in X_u$ the critical points of $F$,
and let $\phi^{(1)}\in
\Omega^{0,\infty}_{\bX/B}(\bpi^{-1}(U)\setminus{C})[\![z]\!]$ be such that
\ben
(zd_{\bX/B}+ d_{\bX/B}F\wedge) \phi^{(1)}(p,z) = \psi^{(1)}(p,z).
\een
Note that the above equation has a solution for all $p\in X\setminus{C}$.

It is easy to check that the pairing induces a pairing on
$\widehat{\H}$ satisfying the following properties.
\begin{enumerate}
\item[(HR1)]
$K(\psi^{(1)},\psi^{(2)})=-K(\psi^{(2)},\psi^{(1)})^*$, where $*$ is
the involution $z\to -z$.
\item[(HR2)]
$a(z) K(\psi^{(1)},\psi^{(2)}) = K(a(z)\psi^{(1)},\psi^{(2)}) =
K(\psi^{(1)},a(-z)\psi^{(2)})$ for all $a\in \O_B[z]$. 
\item[(HR3)] The Leibnitz rule holds
\ben
z v K(\psi^{(1)},\psi^{(2)}) = K(z\nabla_v\psi^{(1)},\psi^{(2)})-K(\psi^{(1)},z\nabla_v\psi^{(2)})
\een
for all $v\in \T_B$ and $v=z\partial_z$. 
\item[(HR4)]
Let $K^{(0)}(\psi^{(1)},\psi^{(2)})\in \O_B$ be the coefficient in
front of $z^1$ in $K(\psi^{(1)},\psi^{(2)})$ and $\psi^{(i)}_0(p)$ be the coefficient
of $\psi^{(i)}(p,z)$ in  front of $z^0$. Then 
\ben
K^{(0)}(\psi^{(1)},\psi^{(2)}) = \sum_{i=1}^N 
\operatorname{res}_{p=p_i(u)}
\frac{\psi^{(1)}_0(p) \psi^{(2)}_0(p)}{d_{X/B}F(p)}.
\een
\item[(HR5)] $K(\psi^{(1)},\psi^{(2)})\in \O_B(U) [\![z]\!]z$.
\end{enumerate}
\begin{proposition}\label{hrp-uniqueness}
The properties (HR1)--(HR5) determine the higher residue pairing uniquely.
\end{proposition}
\proof
Let $\omega\in \pi_*\Omega^{1,\infty}_{\bX/B}(U)$ be a form such that
$\omega(p)\neq 0$ for all $p\in C\cap \pi^{-1}(U)$. The existence of
such form is proved in Lemma \ref{hv-form}. 
We have the following formulas for the Gauss--Manin connection
\beq\label{GM-ui}
z\nabla_{\partial_{u_i}} z\nabla_{\partial_{u_j}} [\omega] =
\sum_{k=1}^N C_{ij}^k(u,z) z\nabla_{\partial_{u_k}}[\omega],
\eeq
where $C_{ij}^k(u,z)=\sum_{n=0}^\infty C_{ij;n}^k(u) z^n$ for some
holomorphic functions $C_{ij;n}^k\in \O_B(U)$ and similarly
\beq\label{GM-z}
z^2\nabla_{\partial_z} z\nabla_{\partial_{u_i}}[\omega] =\sum_{j=1}^N
C_i^j(u,z) z\nabla_{\partial_{u_j}}[\omega].
\eeq
where $C_i^j(u,z)=\sum_{n=0}^\infty C_{i;n}^j(u)z^n$ for some
holomorphic functions $C_{i;n}^j\in \O_B(U)$.

We are going to introduce several matrices. They will all be square
matrices of size $N\times N$. 
Let us denote by $K(u,z)$
the matrix with entries
$K_{ij}(u,z):=K(z\nabla_{\partial_{u_i}},z\nabla_{\partial_{u_j}})$, 
by $C_i(u,z)$ $(1\leq i\leq N)$ 
the matrix with entries $(C_i)_{kj}:= C_{ij}^k(u,z)$, and by $C_0(u,z)$ the
matrix with entries $(C_0)_{ij}= C_i^j(u,z)$. Put
$K(u,z)=\sum_{n=0}^\infty K_n(u) z^{n+1}$ and 
$C_i(u,z)=\sum_{n=0}^\infty C_{i;n}(u)z^n$ $(0\leq i\leq N)$. 
Property (HR4)
fixes $K_0(u)$, while a straightforward computation using formula
\eqref{ui-derivative} yields that $C_{i;0} = E_{ii}$ $(1\leq i\leq N)$
and $C_{0;0}=-\sum_{i=1}^N u_i E_{ii}$, where $E_{ij}$
denotes the matrix with only one non-zero entry, which is on position
$(i,j)$ and it is equal to 1. 

Using (HR2) and (HR3) we get 
\ben
z\partial_{u_i} K(u,z) = C_i(u,z)^T K(u,z) - K(u,z) C_i(u,-z)
\een
and 
\ben
z^2\partial_z K(u,z) = C_0(u,z)^T K(u,z)-K(u,z) C_0(u,-z).
\een
Let us assume that we have two pairings $K',K''$ satisfying the axioms
(HR1)--(HR5), then the matrix $K(u,z):=K'(u,z)K''(u,z)^{-1}$
satisfies the differential equations
\ben
 z\partial_{u_i} K(u,z) & = & [C_i(u,z)^T, K(u,z)],\quad 1\leq i\leq N,\\
 z^2\partial_z K(u,z)  & = & [C_0(u,z)^T K(u,z)].
\een
Arguing by induction on $n$, we are going to prove that $K_n=0$ for
all $n>0$. 
Comparing the coefficients in front of the powers of $z$, we get the
following system of equations
\ben
[K_0 ,E_{ii}]=[K_0 ,C_{0;0}]=0,\quad 1\leq i\leq N,
\een
\ben
\partial_{u_i} K_n  = [E_{ii},K_{n+1} ] + \sum_{m=1}^{n+1}
[C_{i;m} ^T, K_{n+1-m} ]  ,\quad n\geq 0,\quad 1\leq i\leq N,
\een
and 
\ben
(n+2)K_{n+1}  = [C_{0;0} ,K_{n+2} ]+\sum_{m=1}^{n+2}
[C_{0;m} ^T, K_{n+2-m}]  , \quad n\geq -1.
\een
The firs set of equations is trivially satisfied, because $K_0$ is the
identity matrix (here we used axiom (HR4), which implies that the leading terms of $K'$ and $K''$
are fixed and equal. If $K_1=\cdots=K_n=0$, then from
the 2nd set of equations we get that $[K_{n+1},E_{ii}]=0$ for all $i$,
so $K_{n+1}$ must be diagonal. Note that in the last equation $C_{0;0}$
and $K_{n+2-m}$ $(1\leq m\leq n+2)$ are diagonal matrices. Comparing
the diagonal entries, we get $K_{n+1}=0$.
\qed
\begin{remark}
The same argument can be used to prove the uniqueness of K. Saito's
higher residue pairing in the settings of singularity theory. 
\end{remark}
Note that the pairing on $\H$ defined by formula \eqref{hrp} satisfies
all axioms (HR1)--(HR5). Moreover, it can be
extended uniquely to a pairing on the completion $\widehat{\H}$, so
that the axioms  (HR1)--(HR5) still hold. 
\begin{corollary}
The higher residue pairing \eqref{hrp-2} coincides with the pairing on
$\widehat{\H}$ defined by \eqref{hrp}.
\end{corollary}
\proof
The only non-trivial part in the proof is to verify that the pairing
\eqref{hrp} satisfies (HR4). However, the verification reduces to
computing the leading order term of the stationary phase asymptotic of
the oscillatory integrals. A standard computation yields
\beq\label{spa}
\int_{\Gamma_i} e^{F/z} \omega \sim (-2\pi z)^{1/2} e^{u_i/z}
\left.\frac{\omega(p)}{dt_i(p)}\right|_{p=p_i}\Big(1+\cdots\Big), \quad
z\to 0,
\eeq
where $\omega$ is a holomorphic 1-form on $X_u$, $t_i$ is a local
coordinate on $X_u$ near the critical point 
$p_i$, s.t., $F(p)=u_i+\frac{1}{2}t_i(p)^2$. 
The leading order term of the pairing \eqref{hrp} becomes 
\ben
z\ \sum_{i=1}^N
\left.\frac{\omega_1(p)\omega_2(p)}{dt_i(p)^2}\right|_{p=p_i} = 
z\ \sum_{i=1}^N
\operatorname{res}_{p=p_i} \frac{\omega_1(p)\omega_2(p)}{dF(p)}.\qed
\een

\subsection{Construction of the good basis}\label{subsec:gb}

Let us fix a symplectic basis $\{\alpha_i,\beta_i\}_{i=1}^g\subset
H_1(\Sigma,\mathbb{Z})$. Since $B$ is contractible, we can use the
Gauss--Manin connection to construct a symplectic basis in
$H_1({\bX}_u;\mathbb{Z})$ for all $u\in B$. Recall that by construction, the ramification points over infinity
of the branched covering $\overline{X}\to B\times \PP^1$ provide
$d$ holomorphic sections $\infty_i:B\to \overline{X}$, $1\leq i\leq d$ of $\bpi:
\overline{X}\to B$, while the ramifications over the finite branch
points provide $N$ holomorphic sections $p_i:B\to X$, $1\leq i\leq
N$. In particular, the connected component $C_i$ of the relative
critical set $C$ of $F$ coincides with $p_i(B)$. Let us fix a
holomorphic function $t_i$ defined in a tubular neighborhood of $C_i$,
s.t., 
\ben
F(p) =: u_i + \frac{1}{2} t_i(p)^2,\quad 1\leq i\leq N,
\een
where $u_i=F(p_i\circ \pi(p))$ is the critical value of $F$
corresponding to the connected component $C_i$. The choice of each $t_i$ is
unique up to a sign.

For fixed $u\in B$, let us denote by $\omega_p(q)$ the unique meromorphic
1-form of the 3rd kind on $\overline{X}_u$ that has poles only at $q=p$ and
$q=\infty_1(u)$ and such that 
\ben
\oint_{q\in \alpha_i} \omega_p(q) = 0,\quad
\operatorname{res}_{q=p} \omega_p(q) = 
-\operatorname{res}_{q=\infty_1} \omega_p(q) =1.
\een 
Let us introduce the relative differential forms 
\ben
\omega_i(p) =
\operatorname{res}_{q=p_i} t_i(q)^{-1} \omega_p(q) \,
d_{\bX/B}F(p),\quad 1\leq i\leq N.
\een
If $p$ is sufficiently close to $C_i$, then using
$(u_1,\dots,u_N,t_i)$ as local coordinates, we get that 
\beq\label{local-pi}
\omega_i(p) = (-1+O(t_i))dt_i(p).
\eeq
In particular, we get that $\omega_i$ is holomorphic in a neighborhood
of $C_i$. Clearly, $\omega_i(p)$ is holomorphic in a tubular
neighborhood of $C_j$ for $i\neq j$, so $\omega_i\in
\Omega^1_{X/B}(X')$, where $X'$ is a tubular neighborhood of
$C$. 
Let us recall Proposition \ref{per-iso}, part b) with $X'$ being the
tubular neighborhood of $C$ introduced above and $\omega\in
\Omega^{1,\infty}_{\bX/B}(\bX)$ a form whose existence is guaranteed by Lemma
\ref{hv-form}. We get that the forms
$\omega_i\in \Omega^1_{X/B}(X')$ determine global sections $[\omega_i]\in
\widehat{\H}(B). $ 

The 1-forms $\omega_i(p)$ are multivalued for $p\in X$, because
$\omega_p(q)$ is multivalued:
\ben
\omega_p(q) = \int_{\infty_1}^p B(p',q),
\een
where $B(p',p'')$ is the so-called {\em Riemann's 2nd fundamental
  form} or {\em fundamental bi-differential}. It is defined as the
unique symmetric quadratic meromorphic differential on
$\overline{X}_u\times \overline{X}_u$ that has a pole of order 2 along
the diagonal with no-residues normalized by
\ben
B(p',p'') = \frac{dt(p')dt(p'')}{(t(p')-t(p''))^2} + \cdots,
\een 
and 
\ben
\oint_{p'\in \alpha_i} B(p',p'')=0,\quad 1\leq i\leq g,
\een
where in the first condition $p'$ and $p''$ are sufficiently close to
some point $p_0\in \overline{X}_u$, $t$ is a local coordinate in a
neighborhood of $p_0$, and the dots stand for terms that are
holomorphic in a neighborhood of $(p_0,p_0)$ in $\overline{X}_u\times
\overline{X}_u$. Note however, that the cohomology class $[\omega_i]$
is independent of the choice of a holomorphic branch of $\omega_i$
in a neighborhood of $C$. 

\begin{proposition}\label{prop:gb}
The cohomology classes $\{[\omega_i]\}_{i=1}^N\subset \widehat{\H}(B)$
form a good basis in which the higher residue pairing takes the form
\ben
K([\omega_i],[\omega_j])=z\delta_{ij},\quad 1\leq i,j\leq N.
\een 
\end{proposition}
\proof
Recalling Proposition \ref{per-iso}, Part a), we get that condition
(GB1) is equivalent to saying that for every open connected Stein
subset $U\subset V$ the projection of $\omega_i$ to 
\beq\label{H-mod-z}
\widehat{\H}(U)/z \widehat{\H}(U)=
\bpi_*\Omega^{1,\infty}_{\bX/B}(U)/
d_{\bX/B}F\wedge \bpi_*\Omega^{0,\infty}_{\bX/B}(U)
\eeq
is an $\O_B(U)$-basis. This however is obvious, because
$\omega_i(p_j)=-\delta_{ij}dt_j$, so the condition (G1) in the
definition of a good basis holds.

Using formula \eqref{hrp} for the higher residue pairing, we get 
\ben
K([\omega_a],[\omega_b])=\frac{1}{2\pi\sqrt{-1}} \, \sum_{i=1}^N
\int_{\Gamma_i} e^{F(p')/z} \omega_a(p')\, 
\int_{\Gamma^\vee_i} e^{-F(p')/z} \omega_b(p'),
\een
where on the RHS the oscillatory integrals should be identified with
the corresponding stationary phase asymptotic as $z\to 0$. In
particular, only the germs of the integration paths $\Gamma_i$ and
$\Gamma_i^\vee$ near the critical point $p_i$ are relevant. This
justifies why we can replace the holomorphic form in
$\Omega^{1,\infty}_{\bX/B}(\bX)[\![z]\!]$ that represent the cohomology class
$[\omega_a]$ by a holomorphic branch of $\omega_a$ defined
only in a neighborhood of the critical points. 

We are going to proof that condition $(G2)$ from the definition of a
good basis is equivalent to a certain identity for the matrix series
$R_\Sigma(z)=1+R_{\Sigma,1} z+\cdots$ whose $(a,i)$-th entry is defined by
\ben
[R_\Sigma(z)]_i^a:=-(-2\pi z)^{-1/2} 
\int_{p'\in \Gamma_i} e^{(F(p')-u_i)/z} \omega_a(p').
\een
\begin{remark}
In the notation of \cite{DNOPS1}, $R_\Sigma(z) = \widehat{R}^{-1}(-z)$.
\end{remark}
The leading order term of $R_\Sigma(z)$ is the identity matrix as
it can be seen easily from \eqref{spa} and \eqref{local-pi}. 
The higher residue pairing takes the form
\ben
K([\omega_a],[\omega_b]) = z 
\sum_{i=1}^N [R_\Sigma(z)]^a_i\, [R_\Sigma(-z)]^b_i.
\een
We will prove that the series $R_\Sigma(z)$
satisfies the symplectic condition $R_\Sigma(z)
R_\Sigma^T(-z)=1$, therefore
$K([\omega_a],[\omega_b])=z\delta_{ab}$. 

The symplectic condition is stated in \cite{DNOPS1}, Lemma 5.1 and it
is a consequence of a more general identity proved in the
lemma. For the sake of completeness and for the reader's convenience,
let us prove directly the symplectic condition. Our argument follows
the ideas of \cite{DNOPS1}.  It is more convenient to prove that 
\beq\label{symp-id}
\sum_{a=1}^N [R_\Sigma(z)]^a_i\, [R_\Sigma(-z)]^a_j = \delta_{ij}.
\eeq
By definition the LHS of \eqref{symp-id} is 
\ben
\sum_{a=1}^N 
\frac{z^{-1}}{2\pi\sqrt{-1}} \int_{\Gamma_i\times \Gamma^\vee_j}
e^{(F(p')-u_i)/z}e^{-(F(p'')-u_j)/z} 
\omega_a(p')\omega_a(p'').
\een
The key observation is that 
\ben
\omega_a(p')\omega_a(p'') =\Big( \operatorname{res}_{q=p_a}
\frac{\omega_{p'}(q)\omega_{p''}(q)}{d_{X/B} F(q)}\Big)\, d_{X/B}F(p')d_{X/B}F(p'').
\een
Using the residue theorem on $\overline{X}_u$ we get 
\ben
\sum_{a=1}^N \operatorname{res}_{q=p_a}
\frac{\omega_{p'}(q)\omega_{p''}(q)}{d_{X/B} F(q)} = 
-(\operatorname{res}_{q=p'}+\operatorname{res}_{q=p''})
\frac{\omega_{p'}(q)\omega_{p''}(q)}{d_{X/B} F(q)} = 
-\frac{\omega_{p''}(p')}{d_{X/B} F(p')} - \frac{\omega_{p'}(p'')}{d_{X/B} F(p'')},
\een
where we used the Cauchy theorem
\ben
\operatorname{res}_{q=p} \omega_p(q) f(q) = f(p)
\een
for any function $f$ on $\overline{X}_u$ holomorphic in a
neighborhood of $p$. The LHS of \eqref{symp-id} turns into 
\ben
-\frac{z^{-1}}{2\pi\sqrt{-1}} \int_{\Gamma_i\times \Gamma^\vee_j}
e^{(F(p')-u_i)/z}e^{-(F(p''-u_j)/z} 
\Big(\omega_{p''}(p')\, d_{X/B}F(p'') + \omega_{p'}(p'')\, d_{X/B}F(p')\Big).
\een
Let us consider first the case when $i\neq j$. In this case
$\omega_{p'}(p'')$ and $\omega_{p''}(p')$ have no singularities, so
using integration by parts we get
\ben
-\frac{z^{-1}}{2\pi\sqrt{-1}} \int_{\Gamma_i\times \Gamma^\vee_j}
e^{(F(p')-u_i)/z}e^{-(F(p''-u_j)/z} 
\omega_{p''}(p')\, d_{X/B}F(p'') = \\
-\frac{1}{2\pi\sqrt{-1}} \int_{\Gamma_i\times \Gamma^\vee_j}
e^{(F(p')-u_i)/z}e^{-(F(p''-u_j)/z} 
B(p'',p')
\een
and 
\ben
-\frac{z^{-1}}{2\pi\sqrt{-1}} \int_{\Gamma_i\times \Gamma^\vee_j}
e^{(F(p')-u_i)/z}e^{-(F(p''-u_j)/z} 
\omega_{p'}(p'')\, d_{X/B}F(p') = \\
\frac{1}{2\pi\sqrt{-1}} \int_{\Gamma_i\times \Gamma^\vee_j}
e^{(F(p')-u_i)/z}e^{-(F(p''-u_j)/z} 
B(p',p'').
\een
The fundamental bi-differential is symmetric, so the above integrals cancel out,
i.e., formula \eqref{symp-id} holds for $i\neq j$.

If $i=j$, then as discussed above, since we are interested only in the
asymptotic of the integrals, we may assume that $\Gamma_i$ is a
small path defined in a neighborhood of $p_i$. Let us split
$\omega_{p'}(p'')$ and $\omega_{p''}(p')$ into singular and regular
parts:
\ben
\omega_{p'}(p'') = \frac{dt_i(p'')}{t_i(p'')-t_i(p')} + \omega^{\rm reg}_{p'}(p'')
\een
and 
\ben
\omega_{p''}(p') = \frac{dt_i(p')}{t_i(p')-t_i(p'')} + \omega^{\rm reg}_{p''}(p').
\een
The regular parts do not contribute, because they cancel out after
integration by parts just like in the case of $i\neq j$. 
While the singular parts add up to 
\ben
\frac{z^{-1}}{2\pi\sqrt{-1}} \int_{\Gamma_i\times \Gamma^\vee_i}
e^{(F(p')-u_i)/z}e^{-(F(p''-u_i)/z} 
dt_i(p')dt_i(p'').
\een
Again, since we are interested only in the asymptotic as $z\to 0$, we
may assume that 
$
F(p)=u_i+\frac{1}{2}t_i^2
$
and 
\ben
\Gamma_i = \{t_i=\sqrt{-1} s\ |\ -\infty<s<+\infty\},\quad
\Gamma^\vee_i = \{ t_i=s\ |\ -\infty<s<+\infty\}.
\een
The above integral turns into
\ben
\frac{z^{-1}}{2\pi} \Big(\int_{-\infty}^\infty e^{-s^2/(2z)} ds \Big)^2 =
\frac{1}{2\pi} \Big(\int_{-\infty}^\infty e^{-s^2/2} ds \Big)^2 = 1.\qed
\een

\subsection{The Gauss--Manin connection}
Let us derive a formula for the Gauss--Manin connection in the frame
of $\widehat{\H}$ given by the good basis $\{[\omega_i]\}_{i=1}^N$
defined in the previous section. The formulas can be obtained easily
from \cite{DNOPS2}, Theorem 7. However, for the sake of completeness,
let us derive them directly.

The main tool is the Rauch's variation formula (see \cite{Ra})
\beq\label{Rauch:vf}
\frac{\delta \omega_p}{\delta u_i} (q) = \operatorname{res}_{q'=p_i}
\frac{\omega_p(q')B(q',q)}{d_{X/B}F(q')} .
\eeq
Let us explain the meaning of the LHS. For fixed $u\in B$, let us
restrict our family $X\to B$ to a small disc $B_i^*(u)\subset B$ with
center $u$, obtained by varying only the $i$-th branch point $u_i$,
i.e., 
\ben
B_i^*(u)=\{u^*\in B \  |\ u^*_j=u_j \mbox{ for } i\neq j ,\
|u^*_i-u_i|\ll 1\}.
\een
The resulting family $X^*\to B_i^*(u)$ is a deformation of $X_u$
constructed by replacing the holomorphic coordinate $t_i$ on $X_u$ around the
critical point $p_i$ in which $F(p)=u_i+\frac{1}{2} t_i(p)^2$ with a
holomorphic coordinate $t_i^*$, s.t., $F(p)=u^*_i + \frac{1}{2}
t^*(p)^2$. In particular, by removing a tubular neighborhood of
$C_i\cap X^*$ in $X^*$, we get a holomorphically trivial
family. Therefore, we can identify the differentials $\omega_p(q)$,
$q\in X^*$ as a family of differentials on $\overline{X}_u$
holomorphic in the complement of a small disc around $p_i$. The LHS of
\eqref{Rauch:vf} is interpreted as the usual derivative with respect
to $u^*_i$ evaluated at $u^*_i=u_i$.  

The Rauch's variational derivative is compatible with the Gauss--Manin
connection 
\ben
\partial_{u_i}\,  \int_\Gamma e^{F(p)/z} \omega_a(p)=
\int_\Gamma e^{F(p)/z} \frac{\delta \omega_a}{\delta u_i} (p),
\een
where the integrals are interpreted via their asymptotic as $z\to 0$,
so we may assume that the integration cycle $\Gamma$ is supported in a tubular
neighborhood of the critical points, where $\omega_a(p)$ is holomorphic.
Recalling the definition of $\omega_a(p)$ we get 
\ben
\frac{\delta \omega_a}{\delta u_i} (p) = 
\operatorname{res}_{q'=p_a} \Big(
\delta_{ia} \frac{\omega_p(q')}{t_a(q')^3}+
\operatorname{res}_{q''=p_i}
\frac{\omega_p(q'')B(q'',q')}{t_a(q')d_{X/B}F(q'')} \Big)\, d_{X/B}F(p)
\een
where we used the Rauch's variational formula \eqref{Rauch:vf} and 
\ben
\frac{\delta (t_a^{-1})}{\delta u_i} = \delta_{ia}t_a^{-3}.
\een
Let us consider first the case when $i\neq a$, then the above formula
turns into 
\ben
\frac{\delta \omega_a}{\delta u_i} (p) =\beta_{ai}\, \omega_i(p) , 
\een
where
\beq\label{beta_ai}
\beta_{ai} =  
\operatorname{res}_{q'=p_a}
\operatorname{res}_{q''=p_i}
\frac{B(q',q'')}{t_a(q')t_a(q'')}.
\eeq
If $i=a$, then we use the residue theorem to compute the residue with
respect to $q''$. Note that the poles are only at $q''=p_a$, $1\leq
a\leq N$, $q''=p$ and $q''=q'$, so the residue operation with respect
to $q''$ can be replaced by 
\ben
-\sum_{a\neq i}\operatorname{res}_{q''=p_a} -
\operatorname{res}_{q''=p} -\operatorname{res}_{q''=q'}.
\een
The contribution of the sum of the residues is exactly what we have
computed in the case $i\neq a$, i.e., 
\ben
-\sum_{a\neq i} \beta_{ia}\, \omega_a(p).
\een
The contribution of the residue at $q''=p$ is computed via the Cauchy theorem
\ben
-\operatorname{res}_{q'=p_i}\frac{B(p,q')}{t_i(q')} = 
-d_{X/B} \Big(\frac{\omega_i(p)}{d_{X/B}F}
\Big).
\een
The residue at $q''=q'$ can be computed using that
\ben
\operatorname{res}_{q''=q'} f(q'')B(q'',q') =df(q') 
\een
for every meromorphic function $f$. Using the above formula we get
\ben
-\operatorname{res}_{q''=q'}
\frac{\omega_p(q'')B(q'',q')}{t_i(q')d_{X/B}F(q'')} =
-t_i(q')^{-1} d_{q'} \Big(\frac{\omega_p(q')}{d_{X/B}F(q')}\Big).
\een
Since we are interested in the residue at $q'=p_i$, let us use the
local coordinate $t_i(q')$. Then $\omega_p(q') = f_i(p,q')dt_i(q')$
for some function $f_i$ holomorphic in $q'$ near $p_i$  
and $d_{X/B}F(q')=t_i(q')dt_i(q')$. Therefore
\ben
-t_i(q')^{-1} d_{q'} \Big(\frac{\omega_p(q')}{d_{X/B}F(q')}\Big) = 
-t_i(q')^{-3} \omega_p(q') -d_{q'}f_i(p,q')
\een
and we get the following formula
\ben
\sum_{a=1}^N\frac{\delta \omega_a}{\delta u_i} = -d_{X/B} \Big(\frac{\omega_i(p)}{d_{X/B}F}\Big).
\een
For the oscillatory integrals we get the following equations
\ben
\partial_{u_i}\int_\Gamma e^{F(p)/z} \omega_a & = & \beta_{ai}
\int_\Gamma e^{F(p)/z} \omega_i,\quad i\neq a,\\
\sum_{a=1}^N \partial_{u_a} \int_\Gamma e^{F(p)/z} \omega_i & = & 
z^{-1} \int_\Gamma e^{F(p)/z} \omega_i.
\een

Similarly, we can compute 
\ben
E\, \int_\Gamma e^{F(p)/z}\omega_i(p),
\een
where $E=\sum_{a=1}^N u_a\partial_{u_a}$. 
Using the Rauch's variation formula we get
\ben
\int_\Gamma e^{F(p)/z}\sum_{a=1}^N
u_a \operatorname{res}_{q'=p_i}\Big(
\delta_{ia}\frac{\omega_p(q')}{t_i(q')^3} + 
\operatorname{res}_{q''=p_a} 
\frac{ \omega_p(q'')B(q'',q') }{ t_i(q')d_{X/B}F(q'') }
\Big) d_{X/B}F(p).
\een
Exchanging the order of the residues we get
\beq\label{Euler-eq1}
\sum_{a=1}^N
u_a \operatorname{res}_{q''=p_a} 
\operatorname{res}_{q'=p_i} 
\frac{ \omega_p(q'')B(q'',q') }{ t_i(q')d_{X/B}F(q'') }
\ d_{X/B}F(p).
\eeq
Note that for $a\neq i$ the above form has a pole of order 1 at
$q''=p_a$. Therefore we can rewrite the sum of the terms for which
$a\neq i$ as
\ben
\sum_{a\neq i} 
\operatorname{res}_{q''=p_a} 
\operatorname{res}_{q'=p_i} 
\frac{ F(q'') \omega_p(q'')B(q'',q') }{ t_i(q')d_{X/B}F(q'') }
\ d_{X/B}F(p).
\een
Using the residue theorem we  can replace the residue operation
$\sum_{a\neq i}\operatorname{res}_{q''=p_a}$ with
$-\operatorname{res}_{q''=p_i}-\operatorname{res}_{q''=p}$. 
Note that 
\ben
\operatorname{res}_{q'=p_i} 
\frac{ B(q'',q') }{ t_i(q')} = d_{q''}
\Big(\frac{\omega_i(q'')}{d_{X/B} F(q'')}\Big)
.
\een
Therefore, the sum \eqref{Euler-eq1} turns into the sum of 
\beq\label{Euler-eq2}
-\operatorname{res}_{q''=p_i} \frac{(F(q'')-u_i)
  \omega_p(q'')}{d_{X/B} F(q'')}\,
d_{q''}
\Big(\frac{\omega_i(q'')}{d_{X/B} F(q'')}\Big) \ d_{X/B}F(p)
\eeq
and 
\beq\label{Euler-eq3}
-F(p) d_{p}
\Big(\frac{\omega_i(p)}{d_{X/B} F(p)}
\Big).
\eeq
Using that $F(q'')=u_i+\frac{1}{2}t_i(q'')^2$ and integration by parts, we get that
\eqref{Euler-eq2} is 
\ben
\frac{1}{2} \, \operatorname{res}_{q''=p_i}
\frac{\omega_i(q'')}{d_{X/B}F(q'')} \, \omega_p(q'')\, d_{X/B}F(p) = 
- \frac{1}{2} \, \omega_i(p),
\een
where we used that the 1-form $\frac{\omega_i(q'')}{d_{X/B}F(q'')}\,
dt_i(q'')$ has a pole of order 1 and residue $-1$ at $q''=p_i$ (see
formula \eqref{local-pi}). Note also that the contribution of
\eqref{Euler-eq3} to the oscillatory integral is
\ben
-\int_\Gamma e^{F(p)/z}F(p) d_{p}
\Big(\frac{\omega_i(p)}{d_{X/B} F(p)}
\Big) = \int_\Gamma e^{F(p)/z}\Big(\frac{F(p)}{z} +1\Big)\omega_i
\een
We get the following equation 
\ben
(z\partial_z+E)\, \int_\Gamma e^{F(p)/z} \omega_i(p) = 
\frac{1}{2}\, \int_\Gamma e^{F(p)/z} \omega_i(p).
\een
In other words we proved the following Proposition.
\begin{proposition}\label{GM:gb-frame}
In the good basis frame $\{[\omega_i]\}_{i=1}^N$ the Gauss--Manin
connection takes the form
\ben
\nabla_{\partial_{u_i}} [\omega_a ]& = & \beta_{ai}
[\omega_i],\quad i\neq a,\\
\sum_{a=1}^N z\nabla_{\partial_{u_a}} [ \omega_i] & = & 
[ \omega_i], \\
(z\nabla_{\partial_z} + \nabla_E)
[\omega_i] & = &\frac{1}{2}\,[\omega_i],\quad
1\leq i\leq N.
\een
\end{proposition}

\section{Primitive forms}\label{sec:pf-classif}

The goal in this section is to define a primitive form and prove that
their classification reduces to the differential equations that
classify semi-simple Frobenius manifolds.

\subsection{Definition}\label{subsec:pf}

Let $\omega\in \widehat{\H}(V) $ be a cohomology class whose leading
order term is a volume form along $C$. Recall that under this
condition the period map \eqref{V:per-iso} is an isomorphism.
Furthermore, the form $\omega$ is
called a {\em primitive form} on $V$ if the following 5 properties are
satisfied
\begin{enumerate}
\item[(PF1)]
$K^{(n)}(z\nabla_{\partial_{u_i}} \omega ,z\nabla_{\partial_{u_j}} \omega )=0$
for all $1\leq i,j\leq N$ and $n\geq 1$.
\item[(PF2)]
$K^{(n)}(z\nabla_{\partial_{u_i}} \omega 
z\nabla_{\partial_{u_j}} \omega , z\nabla_{\partial_{u_k}} \omega )=0$
for all $1\leq i,j,k\leq N$ and $n\geq 2$.
\item[(PF3)]
$K^{(n)}(z^2\nabla_{\partial_{z}} \omega 
z\nabla_{\partial_{u_i}} \omega , z\nabla_{\partial_{u_j}} \omega )=0$
for all $1\leq i,j\leq N$ and $n\geq 2$.
\item[(PF4)]
The cohomology class $ \omega $ is homogeneous in the following sense 
\ben
(z\nabla_{\partial_z} + \nabla_E) \omega =r\,  \omega ,
\een
where $E=\sum_{i=1}^N u_i\partial_{u_i}$ and $r$ is some constant
independent of $z$ and $u$.
\item[(PF5)]
\ben
\sum_{i=1}^N z\nabla_{\partial_{u_i}} \omega = \omega .
\een
\end{enumerate}
If $\omega$ is a primitive form, then the Gauss--Manin connection
\eqref{GM-ui} and \eqref{GM-z} takes a very simple form. Note that
the residue pairing
\ben
(\partial_{u_i},\partial_{u_j}):=
K^{(0)}(z\nabla_{\partial_{u_i}} \omega ,z\nabla_{\partial_{u_j}} \omega )=\frac{\delta_{ij}}{\Delta_i},
\een
where $\Delta_i$, $1\leq i\leq N$ are some holomorphic functions,
s.t., $\Delta_i(u)\neq 0$ for all $u\in U$. From formula \eqref{GM-ui}
and axiom (PF2) we get 
\ben
C_{ij;n}^k=0,\quad n\geq 2.
\een
The leading order terms $C_{ij;0}^k=\delta_{ij}\delta_{ik}$ coincide
with the structure constants of the Frobenius multiplication, while
the Leibnitz rule for the higher residue pairing shows that
$C_{ij;1}^k$ are the Christophel's symbols for the residue pairing.
Similarly, from \eqref{GM-z} and axiom (PF3) we get that 
\ben
C_{i;n}^j=0,\quad n\geq 2.
\een
The leading order term $C_{i;0}^j=-u_i\delta_{ij}$ is the matrix of the
linear operator of Frobenius multiplication by 
$-E$. Axiom (PF4) yields the following formula 
\ben
\nabla^{\rm L.C.}_{\partial_{u_i}} E = (r+1)\partial_{u_i}-
\sum_{j=1}^N C_{i;1}^j \partial_{u_j},
\een
where $\nabla^{\rm L.C.}$ is the Levi--Civita connection for the
residue pairing. It is straightforward to prove that the gauge
transformation $z^{1/2}\nabla z^{-1/2}$ of the Gauss--Manin
connection,  turns into the deformed flat connection of a semi-simple
Frobenius structure on $U$ of conformal dimension $D=1-2r$ (see \cite{SaT} for
more details).

\subsection{Preliminary notation}\label{sec:pn}

Let $\{[\omega_i]\}_{i=1}^N\subset \widehat{\mathcal{H}}(B)$ be the
good basis constructed in the previous section. Let us also fix an
open subset $V\subset B$. We would like to find
all formal series 
$$
c_i(u,z)\in \CC[\![z]\!], \quad 1\leq i\leq N,
$$
depending analytically on $u\in V$, s.t., 
\ben
\sum_{i=1}^N c_i(u,z) [\omega_i]
\een
is a primitive form in $\widehat{\mathcal{H}}(V)$. 
Put 
\ben
\omega:=([\omega_1],\dots,[\omega_N]).
\een
Then the Gauss--Manin connection (see Proposition \ref{GM:gb-frame})
takes the form
\ben
\nabla_i\omega & = & \omega\, \widetilde{B}_i(u,z),\quad 1\leq i\leq N,\\
\nabla_{z\partial_z+E}\omega  & = & \frac{1}{2}\omega,
\een 
where $\nabla_i:=\nabla_{\partial_{u_i}}$ and the matrix 
\ben
\widetilde{B}_i(u,z) = z^{-1} E_{ii}+ B_i(u),\quad B_i(u):=\sum_{j:j\neq i}
\beta_{ij}(u) (E_{ij}-E_{ji}),
\een
where $E_{ij}$ denotes the square matrix of size $N$ with all entries
0, except for the entry in position $(i,j)$, which is 1.

Let us denote by $c(u,z):=(c_1(u,z),\dots,c_N(u,z))^T$ the column
vector whose entries are the functions that we would like to
classify. Put 
\ben
\widetilde{\omega}=(\widetilde{\omega}_1,\dots,\widetilde{\omega}_N),
\een
where
\ben
\widetilde{\omega}_i=z\nabla_i (\omega \, c) = 
\omega\Big( E_{ii}c + z(B_i c+\partial_{u_i} c)\Big),\quad 1\leq i\leq N.
\een
Note that $\widetilde{\omega}=\omega \widetilde{R}(u,z)$, where
$\widetilde{R}(u,z)$ is a matrix whose $i$-th column is given by 
\ben
E_{ii}c(u,z) + z(B_i (u)c(u,z)+\partial_{u_i} c(u,z)).
\een
If $\omega\, c(u,z)$ is a primitive form for $u\in V$, then since by
construction $\omega_i(p_j)=-\delta_{ij} dt_j$, we get that 
\ben
\omega(p_j)\, c(u,z) = (-c_j(u,0) +O(z))dt_j.
\een
By definition (see Section \ref{sec:pf}) $c_j(u,0)\neq 0$ for all
$u\in U$.  Put 
\ben
C(u):=\operatorname{Diag}(c_1(u,0),\dots,c_N(u,0))
\een
and $R(u,z):=\widetilde{R}(u,z) C(u,z)^{-1}$. Note that the series 
\ben
R(u,z)=1+R_1(u)z+R_2(u) z^2+\cdots,
\een
where $R_i(u)$ are square matrices of size $N$ whose entries depend
analytically on $u\in V$.

\subsection{Differential and algebraic constraints}\label{sec:da-constr}
In this section, we will be using quite frequently the following
notation. If $A$ is a matrix, then $A_{ij}$ will be the entry in raw
$i$ and column $j$. In case $A$ is a raw (resp. column), then we
denote by $A_i$ the $i$-th entry of the raw (resp. column).

Axiom (PF1) is equivalent to
\ben
K(\widetilde{\omega}_i, \widetilde{\omega}_j) = z
K^{(0)}(\widetilde{\omega}_i,\widetilde{\omega}_j),\quad 
1\leq i,j\leq N.
\een
Using that 
\ben
\widetilde{\omega}_i = (\omega R(u,z)C(u))_i = \sum_{k=1}^N \omega_k R_{ki}(u,z)c_i(u,0)
\een
and $K(\omega_k,\omega_\ell) = z\delta_{k\ell}$, we get 
\ben
K^{(0)}(\widetilde{\omega}_i, \widetilde{\omega}_j) = c_i(u,0)c_j(u,0)\delta_{ij}
\een
and 
\beq\label{eq:pf1}
R(u,z)^T R(u,-z)=1.
\eeq
Put
\ben
\gamma_{ij}:= \partial_{u_j}c_i(u,0)\, c_j(u,0)^{-1},\quad 1\leq
i,j\leq N
\een
and 
\ben
\Gamma_i= \sum_{j:j\neq i} \gamma_{ij}(E_{ij}-E_{ji}), \quad 1\leq
i\leq N.
\een
\begin{lemma}\label{le:pf2}
The operator series $R(u,z)$ satisfies the following differential
equations
\ben
z\partial_{u_i}R = [R,E_{ii}]+z(R\Gamma_i-B_i R),\quad 1\leq i\leq N.
\een
\end{lemma}
\proof
According to Axiom (PF2), 
\beq\label{eq:pf2-1}
K(z\nabla_{\partial_{u_i}}
  \widetilde{\omega}_j,\widetilde{\omega}_k)\quad \in \quad
\CC\, z +\CC\, z^2.
\eeq
On the other hand
\ben
z\nabla_{\partial_{u_i}} \widetilde{\omega}_j = 
(z\nabla_{\partial_{u_i}} (\omega R C))_j=
((z\nabla_{\partial_{u_i}} \omega) R C + \omega z\nabla_{\partial_{u_i}}(RC) )_j.
\een
Recalling the differential equations for $\omega$ we get 
\ben
z\nabla_{\partial_{u_i}} \widetilde{\omega}_j = \sum_{\ell=1}^N
[\omega_\ell]\, (E_{ii}\widetilde{R} + z(B_i \widetilde{R}
+\partial_{u_i} \widetilde{R}))_{\ell j}. 
\een
Therefore \eqref{eq:pf2-1} is equivalent to
\ben
\sum_{\ell=1}^N 
(E_{ii}\widetilde{R} + 
z(B_i \widetilde{R}+\partial_{u_i} \widetilde{R}) )_{\ell j}
\widetilde{R}_{\ell k} (u,-z) \quad \in \quad
\CC+ \CC\, z.
\een
The above condition, written in matrix form, becomes
\ben
\widetilde{R}(u,-z)^T ((E_{ii}\widetilde{R} + 
z(B_i \widetilde{R}+\partial_{u_i} \widetilde{R}) ) =
A_0^{(i)}+A_1^{(i)} z,
\een
where $A_\alpha^{(i)}$, $\alpha=0,1$, are some matrices independent of
$z$. Recalling $\widetilde{R}(u,z)=R(u,z)C(u)$ we get 
\beq\label{eq:pf2-2}
E_{ii}R(u,z)+z(B_i(u)R+\partial_{u_i} R(u,z)) = R(u,z) (B_0^{(i)}+
B_1^{(i)}\, z),
\eeq
where
\ben
B_0^{(i)} = C(u)^{-1} A_0^{(i)} C(u)^{-1},\quad 
B_1^{(i)} = C(u)^{-1} A_1^{(i)} C(u)^{-1} -\partial_{u_i} C(u) C(u)^{-1}.
\een
Comparing the coefficients in front of $z^0$ and $z^1$ in
\eqref{eq:pf2-2} we get 
\ben
B_0^{(i)} & = &  E_{ii},\\
B_1^{(i)} & = & [E_{ii},R_1]+ B_i.
\een
The commutator $[E_{ii},R_1]$ can be expressed in terms of $B_i$ and
$C(u)$. By definition $R_1=\widetilde{R}_1 C(u)^{-1}$, where
$\widetilde{R}_1$ is the coefficient in front of $z^1$ in
$\widetilde{R}(u,z)$. Recalling the definition of
$\widetilde{R}(u,z)$, we get
\ben
(\widetilde{R}_1)_{ij} =  (B_j(u)c(u,0))_i + \partial_{u_j}
c_i(u,0) = 
\partial_{u_j} c_i(u,0) -\beta_{ij}(u) c_j(u,0),\quad \mbox{for } i\neq j. 
\een
Therefore
\ben
(R_1)_{ij} = \gamma_{ij}-\beta_{ij}
\een
and for the commutator we get
\ben
[E_{ii},R_1] = \Gamma_i -B_i.
\een
The above formula yields $B_1^{(i)}=\Gamma_i$, so the differential
equation \eqref{eq:pf2-2} turns into the differential equation we
wanted to prove. 
\qed

\begin{lemma}\label{le:pf3-5}
The following differential equations hold
\ben
(z\partial_z + E)R(u,z) & = & 0, \\
(z\partial_z + E)c(u,z) & = & \Big(r-\frac{1}{2}\Big) c(u,z), \\
(\partial_{u_1}+\cdots +\partial_{u_N})\, c(u,z) & = & 0.
\een
\end{lemma}
\proof
The differential equations are consequences of Axioms
(PF2)--(PF5). Let us derive the first one. The computations for the
remaining ones are straightforward. 

According to Axiom (PF2) and (PF3), we have
\beq\label{eq:pf2-3}
K(\nabla_{z\partial_z +E}\, \widetilde{\omega}_i,\widetilde{\omega}_j)
\quad \in \quad 
\CC + \CC\, z.
\eeq
On the other hand
\ben
\nabla_{z\partial_z +E}\, \widetilde{\omega}_i = 
(\nabla_{z\partial_z+E}\, \omega \widetilde{R})_i =
\Big( \omega\Big( \frac{1}{2} \widetilde{R}+(z\partial_z+E) \widetilde{R}\Big)\Big)_i.
\een
Using that $K(\omega_k\omega_\ell)=z\delta_{k\ell}$,
$\widetilde{R}(u,z)=R(u,z)C(u)$, and \eqref{eq:pf1}, we get that
\eqref{eq:pf2-3} is equivalent to 
\ben
R(u,-z)^T \, \Big( z\partial_z + E +\frac{1}{2}\Big)
R(u,z)
\quad \in \quad 
\operatorname{Mat}_{N\times N}(\CC)\, z^{-1} + \operatorname{Mat}_{N\times
  N}(\CC) .
\een
The above condition implies 
\ben
(z\partial_z + E )R(u,z) = R(u,z) \, (A_0 z^{-1} +A_1).
\een
Comparing the coefficients in front of $z^{-1}$ and $z^0$ we get that
$A_0=A_1=0$, which is exactly what we need. 
\qed 

Let $c(u,z)\in \CC^N[\![z]\!]$ be a formal series depending analytically
on $u\in B$, s.t., the $i$-th component $c_i(u,z)$ $(1\leq i\leq N)$ of $c(u,z)$
satisfies $c_i(u,0)\neq 0$ for some $u\in B$. Let us define $V\subset B$ to be
the open subset of those $u\in B$, such that, $c_i(u,0)\neq 0$ for all
$1\leq i\leq N$. According to our assumptions $V\neq \emptyset$.
\begin{proposition}\label{prop:pf}
The cohomology class 
\ben
\omega\, c(u,z) = \sum_{i=1}^N c_i(u,z)\, [\omega_i]
\een
is a primitive form in $\widehat{\mathcal{H}}(V)$ if and only if
the following equations are satisfied  
\beqa
\label{eq:pf-1}
&&
(E_{ii} + z (B_i(u) +\partial_{u_i})) c(u,z) = R(u,z) C(u) e_i ,\\
\label{eq:pf-2}
&&
R(u,z)R(u,-z)^T=1 ,\\
\label{eq:pf-3}
&&
z\partial_{u_i} R(u,z) = [R(u,z),E_{ii}]+ z(R(u,z)\Gamma_i(u)-B_i(u)
R(u,z)) ,\\
\label{eq:pf-4}
&&
(z\partial_z + E )R(u,z) =0, \\
\label{eq:pf-5}
&&
(z\partial_z + E )c(u,z) = (r-1/2) c(u,z) ,\\
\label{eq:pf-6}
&&
(\partial_{u_1}+\cdots +\partial_{u_N}) c(u,z)=0,
\eeqa
where $1\leq i\leq N$ and $e_i$ is the vector column, whose $i$-th entry is 1 and all
other entries are 0.
\end{proposition}

In one direction, the Proposition is already established. Expecting
more carefully the derivation of equations \eqref{eq:pf-1}--\eqref{eq:pf-6}, it is
straightforward to check that the equations in Proposition
\ref{prop:pf} are sufficient to guarantee that $\omega \, c(u,z)$
satisfies all axioms (PF1)--(PF5). 

\subsection{Solving the equations for primitive forms}

We are going to prove that the solutions $c(u,z)$ to
\eqref{eq:pf-1}--\eqref{eq:pf-6} are uniquely determined from $c(u,0)$
and we are going to derive differential equations for $c(u,0)$, which
guarantee that the reconstruction of $c(u,z)$ from $c(u,0)$ is a
solution to \eqref{eq:pf-1}--\eqref{eq:pf-6}. 

\begin{lemma}\label{le:rec-R}
Let $c(u,z)$ be a solution to \eqref{eq:pf-1}--\eqref{eq:pf-6}. Then 

a) We have $\gamma_{ij}=\gamma_{ji}.$

b) The series $R(u,z)$ is uniquely determined from $\gamma_{ij}$,
$1\leq i,j\leq N$. 
\end{lemma}
\proof
a) The symplectic condition $R(u,-z)^TR(u,z)=1$ implies that
$R_1^T=R_1$, so $[E_{ii},R_1]^T=-[E_{ii},R_1]$. On the other hand, 
we already proved (see Lemma \ref{le:pf2}) that 
\ben
[R_1,E_{ii}]+\Gamma_i-B_i=0.
\een
Therefore $\Gamma_i^T=-\Gamma_i$, which is equivalent to
$\gamma_{ij}=\gamma_{ji}$. 

b)
Comparing the coefficients in front of $z^{k+1}$ in \eqref{eq:pf-2}
and \eqref{eq:pf-3} we get 
\beq\label{rec:R-1}
\partial_{u_i} R_k = [R_{k+1},E_{ii}] + R_k\Gamma_i- B_i R_k
\eeq
and
\ben
\sum_{i=1}^N u_i\partial_{u_i} R_{k+1}  = -(k+1)R_{k+1}.
\een
These recursions can be solved uniquely for $R_k$, $k\geq 1$, in terms
of $\gamma_{ij}$. Indeed, we already proved that 
\ben
(R_1)_{ij} = \gamma_{ij}-\beta_{ij},\quad i\neq j.
\een 
To determine the diagonal entries of $R_1$, we combine the above
equations to get
\ben
-R_1= \sum_i u_i [R_2, E_{ii}] + R_1\Big(\sum_i
u_i\Gamma_i\Big)-\Big(\sum_i u_iB_i\Big) R_1.
\een
From here we get 
\ben
(R_1)_{ii} = \sum_{j:j\neq i} (u_i-u_j)(\gamma_{ij}-\beta_{ij})(\gamma_{ij}+\beta_{ij}).
\een
In general, we use equations \eqref{rec:R-1} to determine the
non-diagonal entries of $R_{k+1}$, while for the diagonal entries we
use 
\ben
-(k+1)R_{k+1}= \sum_i u_i [R_{k+2}, E_{ii}] + R_{k+1}\Big(\sum_i
u_i\Gamma_i\Big)-\Big(\sum_i u_iB_i\Big) R_{k+1}.
\qed
\een

\begin{lemma}\label{le:rec-c}
If $c(u,z)$ is a solution to \eqref{eq:pf-1}--\eqref{eq:pf-6}, then
\ben
c(u,z)=R(u,z)\, \sum_{i=1}^N c_i(u,0) e_i.
\een 
\end{lemma}
\proof
Equation \eqref{eq:pf-3} can be written us
\ben
(z(\partial_{u_i}+B_i(u))+E_{ii})\circ R(u,z) = R(u,z)\circ (z(\partial_{u_i}+\Gamma_i(u))+E_{ii}),
\een
i.e., $R(u,z)$ is a gauge transformation intertwining two
connections. 
Using the above relation, we get that $C(u)e_i$ is 
\ben
R(u,z)^{-1}   (z(\partial_{u_i}+B_i(u))+E_{ii}) c(u,z)= 
(z(\partial_{u_i}+\Gamma_i(u))+E_{ii})(R(u,z)^{-1}c(u,z)).
\een
Summing over all $i$, we get
\ben
c(u,0) = \sum_{i=1}^N C(u)e_i = R(u,z)^{-1} c(u,z),
\een
where we used that $\sum_i \Gamma_i = 0$ and $\sum_i \partial_{u_i}
R(u,z)=0$. The former identity follows from the definition of
$\Gamma_i$ and the symmetry $\gamma_{ij}=\gamma_{ji}$, while the
latter is a consequence of \eqref{eq:pf-3}: sum \eqref{eq:pf-3} over
all $i$ and recall that $\sum_i \Gamma_i=\sum_i B_i=0.$
\qed

Using Lemma \ref{le:rec-R} and Lemma \ref{le:rec-c} we get the following
equations for $c(u,0)$:
\beq\label{de:c-1}
\gamma_{ij}:=(\partial_{u_j} c_i(u,0) ) c_j(u,0)^{-1}\quad 
\mbox{is symmetric in $i$ and $j$}, 
\eeq
\beq\label{de:c-2}
[\partial_{u_i}+\Gamma_i,\partial_{u_j}+\Gamma_j]=0,\quad 1\leq
i,j\leq N,
\eeq
\beq\label{de:c-3}
(\partial_{u_1}+\cdots + \partial_{u_N})c(u,0)=0,
\eeq
\beq\label{de:c-4}
E\, c_i(u,0) = (r-1/2) c_i(u,0),\quad 1\leq i\leq N.
\eeq

\begin{proposition}\label{prop:c-sol}
Let $c_i(u,0)$, $1\leq i\leq N$, be a set of functions analytic for
$u\in V$ and such that $c_i(u,0)\neq 0$ for all $i$. The functions $c_i$ solve the
equations \eqref{de:c-1}--\eqref{de:c-4} if and only if they can be
extended to functions $c_i(u,z)$, $1\leq i\leq N$, solving the
equations \eqref{eq:pf-1}--\eqref{eq:pf-6}.
\end{proposition}

\proof
We have already proved that if $c_i(u,z)$, $1\leq i\leq N$, solve
\eqref{eq:pf-1}--\eqref{eq:pf-6}, then $c_i(u,0)$, $1\leq i\leq N$, solve
\eqref{de:c-1}--\eqref{de:c-4}. In the inverse direction, let us
assume that $c_i(u,0)$, $1\leq i\leq N$, solve \eqref{de:c-1}--\eqref{de:c-4}. We need
to check that the series $c_i(u,z)$, $1\leq i\leq N$, defined by the reconstructions of
Lemma \ref{le:rec-R} and Lemma \ref{le:rec-c} solve
\eqref{eq:pf-1}--\eqref{eq:pf-6}. This however is straightforward. 
\qed

\medskip

\emph{Proof of Theorem \ref{t1}.}
According to \cite{Du}, a semi-simple Frobenius structure on $V$ is specified by a set of $N$
functions $c_i(u,0)$, $1\leq i\leq N$, satisfying the system of differential equations defined by 
\eqref{de:c-1}--\eqref{de:c-4} and such that $c_i(u,0)\neq 0$ for all
$u\in V$.  On the other hand, recalling Proposition \ref{prop:c-sol}
we get that there is a one-to-one
correspondence between solutions $c_i(u,0)$ to
\eqref{de:c-1}--\eqref{de:c-4} satisfying $c_i(u,0)\neq 0$ for $u\in
V$ and primitive forms in $\widehat{\mathcal{H}}(V)$. 
\qed

\section{Polynomial primitive forms}\label{sec:pf-poly}

Suppose that $V\subset B$ is an open subset. We say that the
cohomology class $\omega\in \widehat{\H}(V)$ is a {\em
  polynomial class} if for every $p\in V$ we can find an open connected Stein
subset $U\subset V$ such that the restriction $\omega|_U\in \widehat{\H}(U)$ can be
represented by a form 
$$
\sum_{n=0}^{n_0(U)} \omega_U^{(n)} (-z)^n,\quad \omega_U^{(n)} \in
\Omega^{1,\infty}_{\bX/B}(\bpi^{-1}(U)),
$$
depending polynomially on $z$. Note that we do not require the
degrees $n_0(U)$ 
of the polynomials to be uniformly bounded. Nevertheless, we will
prove that we can always choose a covering of $V$ with polynomial representatives
for which the degree of the polynomials are uniformly bounded. If
$\omega$ is both polynomial and  primitive, then we say that $\omega$
is a {\em polynomial primitive  form}.

\subsection{Polynomiality and the sheaf $\H$.}
The main goal in this section is to prove that the polynomial classes
in $\widehat{\H}$ correspond to cohomology classes in $\H$. To make
this statement precise, let us first prove that the map $i: \H\to
\widehat{\H}$ induced by the inclusion
\ben
\bpi_*\Omega^1_{\bX/B}[z]\to \bpi_*\Omega^1_{\bX/B}[\![z]\!]
\een
is injective. 
Recall that the sheaf $\H$ is the sheafification of the presheaf $H_F$
on $B$ defined by 
\ben
H_F(U):= \bpi_*\Omega^1_{\bX/B}(U)[z]/(zd_{\bX/B}+d_{\bX/B}F\wedge)
\bpi_*\Omega^0_{\bX/B}(U)[z].
\een
Similarly, the completion $\widehat{\H}$ is the sheafification of the
presheaf
\ben
\widehat{H}_F(U):=
\bpi_*\Omega^1_{\bX/B}(U)[\![z]\!]/(zd_{\bX/B}+d_{\bX/B}F\wedge)
\bpi_*\Omega^0_{\bX/B}(U)[\![z]\!].
\een
The injectivity of $i$ is a corollary of the following lemma, which
implies that the induced map on the stalks is injective.
\begin{lemma}
Let $U\subset B$ be an open connected subset, then the natural map
$i_F: H_F(U)\to \widehat{H}_F(U)$ is injective.
\end{lemma}
\proof
We have to prove that if 
\ben
(zd_{\bX/B}+d_{\bX/B}F\wedge) \phi = \psi,
\quad \psi\in  \Omega^{1,\infty}_{\bX/B}(\bpi^{-1}(U))[z],
\quad \phi\in  \Omega^{0,\infty}_{\bX/B}(\bpi^{-1}(U))[\![z]\!]
\een
then $\phi$ is polynomial in $z$. Let us write $\phi=\sum_k \phi_k
z^k$ and $\psi=\sum_k \psi_k z^k$ with $\psi_k=0$ for $k>k_0$ for some
positive integer $k_0$. We have to prove that $\phi_k=0$ for $k\gg
0$. 

Since $U$ is connected, the intersection $\bpi^{-1}(U)\cap D_\infty$
has finitely many connected components. Therefore
\ben
H^0(\bpi^{-1}(U), \Omega^{0,\infty}_{\bX/B}) = 
\varinjlim_m H^0(\bpi^{-1}(U), \O_{\bX} (m)) = 
\bigcup_{m=0}^\infty  H^0(\bpi^{-1}(U), \O_{\bX}(m)) .
\een
Therefore, we can choose $m_0$ such that $\phi_{k_0}\in
\bpi_*\O_{\bX}(m_0)(U)$. Note that 
\ben
\phi_{k_0+\ell} = \Big(-\frac{d_{\bX/B}}{d_{\bX/B}F} \Big)^\ell \phi_{k_0}
\een
for all $\ell>0$. Note that the RHS of the above identity is a section
of 
\ben
\bpi_*\O_{\bX}((m_0-\ell)D_\infty) .
\een
But the above pushforward is $0$ for $\ell>m_0$, so $\phi_k=0$ for all
$k>k_0+m_0$. 
\qed

\medskip

For every open subset $U\subset B$ we have the following commutative
diagram
\begin{diagram}
H_F(U) & \rTo^{\rho_F} & \H(U) \\
\dTo^{i_F} &  & \dTo_{i} \\
\widehat{H}_F(U) & \rTo^{\widehat{\rho}_F} & \widehat{\H}(U) 
\end{diagram}
where $\rho_F$ and $\widehat{\rho}_F$ are the natural maps from a
presheaf to its sheafification. If we assume in addition that $U$ is
connected and Stein, then the above lemma implies that $i_F$ is
injective, while Proposition \ref{per-iso} implies that
$\widehat{\rho}_F$ is an isomorphism. We get that both $\rho_F$ and
$i$ must be injective.  The following Proposition is straightforward
to prove by using the properties of the above commutative
diagram.
\begin{proposition}\label{prop:oscil-poly}
Let $V\subset B$ be an open subset. Then a cohomology class $\omega\in
\widehat{\H}(V)$ is polynomial iff $\omega\in \H(V)$.
\end{proposition}

\subsection{Primary differentials}\label{sec:p-dif}
Let us recall the definition of Dubrovin's primary differentials. They are splited into
five types. Following \cite{DNOPS2, Shr} we express them in terms of the
fundamental bi-differential.

\bigskip
{\em Type I.} Normalized Abelian differentials of the second kind on
$\overline{X}_u$, $u\in U$ with poles only at $\infty_i$ ($1\leq i\leq
d$) of order not exceeding the order of the pole of $d_{X/B}F$
\ben
\phi_{t_{i,a}} (p):= \frac{1}{a}\operatorname{res}_{q=\infty_i}
F(q)^{a/m_i} B(q,p),\quad 1\leq i\leq d,\quad 1\leq a\leq m_i-1.
\een

\bigskip
{\em Type II.} Normalized Abelian differentials of the second kind
\ben
\phi_{v_i}(p) := \operatorname{res}_{q=\infty_i}
F(q) B(q,p),\quad 2\leq i\leq d.
\een
Note that $\phi_{v_i}(p)$ has pole only at $p=\infty_i$ and the principal
part is of the form $-d_{X/B}F(p)+\mbox{ regular terms}$.

\bigskip
{\em Type III.}
Normalized Abelian differentials of the third kind 
\ben
\phi_{w_i}(p) := \omega_{\infty_i,\infty_1}(p)= 
\int_{\infty_1}^{\infty_i} B(q,p),\quad 2\leq i\leq d
\een
having poles of order 1 only at $\infty_1$ and $\infty_i$ with
residues respectively $-1$ and $1$. 

\bigskip
{\em Type IV.}
Multi-valued analytic differentials on $\overline{X}_u$ 
\ben
\phi_{r_i}(p) := \frac{1}{2\pi\sqrt{-1}} \, \oint_{q\in \alpha_i} F(q)B(q,p),
\quad 1\leq i\leq g
\een
with increment along the cycle $\beta_j$ equal to 
\ben
\phi_{r_i}(p+\beta_j)-\phi_{r_i}(p) = 
-\delta_{ij} \, d_{X/B}F(p).
\een

\bigskip
{\em Type V.}
The holomorphic 1-forms on $ \overline{X}_u$ 
\ben
\phi_{s_i}(p) := \frac{1}{2\pi\sqrt{-1}} \, \oint_{q\in \beta_i} B(q,p),
\quad 1\leq i\leq g
\een
satisfying the normalization condition 
\ben
\oint_{p\in \alpha_i} \phi_{s_j}(p) = \delta_{ij}.
\een
Every primary differential $\phi$ induces a cohomology class
$[\phi]\in \widehat{\H}$ as follows. Let $X'\subset \bX$ be a tubular
neighborhood of $C$ such that $\phi\in \Omega^1_{\bX/B}(X')$. Then
$[\phi]$ is defined via the isomorphism in Proposition \ref{per-iso},
Part b).
\begin{proposition}\label{prop:pd-poly}
The cohomology classes $[\phi]$ are polynomial in $z$, i.e., $[\phi]\in
\H(B)$. 
\end{proposition}
\proof
All primary differentials, except for Type IV are already in
$\Omega^{1,\infty}_{\bX/B}(\bX)$, so the statement in the lemma is
obvious. Let us assume that $\phi=\phi_{r_i}$ is of type IV. 

Let us fix a positive integer $m>2g-1$. Denote by
$\mathcal{K}$ the coherent sheaf on $\bX$ defined by the exact
sequence
\ben
0\rTo \mathcal{K}\rTo \Omega^1_{\bX/B}(m\infty_1)
\rTo^{\operatorname{res}_{\infty_1}}
\CC_{\infty_1}\rTo 0,
\een
where $\CC_{\infty_1}=(\infty_1)_* \CC_B$ is the pushforward of the
constant sheaf $\CC_B$ on $B$ via the section $\infty_1:B\to \bX$. 
On the other hand we have the following exact sequence of sheaves of
$\bpi^{-1}\O_B$-modules
\ben
0 \rTo \bpi^{-1}\O_B \rTo \Omega^0_{\bX/B}((m-1)\infty_1) \rTo^{d_{\bX/B}} 
\mathcal{K}\to 0.
\een
Pushing forward we get
\ben
0\rTo \O_B\rTo \bpi_* \Omega^0_{\bX/B}((m-1)\infty_1) \rTo \bpi_*
\mathcal{K}\rTo R^1\bpi_* (\bpi^{-1}\O_B)\rTo 0.
\een
We get that $R^1\bpi_* (\bpi^{-1}\O_B)$ is a coherent $\O_B$-modules,
because the remaining sheaves in the above exact sequence are coherent
$\O_B$-modules. Moreover, using the proper base change theorem in sheaf
cohomology (see \cite{Iv}, Theorem 6.2), we get that the stalk of $R^1\bpi_* (\bpi^{-1}\O_B)$ at a
point $u\in B$ is
\ben
H^1(\bX_u,j_u^{-1} \bpi^{-1} \O_B) = H^1(\bX_u,\CC)\otimes \O_{B,u},
\een
where $j_u:\bX_u \to \bX$ is the natural inclusion. Therefore, the
sheaf $R^1\bpi_* (\bpi^{-1}\O_B)$ is the sheaf of holomorphic sections of the vector bundle
on $B$ whose fiber over $u\in B$ is $H^1(\bX_u,\CC)$. Note that this
is a holomorphically trivial bundle because $B$ is Stein and
contractible. Let us construct a trivialization by choosing a basis
of $H^1(\Sigma,\CC)$ (recall that $\bX_{u^\circ}:=\Sigma$ is our
reference fiber) Poincare dual to the basis
$\{\alpha_i,\beta_i\}_{i=1}^g\subset H_1(\Sigma,\CC)$ and using the
parallel transport with respect to the Gauss--Manin connection. On the
other hand, since $B$ is Stein, the above 4-term exact sequence remains exact
when we take global sections. Therefore, we can find global
meromorphic forms (with no residues along $\infty_1$)
$$
\phi_{s_j}^\vee\in
H^0(\bX,\Omega^1_{\overline{X}/B}(m\infty_1)),\quad 1 \leq j\leq g,
$$
that represent the Poincare duals of the cycles $\beta_j$, $1\leq
j\leq g$, i.e., 
\ben
\oint_{\alpha_i} \phi_{s_j}^\vee = 0,\quad \oint_{\beta_i} \phi_{s_j}^\vee = \delta_{ij}.
\een
Let us fix $u\in B$ and define
\ben
\psi_i(p)=\int_{p_1(u)}^p \phi_{s_i}^\vee,\quad 1\leq i\leq g,
\een
where $p_1(u)\in \bX_u$ is the critical point of $F$ corresponding to
the critical value $u_1$. This is a multi-valued analytic function with increment along the
$\beta_j$ cycle 
\ben
\psi_i(p+\beta_j)-\psi_i(p) = \oint_{\beta_j} d\psi_i(p) = \delta_{ij}.
\een
It follows that the 1-form 
\ben
\omega_{r_i}:=\phi_{r_i}(p)+ \psi_i(p) dF(p) + z \phi_{s_i}^\vee(p) =
\phi_{r_i}(p)+(z d+ dF\wedge)\psi_i(p)
\een
is holomorphic for all $p\in \bX_u\setminus{\{\infty_1\}}$ with a
finite order pole at $\infty_1$. Clearly the point-wise construction for
$u\in B$ produces a global form  $\omega_{r_i}\in
\Omega^{1,\infty}_{\bX/B}(\bX)[z] $, which is polynomial in  $z$ of
degree 1 and under the restriction map from Proposition \ref{per-iso},
Part b), the cohomology class of $\omega_{r_i}$ is mapped to the  cohomology
class of the primary differential $\phi_{r_i}$.
\qed

\begin{proposition}\label{primary-pf}
If $\phi$ is a primary differential, then the corresponding cohomology
class $[\phi]\in \widehat{\H}(B)$ is a primitive form on an appropriate open
subset $U\subset B$. 
\end{proposition}
\proof
We claim that the higher-residue pairings
\ben
K([\phi],[\omega_a]) = zc_a,\quad 1\leq i\leq N,
\een
where $\{[\omega_a]\}_{a=1}^N\subset\widehat{\H}(B)$ is the good basis
we have constructed in Section \ref{sec:gb} and 
\ben
c_a = -\operatorname{res}_{p=p_a} \frac{\phi(p)}{t_a(p)}
\een
are holomorphic functions on $B$. The notation in the above formula is the
one we have introduced in the beginning of Section \ref{sec:gb}.  
The proof follows from the following two formulas
\beq\label{primary-unit}
\sum_{a=1}^N z\nabla_{\partial_{u_a}}[\phi]=[\phi]
\eeq
and 
\beq\label{primary-ua}
z\nabla_{\partial_{u_a}} [\phi]=c_a [\omega_a]
\eeq
and Proposition \ref{prop:gb}. Let us prove the above two formulas
when $\phi=\phi_{t_{i,b}}$ is a Type I differential. The argument in
the remaining 4 cases is similar. By definition 
\ben
[\phi] = b^{-1} \int e^{F(p)/z} \operatorname{res}_{q=\infty_i}
\Big(F(q)^{b/m_i} B(q,p)\Big) .
\een
On the other hand, recalling the Rauch's variational formula we get 
\ben
\sum_{a=1}^N \delta_{u_a} B(q,p) = \sum_{a=1}^N \operatorname{res}_{q'=p_a}
\frac{B(q,q')B(q',p)}{d_{X/B}F(q')} = 
-d_p\Big( \frac{B(q,p)}{d_{X/B}F(p)}\Big)
-d_q\Big( \frac{B(q,p)}{d_{X/B}F(q)}\Big),
\een
where in the second equality we have used the residue theorem for
$\overline{X}_u$ and the fact that the integrand has poles only at
$q'=p_a$ ($1\leq a\leq N$), $q'=q$,   and $q'=p$.  Using the above
formula and 
\ben
-\operatorname{res}_{q=\infty_i} 
F(q)^{b/m_i}d_q\Big( \frac{B(q,p)}{d_{X/B}F(q)}\Big)=
\frac{b}{m_i} \operatorname{res}_{q=\infty_i} 
F(q)^{\frac{b}{m_i}-1}B(q,p) =0.
\een
we get  
\ben
\sum_{a=1}^N z\nabla_{\partial_{u_a}}[\phi] = 
 b^{-1} \int e^{F(p)/z} (-zd_p) \operatorname{res}_{q=\infty_i}
\Big(\frac{F(q)^{b/m_i} B(q,p)}{d_{X/B}F(p)}\Big) .
\een
Integrating by parts the RHS of the above formula we get
\eqref{primary-unit}. 

To prove \eqref{primary-ua} we use that 
\ben
\delta_{u_a} B(q,p) = \operatorname{res}_{q'=p_a}
\frac{B(q,q')B(q',p)}{d_{X/B}F(q')} =
\operatorname{res}_{q'=p_a}
\Big(\frac{B(q,q')}{t_a(q')} \Big)\,  
\operatorname{res}_{q'=p_a}
\Big(\frac{B(q',p)}{t_a(q')} \Big)
\een
The above formula implies that 
\ben
z\nabla_{u_a} [\phi] = 
c_a \int e^{F(p)/z} 
(-zd_p) \operatorname{res}_{q'=p_a}
\Big(\frac{\omega_p(q')}{t_a(q')} \Big) = c_a[\omega_a].
\een

We proved that $[\phi]=\sum_{a=1}^N c_a(u) [\omega_a]$, where the
coefficients $c_a$ are independent of $z$. We are going to prove that
the functions $c_a(u)$ satisfy the differential equations
\eqref{de:c-1}--\eqref{de:c-2} and that the matrix $R(u,z)$
reconstructed from $c_a(u)$ via Lemma \ref{le:rec-R} is 1. To this
end, it is enough to prove that $\gamma_{ij}=\beta_{ij}$ and that
$[\phi]$ is homogeneous. Indeed, if we know that
$\gamma_{ij}=\beta_{ij}$, then equation \eqref{de:c-1} is
satisfied. Furthermore, 
\ben
\Gamma_i=\sum_{j:j\neq i} \gamma_{ij}(E_{ij}-E_{ji}) = \sum_{j:j\neq
  i} \beta_{ij}(E_{ij}-E_{ji})  =B_i,
\een
so recalling Proposition \ref{GM:gb-frame} we get that
$
\partial_{u_i}+\Gamma_i 
$
is the leading order term of the Gauss-Manin connection on
$\widehat{\H}$ (written in the frame ${[\omega_i]}_{i=1}^N$):
\ben
\nabla_{\partial_{u_i}} = \partial_{u_i}+B_i + z^{-1} E_{ii}.
\een
The flatness of the Gauss-Manin connection implies equations
\eqref{de:c-2}. The next equation \eqref{de:c-3} is a consequence of
\eqref{primary-unit} and Proposition \ref{GM:gb-frame}.
Finally, it is easy to see that if $\gamma_{ij}=\beta_{ij}$ then in the
reconstruction procedure 
for $R(u,z)$ we have $R_k(u)=0$ for all $k\geq 1$. 

Let us prove that $\gamma_{aj}=\beta_{aj}$ and that $[\phi]$ is
homogeneous. Again, we will do this only for $\phi=\phi_{t_{i,b}}$ a
primary differential of Type I, because the argument in the remaining
cases is similar. First, we need to prove that 
\ben
\partial_{u_a} c_j(u) = \beta_{aj}(u) c_a(u),\quad a\neq j.
\een 
By definition  
\ben
c_j = -b^{-1}
\operatorname{res}_{q=\infty_i}
 \operatorname{res}_{p=p_j} \frac{ F(q)^{b/m_i} B(q,p) }{t_j(p)}.
\een
We are going to use that 
\ben
\partial_{u_a} c_j = -b^{-1}
\operatorname{res}_{q=\infty_i}
 \operatorname{res}_{p=p_j} 
\delta_{u_a}\Big(
\frac{ F(q)^{b/m_i} B(q,p) }{t_j(p)}\Big),
\een
where $\delta_{u_a}$ is the Rauch variational derivative. If $a\neq j$, then we have 
\ben
\delta_{u_a} \Big(\frac{ F(q)^{b/m_i} B(q,p) }{t_j(p)}\Big) = 
\operatorname{res}_{q'=p_a}
\frac{ F(q)^{b/m_i} B(q,q')B(q',p) }{t_j(p)d_{X/B}F(q')}.
\een
On the other hand,
\ben
\operatorname{res}_{q'=p_a}
\frac{B(q,q')B(q',p) }{d_{X/B}F(q')} = 
\operatorname{res}_{q'=p_a}
\Big(\frac{ B(q,q')}{t_a(q')}\Big)  \, 
\operatorname{res}_{q'=p_a}
\Big(
\frac{ B(q',p)}{t_a(q')} \Big).
\een
It remains only to recall the definition \eqref{beta_ai}.

For the homogeneity part, we have to compute 
\ben
(z\nabla_{\partial_z}+\nabla_E )[\phi] =
\Big(z\nabla_{\partial_z}+
\sum_{a=1}^N u_a \nabla_{\partial_{u_a}}\Big)\, [\phi].
\een
Using the Rauch's variational formula we get
that $\sum_a u_a\delta_{u_a} B(q,p)$ is
\ben
\sum_{a=1}^N \operatorname{res}_{q'=p_a}
\frac{F(q')B(q,q')B(q',p)}{d_{X/B}F(q')} = 
-d_p\Big(\frac{ F(p) B(q,p)}{d_{X/B}F(p)}\Big)
-d_q\Big( \frac{F(q) B(q,p)}{d_{X/B}F(q)}\Big).
\een
The contribution of the first term to $z\nabla_E [\phi]$ is
\ben
b^{-1}\int e^{F(p)/z}
\operatorname{res}_{q=\infty_i} F(q)^{b/m_i} (-zd_p)
\Big(
\frac{B(q,p)F(p)}{d_{X/B}F(p)} \Big) = -z^2\nabla_{\partial_z} [\phi],
\een
while the contribution of the second term is
\ben
b^{-1}\int e^{F(p)/z}
\operatorname{res}_{q=\infty_i} F(q)^{b/m_i} (-zd_q)
\Big(
\frac{B(q,p)F(q)}{d_{X/B}F(q)} \Big) =(b/m_i)\, [\phi].
\een
It follows that the form $[\phi_{t_{i,b}}]$ is homogeneous of degree
$r_{t_{i,b}}:=b/m_i$. 
\qed

Theorem \ref{t2} is a direct consequence of Propositions
\ref{prop:pd-poly} and \ref{primary-pf}. For future references, 
let us list the homogeneous degrees of the primary differentials
\ben
&& r_{t_{i,b}}=b/m_i,\quad 1\leq i\leq d,\quad 1\leq b\leq m_i-1, \\
&&
r_{v_i}=r_{r_j}=1,\quad 
r_{w_i}=r_{s_j}=0,\quad
2\leq i\leq d,\quad 1\leq j\leq g.
\een

\subsection{Polynomiality of the $R$-matrix}
Let $V\subset B$ be an open connected subset. Note that every
cohomology class $\omega\in \widehat{\H}(V)$ can be written as 
\ben
\omega=\sum_{i=1}^N c_i [\omega_i],\quad c_i  := K(\omega,[\omega_i])z^{-1},
\een
where the coefficients $c_i\in \O_B(V)[\![z]\!]$. 
\begin{proposition}\label{poly-ci}
The cohomology class $\omega$ is polynomial, if and only if the coefficients
$c_i(u,z)$ $(1\leq i\leq N)$ depend polynomially on $z$.
\end{proposition}
\proof
For fixed $u\in B$ the differentials $\omega_i(p)$ are multi-valued analytic on $X_u$ with
increment along the cycle $\beta_j$ given by 
\ben
\omega_i(p+\beta_j)-\omega_i(p) = c_{ij}(u) dF(p),
\een
where 
\ben
c_{ij}(u) = \operatorname{res}_{q=p_i}
\frac{\theta_j(q)}{t_i(q)},\quad \theta_j(q):=\oint_{p\in\beta_j} B(p,q). 
\een
Recall that in the proof of Proposition
\ref{prop:pd-poly} we have constructed multivalued analytic functions
$\psi_i(p)$ on $\bX$ with finite order poles along the divisor $D_\infty$, such that
the 1-forms $d_{\bX/B}\psi_i\in H^0(\bX, \Omega^{1,\infty}_{\bX/B})$,
have vanishing residues along $D_\infty$, and 
\ben
\oint_{p\in \alpha_i} d\psi_j(p) = 0,\quad \oint_{p\in \beta_i}
d\psi_j(p)=\delta_{ij}, \quad 1\leq i,j\leq g.
\een
We get that under the restriction map defined by Proposition
\ref{per-iso}, Part b), the cohomology class of $\omega_i(p)$ is
identified with the cohomology class of 
\ben
\omega_i - \sum_{j=1}^g c_{ij}(u)  \, (zd_{\bX/B} +d_{\bX/B}F\wedge
)\, \psi_j(p).
\een
Note that the form on the RHS is analytic and polynomial in $z$ of
degree at most 1, so the cohomology classes $[\omega_i]$ are
polynomial. Therefore, we need only to prove that if $\omega$ is
polynomial, then $c_i(u,z)$ are polynomial in $z$.

Let us fix $u_0\in V$ and an open connected Stein neighborhood $U\subset
V$ of $u_0$, s.t., $\omega$ can be represented by a polynomial form
\ben
\omega_U=\sum_{n=0}^{n_0(U)} \omega^{(n)}_{U} z^n,\quad 
\omega_U^{(n)} \in H^0(\bpi^{-1}(U), \Omega^1_{\bX/B}(m\,D_\infty)),
\een
where $m$ is a sufficiently large integer.  
Note that if $\theta\in H^0(\bpi^{-1}(U),
\Omega^1_{\bX/B}(m\,D_\infty))$, then since $\theta$ has a finite
order pole along $D_\infty$, we can express the singular part of
$\theta$ along $D_\infty$ in terms of linear combinations of finitely
many expressions of the type $F(p)^n\phi(p)$, where $n\geq 0$ and
$\phi$ is a primary differential of type I, II, or III. Therefore,
there are polynomials  $a_\phi(\lambda)\in \O_B(U) [\lambda]$, s.t., 
\ben
\theta(p)-\sum_{\phi\in I\cup II\cup III} a_\phi(F(p)) \phi(p),
\een
is a holomorphic 1-form on $\overline{X}_u$ and hence it can be 
written as a linear combination of primary differentials of type V. 
Using this observation, we get that 
\ben
[\omega_U] = \sum_\phi [b_\phi(z,F) \phi] = \sum_\phi b_\phi(z,-z^2\nabla_{\partial_z})[\phi],
\een
where the sum is over primary differentials of Type I, II, III, and V,
$b_\phi\in \O_B(U)[z,\lambda]$, and for the second equality we used that 
\ben
-z^2\nabla_{\partial_z} \int e^{F(p)/z}\phi(p) = \int
e^{F(p)/z}F(p)\phi(p).
\een
We get that
\ben
[\omega_U]= 
\sum_\phi\sum_{i=1}^N 
c_{\phi,i} b_\phi(z,-z^2\nabla_{\partial_z}) [\omega_i] ,
\een 
where the coefficients $c_{\phi,i}$ are independent of $z$. On the
other hand, according to Proposition \ref{GM:gb-frame} (see the
notation in Section \ref{sec:pn}) we have 
\ben
-z^2\nabla_{\partial_z} \omega^{\rm gb} = 
\omega^{\rm gb} \Big(\Big(\frac{1}{2}-\sum_{i=1}^N u_i
B_i(u)\Big)z-\sum_{i=1}^N u_iE_{ii}\Big),
\een 
where $\omega^{\rm gb}:=([\omega_1],\dots,[\omega_N])$. 
Therefore, $b_\phi(z,-z^2\nabla_{\partial_z}) [\omega_i]$ is a linear
combination of $[\omega_j]$, $1\leq j\leq N,$ with coefficients
depending polynomially on $z$. We get that 
\ben
\omega_U=\sum_{i=1}^N c_i^U \,[\omega_i],
\een
where $c_i^U\in \O_B(U)[z].$ Therefore, $c_i|_U=c_i^U$ is polynomial
in $z$ of some degree $n(U)$. Writing $c_i=\sum_{n=0}^\infty c_{i,n}
z^n$, we get that $c_{i,n}=0$ for all $n>n(U)$. 
\qed

Note that Proposition \ref{poly-ci} can be reformulated as
follows. The map 
\ben
\H\to \O_B[z]^{\oplus N},\quad 
\omega\mapsto (c_1,\dots,c_N),\quad c_i:=K(\omega,\omega_i)
\een
is an isomorphism of $\O_B[z]$-modules and $\{[\omega_i]\}_{i=1}^N$ is
an $\O_B[z]$-basis. Therefore, we proved Corollary \ref{c1}.

\smallskip
Let $c(u)=(c_1(u),\dots,c_N(u))$ be an arbitrary solution to the
equations \eqref{de:c-1}--\eqref{de:c-4}. Let us denote by $R(u,z)$
the matrix reconstructed from $c(u)$ via the algorithm of Lemma
\ref{le:rec-R}. Then we have the following corollary  of Proposition
\ref{poly-ci}. 
\begin{corollary}
The primitive form corresponding to $c(u)$ is polynomial if and only
if the entries of the matrix $R(u,z)$ depend  polynomially on $z.$ 
\end{corollary}

In the next section we are going to formulate the problem of relating
the correlator forms defined by the Eynard--Orantin recursion with the
higher-genus invariants of a 
semi-simple Frobenius manifold. The case for which $R(u,z)=1$ plays a
key role. Note that if $R(u,z)=1$, then according to Lemma
\ref{le:rec-c} the coefficients $c_i(u,z)=c_i(u,0)$ are independent of
$z$. Polynomial primitive forms for which $R(u,z)=1$ can be classified
in terms of primary differentials. We have the following corollary.
\begin{corollary}\label{cor:z-indep}
Let $\omega\in \widehat{\H}$ be a cohomology class of 
homogeneous degree $r$. Then $\omega$ is a primitive
form whose matrix $R(u,z)=1$ if and only if 
\ben
\omega=\sum_{\phi:\operatorname{deg}(\phi)=r} a_\phi [\phi],
\een
where the sum is over all primary differentials of homogeneous degree
$r$ and $a_\phi$ are constants independent of $u$ and $z$.
\end{corollary}
\proof
If $R(u,z)=1$ then $R_1(u)=0$, so $\gamma_{ij}(u)=\beta_{ij}(u)$. 
It remains only to use that the space of solutions $c(u)=(c_1(u),\dots,c_N(u))$ to
the system of differential 
equations \eqref{de:c-1}--\eqref{de:c-3} satisfying the additional
constraint $\gamma_{ij}=\beta_{ij}$ is $N$-dimensional and a basis is
given by 
\ben
c_\phi(u)=(c_{\phi,1}(u),\dots,c_{\phi,N}(u)),\quad 
c_{\phi,i}(u) = -\operatorname{res}_{p=p_i}\frac{\phi(p)}{t_i(p)},
\een 
where $\phi$ runs through the set of primary differentials.
\qed 

\subsection{Two-dimensional Frobenius manifolds}

We would like to classify all two dimensional semi-simple Frobenius
manifolds that correspond to the polynomial primitive forms of the
Hurwitz spaces. To begin with, we need to classify all coverings  
$f:\Sigma\to \mathbb{P}^1$ for which 
\ben
N=2g-2+d+\sum_{i=1}^d m_i = 2.
\een
It is easy to see that $g=0$, so $\Sigma=\PP^1$ and that for the
covering map $f$, up to isomorphism there are only two
cases. Either $d=1$, $m_1=3$, and
\ben
f(x)=\frac{x^3}{3} +s_1 x +s_2
\een 
or $d=2$, $m_1=m_2=1$, and 
\ben
f(x)=x+s_1 x^{-1}+s_2,
\een
where $x$ is the coordinate function on
$\CC=\Sigma\setminus{\{\infty\}}$ and  
$(s_1,s_2)\in \CC^*\times \CC$ are some parameters. 

\subsubsection{Case 1}
The universal cover $B=\CC^2$, the family $X=B\times \CC$ and 
\ben
F:X\to \CC,\quad F(t,x) = \frac{x^3}{3}+e^{t_1} x + t_2.
\een
The critical points and the corresponding critical values of $F$ are
respectively 
\ben
p_{1,2}=\pm \sqrt{-1}\, e^{t_1/2}\quad
\mbox{and}\quad
u_{1,2} = t_2\pm \frac{2\sqrt{-1}}{3}\, e^{3t_1/2}.
\een
The good basis (see Section \ref{sec:gb}) takes the form
\ben
\omega_i =-(x+p_i)\frac{dx}{\sqrt{2p_i}},\quad i=1,2,
\een
where we have to make an additional choice of a sign of
$\sqrt{p_i}$. Let us write  
\ben
p_i=e^{(t_1+(2i-1)\pi\sqrt{-1})/2},\quad 
i=1,2, 
\een
then we define
\ben
\sqrt{p_i} := e^{(t_1+(2i-1)\pi\sqrt{-1})/4},\quad i=1,2.
\een
Let us assume that 
\ben
c_1(u,z)\omega_1+c_2(u,z)\omega_2
\een
is a polynomial primitive form of homogeneous degree $r$. Equations
\eqref{de:c-3} and \eqref{de:c-4} imply 
that 
\ben
c_i(u,0) = c_i^\circ (u_1-u_2)^{r-1/2},\quad i=1,2,
\een
where $c_i^\circ$ are some constants. Since 
\ben
\gamma_{1,2}(u) & = & -\frac{r-1/2}{u_1-u_2}\, (c_1^\circ/c_2^\circ)\\
\gamma_{2,1}(u) & = & -\frac{r-1/2}{u_1-u_2}\, (c_2^\circ/c_1^\circ)
\een
we get that $c_1^\circ /c_2^\circ = \pm \sqrt{-1}$. It is enough to
classify primitive forms up to a constant factor, so we may assume
that $c_1^\circ =\sqrt{-1}$ and $c_2^\circ=1$. 

It remains only to find the $R$-matrix reconstructed from $c_i(u,0)$
according to the algorithm of Lemma \ref{le:rec-R}. We already know
that 
\ben
\Gamma_1=-\Gamma_2= 
\frac{a}{u_1-u_2}
\begin{bmatrix}
0 & 1\\
-1 & 0
\end{bmatrix},\quad 
a := -\sqrt{-1} \Big(r-\frac{1}{2}\Big).
\een
We need to compute the matrices $B_1$ and $B_2$.  
The Riemann's second fundamental form is 
\ben
B(x_1,x_2)=\frac{dx_1dx_2}{(x_1-x_2)^2}.
\een
After a straightforward computation we get that the coefficients
\eqref{beta_ai} are 
\ben
\beta_{1,2}(u) = 
\frac{b}{u_1-u_2},\quad b:=\sqrt{-1}/6.
\een
Therefore,
\ben
B_1=-B_2 = 
\frac{b}{u_1-u_2}
\begin{bmatrix}
0 & 1\\
-1 & 0
\end{bmatrix}.
\een
After a straightforward computation we get 
\ben
R(u,z) = 1+\sum_{k=1}^\infty R_k^\circ z^k (u_1-u_2)^{-k},
\een
where 
\ben
R_k^\circ  =a_k 
\begin{bmatrix}
(a+b(-1)^{k-1} & k \\
(-1)^{k-1} k & -b+a(-1)^{k-1} 
\end{bmatrix}
\een
where the numbers $a_k$ $(k\geq 1)$ are defined recursively by
$a_1=a-b$ and 
\ben
a_{k+1} = \frac{1}{k+1}(a^2+b^2+k^2+2ab (-1)^{k-1}) a_k,\quad k\geq 1.
\een
Since the matrix $R(u,z)$ is polynomial in $z$, there exists an
integer $m>0$, such that $a_k=0$ for all $k\geq m$ and $a_{m-1}\neq
0$. Therefore,
\ben
a=\Big(\frac{(-1)^m}{6} \pm m\Big)\sqrt{-1}
\een
Recalling also that $D=1-2r$ is the conformal dimension, we get that
all two-dimensional semi-simple Frobenius manifolds of conformal
dimension 
\ben
D=\frac{(-1)^n}{3} + 2n,\quad n\in \ZZ,
\een
correspond to polynomial primitive forms. 

\subsubsection{Case 2}
The universal cover $B=\CC^2$ and the family $X=B\times \CC^*$, where
$\CC^*=\CC\setminus{\{0\}}$. The function 
\ben
F(t,x) = x+e^{t_1}/x + t_2.
\een
The critical points and the corresponding critical values are
\ben
p_{1,2} = \pm e^{t_1/2},\quad u_{1,2} = t_2\pm 2 e^{t_1/2}.
\een 
The good basis takes the form
\ben
\omega_i=-\frac{p_i(x+p_i)}{x^2} \, \frac{dx}{\sqrt{2p_i}},
\een
where we define
\ben
\sqrt{p_i} :=e^{(t_1+2\pi(i-1)\sqrt{-1})/4},\quad i=1,2.
\een
The coefficient \eqref{beta_ai} is
\ben
\beta_{1,2}=\frac{b}{u_1-u_2},\quad b:=\sqrt{-1}/2.
\een
The formulas for the matrix $R(u,z)$ are the same as in the previous
case, except that the value of the constant $b$ now is $\sqrt{-1}/2$
(instead of $\sqrt{-1}/6$). We get that the polynomial primitive forms
in this case correspond to semi-simple Frobenius manifolds of
conformal dimension 
\ben
D=(-1)^n +2n,\quad n\in \ZZ,
\een
i.e., $D$ is an odd integer. 
\begin{remark}
The Frobenius structures of $A_2$-singularity and of quantum cohomology
of $\PP^1$ have conformal dimensions respectively $\frac{1}{3}$ and
$1$, so they do correspond to polynomial primitive forms. On the other
hand, the Frobenius structures on the orbit space of a finite reflection
group of type $I_2(k)$ has conformal dimension
$D=1-\frac{2}{k}$. Therefore, if $k>3$ (note that type $I_2(3)$
coincides with type $A_2$), the
Frobenius structure does not correspond to a polynomial primitive form. 
\end{remark}

\section{Topological recursion and semi-simple Frobenius structures}\label{sec:tr-pf}

Let us recall the notation of Section \ref{sec:br_cov}. 
Suppose that $\omega=\sum_{n\geq 0} \omega^{(n)}(-z)^n$ is a
polynomial primitive form on $\H(U)$ for some contractible open subset
$U\subset B$ containing the point $u^\circ$. In particular, $U$ is
equipped with a semi-simple Frobenius structure. Using the
Kodaira--Spencer isomorphism (see Section \ref{sec:KS-iso}) we identify the tangent space
$T_{u^\circ}U$ with the algebra of functions on the critical scheme of $f$
\ben
H:=\Gamma(X_{u^\circ}, 
\mathcal{O}_{X_{u^\circ}}/
\mathcal{O}_{X_{u^\circ}}(-p^\circ_1-\cdots -p^\circ_N)
),
\een
where $p_i^\circ$ are the zeros of $df$ in $X_{u^\circ}=\Sigma\setminus f^{-1}(\infty)$. 
Let us trivialize the tangent and the co-tangent bundles 
\ben
T^*U\cong TU\cong U\times T_{u^\circ}U\cong U\times H,
\een
where the first isomorphism uses the Frobenius pairing, the second
one uses the Levi--Civita connection of the Frobenius pairing, and the
last one is the Kodaira--Spencer isomorphism.

\subsection{The periods of the Frobenius structure}\label{sec:periods}

Given $(u,\lambda)\in B\times \CC$, we denote 
\ben
X_{u,\lambda}=\{p\in X_u\ |\ \varphi(p)=(u,\lambda)\}.
\een
The set of $(u,\lambda)$ such that the number of points in
$X_{u,\lambda}$ is not $m_1+\cdots +m_d$ (the degree of the covering
$X_u\to \CC$) form an analytic hypersurface in $B\times \CC$, called
the {\em discriminant}. For every open subset $U\subseteq B$, we put
$(U\times \CC)'$ for the complement to the discriminant in $U\times
\CC$. The relative homology groups  
\ben
H_1(X_u,X_{u,\lambda};\CC),\quad (u,\lambda)\in (B\times \CC)'
\een 
form a rank $N$ vector bundle on $(B\times \CC)'$ equipped with a flat
Gauss--Manin connection.

Let us fix a reference point $(u^\circ,\lambda^\circ)\in (U\times
\CC)'$. For every relative cycle 
\ben
\alpha \in \mathfrak{h}:=H_1(X_{u^\circ},X_{u^\circ,\lambda^\circ};\CC)
\een 
and every integer $n$ we define the multi-valued analytic function $I^{(n)}_\alpha:(U\times
\CC)'\to H$ as follows. First, for $\ell\geq 0$ we define
\ben
I^{(-\ell)}_\alpha(u,\lambda):= 
-d_u \sum_{n\geq 0} \int_{\alpha_{u,\lambda}}
\frac{(\lambda-F(p))^{n+\ell} }{(n+\ell)!} \omega^{(n)}\quad \in
T^*_uU\cong H,
\een
where the value of the RHS depends on the choice of a reference path
in $(U\times \CC)'$ and $\alpha_{u,\lambda}\in
H_1(X_u,X_{u,\lambda};\CC)$ is the relative cycle obtained from
$\alpha$ via a parallel transport along the reference path. Note that 
\ben
\partial_\lambda I^{(-\ell)}_\alpha(u,\lambda) =
I^{(-\ell+1)}_\alpha(u,\lambda),\quad \forall \ell>0.
\een
For $n\geq 0$, we define
\ben
I^{(n)}_\alpha(u,\lambda) = \partial_\lambda^{n} I^{(0)}_\alpha(u,\lambda).
\een
We will refer to $I^{(n)}_\alpha$ as {\em periods}. Their relation to
the oscillatory integrals (see Section \ref{sec:Saito_str}) is the following.
Put
\ben
J_\Gamma(u,z):= (-2\pi z)^{-1/2} \, (z d_u)\, \int_\Gamma e^{F(p)/z} \omega
\quad \in T^*_uU\quad \cong H.
\een
Let us choose the cycle $\Gamma$ to be the Lefschetz thimble
$\Gamma_i$ consisting of points $p\in X_u$, such that the gradient
trajectory through $p$ of the Morse function
$-\operatorname{Re}(F(p)/z)$ flows out of the critical point $p_i$. The
image of $\Gamma_i$ via the map $F:X_u\to \CC$ is a smooth path
starting at the critical value $u_i$ and approaching $\infty$ in such
a way that $\operatorname{Re}(\lambda/z)\to -\infty$ as $\lambda$
approaches $\infty$ along $F(\Gamma_i)$. If $\lambda\in F(\Gamma_i)$
then let us denote by $\gamma^{(i)}_\lambda\in
H_1(X_u,X_{u,\lambda};\ZZ)$ the cycle obtained from $\Gamma_i$ by
truncating all points $p\in \Gamma_i$, such that, 
\ben
\operatorname{Re}(F(p)/z) < \operatorname{Re}(\lambda/z).
\een
We have
\ben
J_{\Gamma_i}(u,z) = (-2\pi z)^{-1/2} \int_{u_i}^\infty 
e^{\lambda/z} I^{(0)}_{\gamma^{(i)}}(u,\lambda)d\lambda. 
\een
Using integration by parts, we also get that 
\ben
J_{\Gamma_i}(u,z)=
\frac{1}{\sqrt{2\pi}}\, (-z)^{\ell-\frac{1}{2}} \, \int_{u_i}^\infty 
e^{\lambda/z} I^{(-\ell)}_{\gamma^{(i)}}(u,\lambda)d\lambda 
\een
for all $\ell\geq 0$.

Let us fix a basis $\{\phi_i\}_{i=1}^N\subset H$ and denote by
$t=(t_1,\dots,t_N)$ the flat coordinate system on $U$, such that
$\partial/\partial t_i=\phi_i$ via the identification $T_uU\cong H$
and $t_i(u^\circ)=t_i^\circ$, $1\leq i\leq N$, where the choice of
$t_i^\circ$ will be specified later on. The Frobenius multiplication
$\bullet$ in $T_uU$ gives rise (via $T_uU\cong H$) to a Frobenius multiplication
$\bullet_u$ in $H$ for every $u\in U$, while the operator $\theta$
(see \eqref{def:theta}) gives rise to a linear operator in $H$
independent of $u$. The oscillatory integrals satisfy the following
system of differential equations 
\beqa\label{qde:i}
z\partial_{t_i} J(u,z)  & = & \phi_i\bullet_u J(u,z) ,\quad 1\leq
i\leq N,\\
\label{qde:z}
(z\partial_z +E) J(u,z) & = & \theta J(u,z).
\eeqa
The periods $I^{(n)}_\alpha(u,\lambda)$ satisfy the system of
differential equations obtained from the above one via the Laplace
transform
\ben
\partial_{t_i}I^{(n)}_\alpha(u,\lambda) & = & 
-\phi_i\bullet_u I^{(n+1)}_\alpha(u,\lambda),\quad 1\leq i\leq N,\\
(\lambda-E\bullet_u) \partial_\lambda I^{(n)}_\alpha(u,\lambda) & = & 
\Big(\theta-\frac{1}{2}-n\Big)I^{(n)}_\alpha(u,\lambda).
\een
The above formulas define a flat connection known as the {\em second
  structure connection}. 

\subsection{Stationary phase asymptotic}\label{sec:spa}

The primitive form can be written as 
\ben
\omega = \sum_{i=1}^N c_i(u,z) [\omega_i],
\een
where $c_i(u,z)$ depend polynomially on $z$ (see Proposition
\ref{poly-ci}). Let $c(u,z)$ be the vector column with entries
$c_i(u,z)$ and $R_\omega(u,z)$ be the matrix uniquely determined from
$c_i(u,0)$ according to Lemma \ref{le:rec-R}. In particular,
$c(u,z)=R_\omega(u,z)c(u,0)$. 

Recall the stationary phase asymptotic
\ben
J_{\Gamma_i}(u,z) \sim \Big(J_{i,0}(u)+J_{i,1}(u) z+\cdots \Big)
e^{u_i/z}\quad z\to 0,
\een 
where $J_{i,k}(u)=\sum_{a=1}^N J_{i,k}^a(u)\phi_a$. Let $J_k(u)$ be the
matrix whose $(a,i)$-entry is $J_{i,k}^a(u)$. The matrix
$\Psi(u):=J_0(u)$ is essentially the Jacobian of the change of flat to
canonical coordinates, i.e., 
\ben
\Psi_{ai}(u) = \sum_{b=1}^N 
-c_i(u,0) \eta^{ab}\frac{\partial u_i}{\partial t_b},
\quad 1\leq a,i\leq N,
\een
where $\eta^{ab}:=(\phi^a,\phi^b)$ and $\{\phi^a\}_{a=1}^N\subset H$
is the basis dual to $\{\phi_a\}_{a=1}^N\subset H$ with respect to the
Frobenius pairing. The matrices $J_k(u)$ with $k>0$ can be recovered
recursively from $\Psi$ by using the differential equations for the
oscillatory integrals. Put 
\ben
R(u,z)=1+R_1(u)z+R_2(u)z^2+\cdots,\quad R_k(u):=\Psi(u)^{-1} J_k(u).
\een
This matrix satisfies the symplectic condition $R(u,z)R(u,-z)^T=1$ and
it is the matrix that defines Givental's total ancestor potential of
the semi-simple Frobenius manifold (see \cite{G1,G2}). 
\begin{proposition}\label{prop:R}
Let $R_\Sigma(u,z)$ be the matrix whose $(a,i)$ entry is
\ben
-(-2\pi z)^{-1/2} \int_{p\in \Gamma_i}
e^{(F(p)-u_i)/z} \omega_a(p).
\een
Then $R(u,z)=R_\omega(u,z)^T\, R_\Sigma(u,z)$.
\end{proposition}
\proof
We follow the notation in Section \ref{sec:pn} and
\ref{sec:da-constr}.  
By definition 
\ben
(J_i(u,z),\partial_{u_j}) = -z\partial_{u_j} \Big( c(u,0)^T
R_\omega(u,z)^T R_\Sigma(u,z)e^{u_i/z}\Big)_i,
\een
where we used Lemma \ref{le:rec-c}. Let us recall the differential
equations for $c(u,0)$, $R_\omega(u,z)$ (see Proposition
\ref{prop:pf}), and $R_\Sigma(u,z)$ (see Section \ref{sec:pn})
\ben
\partial_{u_j} c(u,0) & = &  -\Gamma_j(u)c(u,0)\\
\partial_{u_j}R_\omega(u,z) & = & z^{-1}[R_\omega(u,z),E_{jj}]+
(R_\omega(u,z)\Gamma_j(u)-B_j(u) R_\omega(u,z)),\\
\partial_{u_j}R_\Sigma(u,z) & = &  z^{-1}[E_{jj},R_\Sigma(u,z)] -B_j(u) R_\Sigma(u,z).
\een
Using these differential equations we get
\ben
(J_i(u,z),\partial_{u_j}) = -c_j(u,0)(R_\omega^T\, R_\Sigma)_{ji}\, e^{u_i/z}.
\een
On the other hand
\ben
J_i(u,z) = \sum_{a=1}^N (\Psi(u) R(u,z))_{ai} \phi_a e^{u_i/z} = \sum_{a,k=1}^N
  \phi_a \Psi_{ak}(u) R_{ki}(u,z)e^{u_i/z}.
\een
It remains only to observe that $\sum_a \phi_a \Psi_{ak}(u) = -c_k(u,0)
du_k$ under the identification $T^*U\cong U\times H$.
\qed

Following \cite{Mi1}, we define the
total ancestor potential using the local Eynard--Orantin recursion
(see also \cite{BOSS}). The total ancestor potential is a formal
series of the type
\ben
\mathcal{A}_u(\hbar,\mathbf{t}) = 
\exp\Big(\sum_{g,n=0}^\infty \frac{\hbar^{g-1}}{n!}
\langle \mathbf{t}(\overline{\psi}),\dots,\mathbf{t}(\overline{\psi})\rangle_{g,n}\Big),
\een
where $\mathbf{t}(\overline{\psi})=\sum_{k=0}^\infty \sum_{a=1}^N
t_{k,a}\phi_a\overline{\psi}^k$ with $t_{k,a}$ formal variables and the
correlators
\ben
\langle
\phi_{a_1}\overline{\psi}^{k_1},\dots,\phi_{a_n}\overline{\psi}^{k_n}\rangle_{g,n}
\een
are non-zero only if they are stable (i.e. $2g-2+n>0 $) and are
defined by the following recursion  
\ben
&&
\langle \phi_a\overline{\psi}^k,\mathbf{t},\dots,\mathbf{t}\rangle_{g,n+1} = 
\frac{1}{4} \sum_{i=1}^N 
\operatorname{res}_{\lambda=u_i}
\frac{(I^{(-k-1)}_{\gamma^{(i)}}(u,\lambda),\phi_a)}{
(I^{(-1)}_{\gamma^{(i)}}(u,\lambda),1)}\times \\
&&
\Big( \langle
\phi^+_{\gamma^{(i)}}(u,\lambda;\overline{\psi}),\phi^+_{\gamma^{(i)}}(u,\lambda;\overline{\psi}),
\mathbf{t},\dots,\mathbf{t}\rangle_{g-1,n+2} + 
\sum_{\substack{ g'+g''=g \\ n'+n''=n}}   \\
&&
{n\choose n'}
\langle
\phi^+_{\gamma^{(i)}}(u,\lambda;\overline{\psi}),
\mathbf{t},\dots,\mathbf{t}\rangle_{g',n'+1} 
\langle
\phi^+_{\gamma^{(i)}}(u,\lambda;\overline{\psi}),
\mathbf{t},\dots,\mathbf{t}\rangle_{g'',n''+1} 
 \Big)
\een
where the insertion $\mathbf{t}$ should be understood as $\mathbf{t}(\overline{\psi})$, 
\beq\label{bosonic-field}
\phi_\alpha(u,\lambda;z):=\sum_{n\in \ZZ} I^{(n+1)}_\alpha(u,\lambda)(-z)^n,
\eeq
$\phi_\alpha^+(u,\lambda;z)$ is obtained from
$\phi_\alpha(u,\lambda;z)$ by truncating all terms containing negative
powers of $z$, and all unstable correlators are by definition 0, except
for 
\ben
\langle\phi_\alpha^+(u,\lambda;\overline{\psi}),\mathbf{t}\rangle_{0,2}:=
\sum_{m=0}^\infty \sum_{a=1}^N
(I^{(-m)}_\alpha(u,\lambda), \phi_a) t_{m,a}.
\een
and 
\ben
\langle
\phi_\alpha^+(u,\lambda;\overline{\psi}),
\phi_\alpha^+(u,\lambda;\overline{\psi})
\rangle_{0,2} = 
\frac{1}{2}((\lambda-E\bullet_u) I^{(1)}_{\gamma^{(i)}}(u,\lambda),
I^{(1)}_{\gamma^{(i)}}(u,\lambda)).
\een
We have the following formula for the Laurent series expansion of the
periods
\beq\label{period-exp}
I^{(-n)}_{\gamma^{(i)}}(u,\lambda) =
\sqrt{2\pi}\, 
\sum_{k=0}^\infty (-1)^k \Psi R_k(u) e_i\,  
\frac{(\lambda-u_i)^{k+n-1/2}}{\Gamma(k+n+1/2)},
\eeq
where $e_i$ is the column vector with 1 on the i-th position and 0s
elsewhere and the period on the LHS is identified with a column vector
whose entries are the coordinates of the period with respect to the
basis $\{\phi_i\}_{i=1}^N\subset H$. Note that the recursion 
determines the correlators in terms of the matrix $R(u,z).$

\subsection{The local EO recursion}
Let us fix a point $u^\circ\in B$ and a small neighborhood $U\subset
B$ of $u^\circ$. Let us denote by $X_U:=\pi^{-1}(U)$ the restriction
to $U$ of the family $\pi:X\to B$. Let $X_U^{\rm loc}\subset X_U$ be a
tubular neighborhood of the relative critical set $C\cap X_U$. Using
the local Morse coordinates $t_i$ ($1\leq i\leq N$) of $F:X\to \CC$ we identify $X_U^{\rm
  loc}\cong U\times \Big(\widetilde{\Delta}_1\sqcup\cdots \sqcup \widetilde{\Delta}_N\Big)$ where each
$\widetilde{\Delta}_i=\{|t_i|<\epsilon_i\}\subset \CC$ is a 
sufficiently small disk. The restriction of $F$ to $U\times
\widetilde{\Delta}_i$ takes the form
\ben
F(u,t_i) = u_i+\frac{1}{2} t_i^2,\quad (u,t_i)\in U\times \widetilde{\Delta}_i.
\een
Following \cite{Mi1} we introduce the local EO recursions. It is defined in terms of a
set of symmetric holomorphic forms
\ben
\omega^i\in \Omega^1_{U\times \widetilde{\Delta}_i/U}(U\times
\widetilde{\Delta}_i),\quad 
\omega^{ij}\in 
\Omega^1_{U\times \widetilde{\Delta}_i/U}\boxtimes
\Omega^1_{U\times\widetilde{\Delta}_j/U}(  U\times (\widetilde{\Delta}_i\times
\widetilde{\Delta}_j-\widetilde{\Delta}_{ij})),
\een
where $\widetilde{\Delta}_{ij}=\emptyset$ if $i\neq j$ and
$\widetilde{\Delta}_{ii}=\{(s,t)\in \widetilde{\Delta}_i\times
\widetilde{\Delta}_i\ |\ s^2=t^2\}$
such that 
\ben
P^i(u,\lambda) = \sum_{k=0}^\infty P^i_k(u)(\lambda-u_i)^{k+1/2} 
\een
and
\ben
P^{ij}(u,\lambda_1,\lambda_2) & = &
\frac{\delta_{ij}}{(\lambda_1-\lambda_2)^2} \, \left(
\frac{(\lambda_1-u_i)^{1/2}}{(\lambda_2-u_j)^{1/2}} +
\frac{(\lambda_2-u_j)^{1/2}}{(\lambda_1-u_i)^{1/2}} \right) + \\
&&
\sum_{k,\ell=0}^\infty P^{ij}_{k,\ell}(u)
(\lambda_1-u_i)^{k-1/2}(\lambda_2-u_j)^{\ell-1/2},
\een 
where the above series are defined by
$$
\omega^i(u,t_i) = P^i(u,\lambda) d\lambda,\quad \lambda=u_i+t_i^2/2
$$
and 
$$
\omega^{ij}(u,t_i,t_j) = P^{ij}(u,\lambda_1,\lambda_2)d\lambda_1\cdot
d\lambda_2,\quad 
\lambda_1=u_i+t_i^2/2,\quad \lambda_2 = u_j+t_j^2/2,
$$
where if $i=j$, then $t_j$ should be interpreted as a second copy of $t_i$.
Note that the requirement that $\omega^{ij}$ is symmetric is
equivalent to $P^{ij}(u,\lambda_1,\lambda_2) =
P^{ji}(u,\lambda_2,\lambda_1)$. We will need also the expansion 
\ben
P^{ij}(u,\lambda_1,\lambda_2)=\frac{2\delta_{ij}}{(\lambda_1-\lambda_2)^2}+
\sum_{k=0}^\infty P^{ij}_k(u,\lambda_1) (\lambda_2-\lambda_1)^k.
\een

Put $\Delta_i:=\{(u,\lambda)\in U\times \CC\ |\
|\lambda-u_i|<\epsilon_i^2/2\}$ and let $\Delta_i(u)\subset \CC$ be
the fiber over $u\in U$ of the projection map $\Delta_i\to U$. 
Both $P^i(u,\lambda)$ and $P^{ij}(u,\lambda_1,\lambda_2)$ are
multi-valued holomorphic functions respectively on $\Delta_i$ and
$\Delta_i\times_U \Delta_j$. In order to keep track of their values we
introduce local systems $\cL_i$ on 
$$
\Delta_i^*:= \{(u,\lambda)\in U\times \CC\ |\
0<|\lambda-u_i|<\epsilon_i^2/2\},\quad 1\leq i\leq N.
$$ 
The sections of $\cL_i$ over some open neighborhood $V\subset \Delta^*_i$ are
given by the holomorphic branches of
$(\lambda-u_i)^{1/2}$ on  $V$. Alternatively, we 
introduce a line bundle $L_i\to \Delta_i^*$ whose fiber over a point
$(u,\lambda)\in \Delta_i^*$ is the relative homology group
\ben
H_1(\widetilde{\Delta}_i,F(u,\cdot)^{-1}(\lambda)\cap \widetilde{\Delta}_i;\CC)\cong \CC.
\een 
The line bundle $L_i$ is equipped with a flat Gauss-Manin connection
and intersection pairing $(\alpha|\beta):=\partial \alpha \circ \partial\beta$.
The local system $\cL_i$ is isomorphic to the sheaf of flat sections
$\beta_i$ such that $(\beta_i|\beta_i)=2$ via
\beq\label{Li-iso}
\beta_i\mapsto \frac{1}{2\sqrt{2}}\, \int_{\beta_i} dt_i
= \pm (\lambda-u_i)^{1/2}.
\eeq
The local EO recursion produces a set of multivalued {\em correlator forms}
\ben
\omega_{g,n}^{\alpha_1,\dots,\alpha_n}(u;\lambda_1,\dots,\lambda_n),\quad
\alpha_i\in \cL_{m_i},\quad (u,\lambda_i)\in \Delta_{m_i}^*,
\quad
1\leq m_i\leq N.
\een
Let us fix a base point in each
$\Delta_i^*$ ($1\leq i\leq N$). Then the values of the correlator forms depend on the choice of reference
paths -- one for each point $(u,\lambda_i)\in \Delta_{m_i}^*$. The recursion kernel is defined by 
\ben
K^{\beta_i,\beta_j}(u,\lambda_1,\lambda_2) = 
\frac{\frac{1}{2}\oint_{\lambda\in    C_{\lambda_2}}
P^{ij}(u,\lambda_1,\lambda)d\lambda
}{
P^j(u,\lambda_2)} 
\, \frac{d\lambda_1}{d\lambda_2},
\een
where $C_{\lambda_2}$ is a small loop around $u_j$ based at
$\lambda_2$. 
\begin{remark}\label{contour-antider}
The branch of the integrand in the contour integral $\oint_{
C_{\lambda_2}}$ is fixed in such a way that  the operation
$$
(\lambda_2-u_i)^{a}\mapsto 
\frac{1}{2}\oint_{\lambda\in C_{\lambda_2}}
(\lambda-u_i)^{a} = (\lambda_2-u_i)^{a+1}/(a+1)
$$ 
computes the anti-derivative (or primitive). In other words, using the
reference path to $(u,\lambda_2)$ and the contour $C_{\lambda_2}$ we
first 
specify the branch of $P^{ij}(u,\lambda_1,\lambda)$ by moving
continuously $\lambda$ along $C_{\lambda_2}$ and then we integrate the
resulting function backwards, i.e. the orientation of the cycle
$C_{\lambda_2}$ used in the contour integral is the opposite to the
orientation used to specify the branch of the integrand. 
\end{remark}
The local EO recursion takes the form
\ben
\omega_{0,2}^{\beta_i,\beta_j}(u,\lambda_1,\lambda_2) := 
\begin{cases}
P^{ij}(u,\lambda_1,\lambda_2) d\lambda_1d\lambda_2 & \mbox{ if } (\lambda_1,i)\neq
(\lambda_2,j),\\
P^{ii}_0(u,\lambda_1)d\lambda_1\, d\lambda_1 & \mbox{ otherwise } ,
\end{cases}
\een
and
\ben
&&
\omega_{g,n+1}^{\alpha_0,\dots,\alpha_n}(\lambda_0,\dots,\lambda_n) =
\\
&&
\sum_{i=1}^N 
\operatorname{res}_{\lambda=u_i}
K^{\alpha_0,\beta_i}(\lambda_0,\lambda) \Big( 
\omega_{g-1,n+2}^{\beta_i,-\beta_i,\alpha_1,\dots,\alpha_n}
(\lambda,\lambda,\lambda_1,\dots,\lambda_n) + \\
&&
\sum_{g'+g''=g}\sum_{i_1',\dots,i_{n'}'} 
\omega_{g',n'+1}^{\beta_i,\alpha_{i_1'},\dots,\alpha_{i'_{n'}}}
(\lambda,\lambda_{i_1'},\dots,\lambda_{i'_{n'}})
\omega_{g'',n''+1}^{-\beta_i,\alpha_{i_1''},\dots,\alpha_{i''_{n''}}}
(\lambda,\lambda_{i_1''},\dots,\lambda_{i''_{n''}})
\Big),
\een
where $\beta_i\in \cL_i$, the second and the third sums are over
all splitting $g'+g''=g$ and all subsets
$\{i_1',\dots,i_{n'}'\}\subset \{1,2,\dots,n\}$, and 
\ben
\{i_1'',\dots,i_{n''}''\} := \{1,2,\dots,n\}-\{i_1',\dots,i_{n'}'\}.
\een 

\begin{example}\label{ex:Frob-local-EO}
Let us assume that $U$ is equipped with a semi-simple Frobenius
structure such that $u_1,\dots,u_N$ are canonical
coordinates. Although, our definition of the periods
$I^{(-n)}_{\gamma^{(i)}}(u,\lambda)$ was given only for Frobenius structures
  corresponding to polynomial primitive forms, we define the period
  $I^{(-n)}_{\gamma^{(i)}}(u,\lambda)$ ($1\leq i\leq N$, $n\in \ZZ$)  in general
  as the unique solution to the second
  structure connection that has an expansion of the type
  \eqref{period-exp} for all $(u,\lambda)\in \Delta_i$. Note that in
  general we think of $\gamma^{(i)}$ us a section of the local system
  $\cL_i$, i.e., a choice of a holomorphic branch of
  $(\lambda-u_i)^{1/2}$.   Moreover, the recursion that we used to
  define the total ancestor potential still makes sense. Finally,
  substituting the expansion \eqref{period-exp} in the differential
  equations for the second structure connection yields a recursion
  that uniquely determines the matrices $R_k(u)$, $k\geq 0$ starting
  with $R_0(u)=1$. The main result of \cite{Mi1} is that the
  correlator forms
\ben
{}^{\rm
  Frob}\omega_{g,n}^{\alpha_1,\dots,\alpha_n}(u;\lambda_1,\dots,\lambda_n):= 
\langle
\phi_{\alpha_1}^+(u,\lambda_1;\overline{\psi}),\dots,\phi_{\alpha_n}^+(u,\lambda_n;\overline{\psi})\rangle_{g,n}
d\lambda_1\cdots d\lambda_n
\een
where $(u,\lambda_i)\in \Delta_{m_i}$, $\alpha_i\in \cL_{m_i}$, and
the insertions $\phi_{\alpha_i}$ are defined by \eqref{bosonic-field},
satisfy the local EO recursion with 
\ben
P^j(u,\lambda) = 4(I^{(-1)}_{\gamma^{(j)}}(u,\lambda),1)
\een 
and 
\ben
P^{ij}_{k,\ell}(u) = \frac{2^{k+\ell+1} V_{k,\ell}^{ij} (u)}{(2k-1)!! (2\ell-1)!!},
\een
where $V_{k,\ell}^{ij}(u)$ is the $(i,j)$-entry of the matrix
$V_{k,\ell}$ defined by 
\ben
\frac{R(u,z_1)^T R(u,-z_2)-1}{z_1+z_2} = \sum_{k,\ell=0}^\infty
V_{k,\ell}(u) (-z_1)^k(-z_2)^\ell. 
\een
\end{example}

\begin{example}\label{ex:EO-local-EO}
Let $\phi$ be a relative meromorphic differential on
$\overline{X}_U:=\overline{\pi}^{-1}(U)$  such that $\phi|_{X_U^{\rm
    loc}}$ is holomorphic and $\phi(p)\neq 0$ for all $p\in C\cap X_U$. 
Let $\omega_{g,n}(u;q_1,\dots,q_n)$ be the correlator forms defined by 
the following EO recursion (see \cite{DNOPS2}): all unstable
correlators (i.e. $2g-2+n\leq 0$) are $0$ except for 
\ben
\omega_{0,2}(u;p,q):=B_u(p,q),\quad p\neq q,
\een
where $B_u$ is the fundamental bi-differential on $\overline{X}_u$ and 
\ben
&&
\omega_{g,n+1}(u;q_0,q_1,\dots,q_n) = \\
&&
\sum_{i=1}^N
\frac{1}{2}\operatorname{res}_{p=p_i} 
\frac{\int_p^{\tau_i(p)} B_u(q_0,p') }{ dF(p)\, \int_p^{\tau_i(p)}
  \phi(p')} 
\Big( \omega_{g-1,n+1}(u;p,\tau_i(p),q_1,\dots,q_n)
+\\
&&
\sum_{g'+g''=g}\sum_{i_1',\dots,i_{n'}' }
\omega_{g',n'+1}(u;p,q_{i_1'},\dots,q_{i_{n'}'} )
\omega_{g'',n''+1}(u;\tau_i(p),q_{i_1''},\dots,q_{i_{n''}''} )\Big),
\een
where $\tau_i$ is the local involution defined via the local
coordinate $t_i(p)$ in a neighborhood of
the ramification point $p_i$ as $t_i(\tau_i(p)) := -t_i(p)$, the
last sum is over all subsets $\{i_1',\dots,i_{n'}'\}\subset
\{1,2,\dots,n\}$ and 
\ben
\{i_1'',\dots,i_{n''}''\}:= \{1,2,\dots,n\}\setminus \{i_1',\dots,i_{n'}'\}.
\een
Note that the above recursion is more general then the recursion
proposed originally in \cite{EO}. Namely, if $\phi(p)=-dy(p)$ for
some meromorphic function $y$ on $\overline{X}_u$, then the above
recursion coincides with the EO recursion for the spectral data
$(\overline{X}_u,x_u,y_u)$, where $x_u=F|_{X_u}$ and $y_u=y|_{X_u} $.
 
Let us define the following set of multi-valued symmetric
differentials 
\ben
{}^{\rm EO}\omega_{g,n}^{\alpha^1,\dots,\alpha^n}(u;\lambda_1,\dots,\lambda_n) :=
d_{\lambda_1}\cdots d_{\lambda_n} 
\int_{q_1\in \alpha_{u,\lambda_1}^1} \cdots
\int_{q_n\in \alpha_{u,\lambda_n}^n} 
\omega_{g,n}(u;q_1,\dots,q_n),
\een
where $\alpha^i\in \mathfrak{h}$, $1\leq i\leq n$ are given relative
cycles and the pair $(g,n)$ is stable, i.e., $2g-2+n>0$. Here  $\alpha^i_{u,\lambda}\in
H_1(X_u,X_{u,\lambda};\CC)$ is the relative cycle obtained from
$\alpha^i$ via a parallel transport along some reference path in
$(U\times \CC)'$ from $(u^\circ,\lambda^\circ)$ to $(u,\lambda)$. 
Suppose that $(u,\lambda_i)\in \Delta_{m_i}$ and that the reference
path from $(u^\circ,\lambda^\circ)$ to $(u,\lambda_i)$ and the cycle
$\alpha^i$ are such that
the relative cycle $\alpha^i_{u,\lambda} $ is  the vanishing cycle,
i.e., it can be represented by an arc $\beta^{m_i}_{u,\lambda} \subset
\Delta_{m_i}^*$. In
particular, $\alpha^i_{u,\lambda}\in \cL_{m_i}$ via the isomorphism \eqref{Li-iso}.
\begin{proposition}\label{EO:local-EO}
The correlator forms
${}^{EO}\omega_{g,n}^{\alpha^1,\dots,\alpha^n}(u;\lambda_1,\dots,\lambda_n)$
satisfy local EO recursion with 
\ben
P^{ij}(u,\lambda_1,\lambda_2) :=
d_{\lambda_1} \int_{p\in \beta^i_{u,\lambda_1}} 
d_{\lambda_2} \int_{q\in \beta^j_{u,\lambda_2}} 
B_u(p,q)
\een
and 
\ben
P^j(u,\lambda) = 4 \int_{q\in \beta^j_{u,\lambda}} \phi(q) .
\een
\end{proposition}
\proof
In order to avoid cumbersome notation we drop the superscript ``EO''
and we suppress the dependence on $u$ in the correlators, i.e., $u\in
U$ will be fixed throughout the proof and we put
\ben
\omega_{g,n}^{\alpha^1,\dots,\alpha^n}(\lambda_1,\dots,\lambda_n):=
{}^{EO}\omega_{g,n}^{\alpha^1,\dots,\alpha^n}(u;\lambda_1,\dots,\lambda_n).
\een
Let us apply to the EO recursion defining the forms $\omega_{g,n}$
the operations 
\beq\label{cycle-pullback}
d_{\lambda_i} \int_{q_i\in \alpha_i^{u,\lambda_i}},\quad 0\leq i\leq n.
\eeq
By definition the LHS becomes ${}^{\rm
  EO}\omega_{g,n}^{\alpha^0,\dots,\alpha^n}(\lambda_0,\dots,\lambda_n)$. 
On the other hand on the RHS we get some hybrid
correlators that have insertions both from $\overline{X}_u$ and
$\PP^1$, i.e., we have expressions of the  form
\ben
\omega_{g-1,n+2}^{\alpha^1,\dots,\alpha^n}(p,\tau_i(p),\lambda_1,\dots,\lambda_n)
\een
and 
\ben
\omega_{g',n'+1}^{\alpha^{i_1'},\dots,\alpha^{i_{n'}'}}
(p,\lambda_{i_1'},\dots,\lambda_{i_{n'}'}) 
\omega_{g'',n''+1}^{\alpha^{i_1''},\dots,\alpha^{i_{n''}''}}
(\tau_i(p),\lambda_{i_1''},\dots,\lambda_{i_{n''}''}) ,
\een
where each insertion $\lambda_i$ is paired with a corresponding cycle
$\alpha^i$ (appearing in the superscript) and the pair
$(\alpha^i,\lambda_i)$ means that we have applied to the corresponding
EO-correlator
form the operation \eqref{cycle-pullback}. However, note that 
\ben
\omega^\alpha_{g-1,n+2}(p,\tau_i(p),\dots) = \frac{1}{4}
\omega_{g-1,n+2}^{\beta_i,-\beta_i,\alpha}(\lambda,\lambda,\dots ) +\cdots,
\een
where $\alpha=(\alpha^1,\dots,\alpha^n)$ and 
\ben
\omega^{\alpha'}_{g',n'+1}(p,\dots)\omega^{\alpha''}_{g'',n''+1}(\tau_i(p),\dots) =
\frac{1}{4}
\omega^{\beta_i,\alpha'}_{g',n'+1}(\lambda,\dots)\omega^{-\beta_i,\alpha''}_{g'',n''+1}(\lambda,\dots)
+ \dots
\een
where $\lambda=F(p)$,
$\alpha'=(\alpha_{i'_1},\dots,\alpha_{i'_{n'}})$, 
$\alpha''=(\alpha_{i''_1},\dots,\alpha_{i''_{n''}})$
and the dots that follow the plus sign on the RHS stand for terms holomorphic at
$p=p_i$ (here one has to prove by induction on $(g,n)$ that
$\omega_{g,n}(p,\dots)+\omega_{g,n}(\tau_i(p),\dots)$ is analytic at
$p=p_i$). Using also that 
\ben
\frac{1}{2}\operatorname{res}_{p=p_i} = \operatorname{res}_{\lambda=u_i}
\een 
we get that the correlator forms
$\omega_{g,n}^{\alpha^1,\dots,\alpha^n}$ satisfy a recursion that has
the same form as local EO recursion with 
recursion kernel
\ben
K^{\beta^i,\beta^j}(u,\lambda_1,\lambda_2) = 
\frac{
d_{\lambda_1}\int_{p\in \beta^i_{u,\lambda_1} }
\int_{q\in \beta^j_{u,\lambda_2}}
B_u(p,q) 
}{
 4\int_{q\in\beta^j_{u,\lambda_2}} \phi(q) d\lambda_2},
\een
except that the initial conditions are given by 
\ben
\omega_{0,2}^{\beta^i,\beta^j}( \lambda_1,\lambda_2) = 
d_{\lambda_1}d_{\lambda_2} \int_{p\in \beta^i_{u,\lambda_1} }
\int_{q\in \beta^j_{u,\lambda_2}}
B(p,q) 
\een
for $(i,\lambda_1)\neq (j,\lambda_2)$ and 
\ben
\widetilde{\omega}_{0,2}^{\beta^i,-\beta^i}( \lambda,\lambda) :=-
\widetilde{\omega}_{0,2}^{\beta^i,\beta^i}( \lambda,\lambda) = 4B_u(q,\tau_i(q)),
\een
where $q$ is sufficiently close to $p_i$ and $\lambda=F(q)$. 

By definition (see Remark \ref{contour-antider}) 
\ben
\frac{1}{2}\oint_{\lambda\in C_{\lambda_2}} d_\lambda \int_{q\in \beta^j_{u,\lambda}}
B_u(p,q) = \int_{q\in \beta^j_{u,\lambda}} B_u(p,q),
\een
Therefore,
$K^{\beta^i,\beta^j}(u,\lambda_1,\lambda_2)$ coincides with the
recursion kernel of the local recursion introduced in the statement of
the proposition. The only issue that we have to resolve is that the
initial condition $\omega_{0,2}^{\beta^i,\beta^i}( \lambda,\lambda)$
of the local  EO recursion is supposed to be $P_0^{ii}(u,\lambda)$. On
the other hand  the correlator
$\omega_{0,2}^{\beta^i,\beta^i}( \lambda,\lambda)$  contributes to the recursion only
when we evaluate 
\ben
\omega_{1,1}^\alpha(u;\lambda)=
-\sum_{i=1}^N \operatorname{res}_{\mu=u_i} 
K^{\alpha,\beta^i}(\lambda,\mu)\omega_{0,2}^{\beta^i,\beta^i}(\mu,\mu).
\een
Since $K^{\alpha,\beta_i}(\lambda,\mu)d\mu$ is analytic at $\mu=u_i$
it is sufficient to prove that the difference 
\ben
\omega_{0,2}^{\beta^i,\beta^i}( \lambda,\lambda) -
\widetilde{\omega}_{0,2}^{\beta^i,\beta^i}( \lambda,\lambda) 
\een
is analytic at $\lambda=u_i$. This is a straightforward local
computation. Indeed, using that for $p,q\in
U\times \widetilde{\Delta}_i$ we have
\ben
B(p,q) = dt_i(p)dt_i(q)\Big(\frac{1}{(t_i(p)-t_i(q))^2}
+\sum_{m,n=0}^\infty B^{ii}_{m,n}(u) t_i(p)^m t_i(q)^n\Big)
\een
we get that 
\ben
P^{ii}_0(u,\lambda) = \frac{1}{4}(\lambda-u_i)^{-2} +
\sum_{m,n=0}^\infty 2^{m+n+1} B_{2m,2n}^{ii}(\lambda-u_i)^{m+n-1}
\een
and
\ben
\widetilde{\omega}_{0,2}^{\beta^i,\beta^i}( \lambda,\lambda) = 
d\lambda\cdot d\lambda \Big(
\frac{1}{4}(\lambda-u_i)^{-2} + \frac{2 B_{0,0}^{ii}}{\lambda-u_i}+\cdots\Big), 
\een
where the dots stand for terms analytic at $\lambda=u_i$. The
analyticity that we wanted to prove follows.
\qed
\end{example}

\begin{definition}\label{EO-Frob}
We say that a semi-simple Frobenius structure on $U$ is a solution to an EO
recursion (defined by a relative meromorphic differential on
$\overline{X}_U$) if the correlator forms defined by the corresponding
local EO recursions coincide, i.e., 
\ben
{}^{\rm
  Frob}\omega_{g,n}^{\alpha_1,\dots,\alpha_n}(u;\lambda_1,\dots,\lambda_n)
=
{}^{\rm
  EO}\omega_{g,n}^{\alpha_1,\dots,\alpha_n}(u;\lambda_1,\dots,\lambda_n),\quad
\forall u\in U.
\een
\end{definition}

\subsection{Frobenius manifolds and EO recursion}

Let $B_{m,n}^{i,j}$, $1\leq i,j\leq N$, $m,n\geq 0$ be the
coefficients defined via the expansion of the fundamental
bi-differential  
\beq\label{B-series}
B(p,q) = dt_i(p) dt_j(q) \Big(
\frac{\delta_{ij} }{(t_i(p)-t_j(q))^2} + \sum_{m,n=0}^\infty B_{m,n}^{i,j}
t_i(p)^m t_j(q)^n\Big),
\eeq
where the point $(p,q)$ is sufficiently close to $(p_i,p_j)\in
X_u\times X_u$. Let us denote by $B_{m,n}$ the $N\times N$ matrix
whose $(i,j)$-entry is $B_{m,n}^{i,j}$. The key result in comparing
the topological recursion with the Dubrovin's Frobenius structure is
the following identity.
\begin{lemma}\label{V=B}
The following identity holds
\ben
\frac{R_\Sigma(u,z_1)^T \, R_\Sigma(u,z_2)-1}{z_1+z_2} =
\sum_{m,n=0}^\infty B_{2m,2n} (2m-1)!! (2n-1)!!\, (-z_1)^m (-z_2)^n. 
\een
\end{lemma}
\noindent
This result is due to Eynard \cite{Ey} (see also Lemma 5.1 in
\cite{DNOPS2}). It can be proved using the same technique that we
used in the proof of Proposition \ref{prop:gb}.

\medskip
{\em Proof of Theorem \ref{t3}.} 
Let us assume first that $\omega=[\phi]\in \H(U)$ where $\phi$ is a sum of
homogeneous primary differentials of the same degree. According to
Corollary \ref{cor:z-indep} $\omega$ is a primitive form and the matrix
$R_\omega(u,z)=1$ (see the notation in Section
\ref{sec:spa}). Recalling Proposition \ref{prop:R}, we get that the
R-matrix of the Frobenius structure on $U$ corresponding to the
primitive form $\omega$ is $R(u,z)=R_\Sigma(u,z)$. We claim that both
sets of correlator forms ${}^{\rm Frob}\omega_{g,n}$ and ${}^{\rm
  EO}\omega_{g,n}$ are defined by the same local EO recursion. Indeed,
on the Frobenius side we have (see Example \ref{ex:Frob-local-EO})
\ben
{}^{\rm Frob} P^{ij}_{k,\ell}(u) = \frac{2^{k+\ell+1} V_{k,\ell}^{ij}(u)}{(2k-1)!! (2\ell-1)!!}
\een
and 
\ben
{}^{\rm Frob} P^j(u,\lambda) = 4(I^{(-1)}_{\beta^j}(u,\lambda),1).
\een
While on the EO-recursion side we have (see Example \ref{ex:EO-local-EO})
\ben
{}^{\rm EO}P^{ij}_{k,\ell}(u) = 2^{k+\ell+1} \, B^{ij}_{2k,2\ell}(u)
\een
and 
\ben
{}^{\rm EO} P^j(u,\lambda)=4\int_{q\in \beta^j_{u,\lambda}} \phi(q),
\een
where we used that 
\ben
&&
d_{\lambda_1}d_{\lambda_2}\, \int_{p\in \beta_i^{u,\lambda_1} }
\int_{q\in \beta_j^{u,\lambda_2}}
B(p,q)  = d\lambda_1 d\lambda_2\, \Big( 
\frac{\delta_{ij}}{(\lambda_1-\lambda_2)^2} \Big(
\frac{\sqrt{\lambda_2-u_i}}{\sqrt{\lambda_1-u_i}} +  
\frac{\sqrt{\lambda_1-u_i}}{\sqrt{\lambda_2-u_i}}
\Big)+\\
&&
+
\sum_{m,n=0}^\infty B_{2m,2n}^{i,j} 2^{m+n+1} 
(\lambda_1-u_i)^{m-1/2} \, (\lambda_2-u_i)^{n-1/2}
\Big),
\een
Using Eynard's identity Lemma \ref{V=B} we get $ {}^{\rm
  Frob}P^{ij}_{k,\ell}(u)={}^{\rm EO}P^{ij}_{k,\ell}(u)$. Recalling
the definition of the periods and using that the primitive form is
$[\phi]$ we get
\ben
{}^{\rm Frob}P^{j}(u,\lambda) = 4(I^{(-1)}_{\beta_j}(u,\lambda),1) =
4\int_{q\in \beta_j^{u,\lambda}} \phi(q) =
{}^{\rm EO}P^{j}(u,\lambda).
\een

In the opposite direction. Let us assume that there is a local EO
recursion defined by $P^{ij}(u,\lambda)$ ($1\leq i,j\leq N$) and
$P^j(u,\lambda)$ ($1\leq j\leq N$) such that the corresponding
correlator forms coincide with ${}^{\rm Frob}\omega_{g,n}$ and with
${}^{\rm EO}\omega_{g,n}$ for some semi-simple Frobenius structure on
$U$ and some EO-recursion defined by a relative meromorphic form
$\phi$ on $\overline{X}_U$. It is sufficient to prove that 
$R(u,z)=R_\Sigma(u,z)$. Indeed, let $\omega$ be the primitive form
corresponding to the Frobenius structure on $U$. According to
Proposition \ref{prop:R} if $R(u,z)=R_\Sigma(u,z)$, then
$R_\omega(u,z)=1$. Therefore, according to Corollary
\ref{cor:z-indep} the primitive form is represented by a linear
combination of homogeneous primary differentials. 

Both $R$-matrices can be extracted from the 3-point genus-0
correlators. Let us denote by $P^{ij,j}_{k}(u,\lambda)$ ($1\leq
i,j\leq N$, $k\geq 0$) the coefficients that appear in the expansion 
\ben
P^{ij}(u,\lambda_1,\lambda_2) = \sum_{k=0}^\infty
P^{ij,j}_k(u,\lambda_1) (\lambda_2-u_j)^{k-1/2}.
\een 
In particular, 
\ben
P^{ij,j}_0(u,\lambda)= \delta_{i,j}
(\lambda-u_i)^{-1-1/2}+\sum_{k=0}^\infty P^{ij}_{k,0}(u) (\lambda-u_i)^{k-1/2}.
\een
Using the local EO recursion we get
\ben
\omega_{0,3}^{\alpha_1,\alpha_2,\alpha_3}(u;\lambda_1,\lambda_2,\lambda_3)
= -4\sum_{j=1}^N 
P^{i_1j,j}_0(u,\lambda_1) P^{i_2j,j}_0(u,\lambda_2) P^{i_3j,j}_0(u,\lambda_3)
\frac{d\lambda_1\,d\lambda_2\,d\lambda_3}{P_0^j(u)},
\een
where $\alpha_a\in \cL_{i_a}$ ($a=1,2,3$). Suppose first that
$i_1=i_2=i_3=:i$. Then the coefficient in front of  
$$
-4(\lambda_1-u_i)^{-1/2}(\lambda_2-u_i)^{-1/2}(\lambda_3-u_i)^{-1/2}
d\lambda_1d\lambda_2d\lambda_3
$$
in the above correlator is $1/P_0^i(u)$. Therefore, ${}^{\rm Frob}P_0^i(u) ={}^{\rm
  EO}P_0^i(u) $ for all $i=1,2,\dots,N$. 
Suppose now that $i_2=i_3=:i$. Then the coefficient in front of 
$$
-4(\lambda_2-u_i)^{-1/2}(\lambda_3-u_i)^{-1/2} 
d\lambda_1\, d\lambda_2\, d\lambda_3
$$ 
is $P_0^{i_1i,i}(u,\lambda_1) /P_0^i(u)$. Therefore
\ben
{}^{\rm Frob}P_0^{i_1i,i}(u,\lambda_1) /{}^{\rm Frob} P_0^i(u) =
{}^{\rm EO}P_0^{i_1i,i}(u,\lambda_1) /{}^{\rm EO} P_0^i(u) .
\een 
Since we already proved that ${}^{\rm Frob} P_0^i(u) ={}^{\rm EO}
P_0^i(u)$ we get ${}^{\rm Frob}P_0^{i_1i,i}(u,\lambda_1) ={}^{\rm
  EO}P_0^{i_1i,i}(u,\lambda_1) .$ Comparing the coefficients in front
of $(\lambda_1-u_i)^{k-1/2}$ we get 
\ben
{}^{\rm Frob}P^{i_1i}_{k,0}={}^{\rm EO} P^{i_1i}_{k,0},
\een
i.e., $V^{ij}_{k,0}= B^{ij}_{2k,0} (2k-1)!!$. Finally we get
\ben
(R(u,z)^T-1)/z = \sum_{k=0}^\infty V^{ij}_{k,0} (-z)^k = 
\sum_{k=0}^\infty B^{ij}_{2k,0} (2k-1)!! (-z)^k = (R_\Sigma(u,z)^T-1)/z,
\een
where the last identity follows from Eynard's identity (Lemma
\ref{V=B}).
\qed

\section{Topological recursion for twisted de Rham cohomology}
\label{sec:desc}

Semi-simple Frobenius manifolds admit also the notion of
descendants. Using the forms defined by the topological recursion we will express the
descendant correlation functions for Hurwitz Frobenius manifolds as oscillatory integrals. Motivated
by this result, we propose to think of the forms defined by the
topological recursion as twisted de Rham cohomology classes. This
point of view allows us to find a generalization of the EO recursion
corresponding to the Frobenius manifolds defined by polynomial
primitive forms. 

\subsection{Calibration of the Frobenius structure}

Near $z=\infty$ the system \eqref{qde:i}--\eqref{qde:z} admits a {\em
  weak Levelt solution}, i.e., a fundamental solution of the form
\ben
\Phi(t,z) = S(t,z) z^\delta z^\nu, 
\een
where the matrices $S$, $\delta$, and $\nu$ have the following
properties. We have an expansion
$S(t,z)=S_0+S_1(t)z^{-1}+S_2(t)z^{-2}\cdots$ with $S_0$ constant
(independent of $t$ and $z$) {\em invertible} matrix. The matrices
$\delta=\operatorname{Diag}(\delta_1,\dots,\delta_s)$ and 
$\nu=\operatorname{Diag}(\nu_1,\dots,\nu_s)$ are block-diagonal with
blocks $\delta_i$ and $\nu_i$ constant matrices of the same size. 
The block $\nu_i$ is an upper-triangular nilpotent
matrix and the  block $\delta_i$ is a diagonal matrix
$$
\delta_i=\operatorname{Diag}(m_{i,1}+r_i,\dots,m_{i,p_i}+r_i)
$$ 
such that
$-1<\operatorname{Re} \,  r_i \leq 0$ and
$m_{i,1}\leq m_{i,2}\leq \cdots \leq m_{i,p_i}$ is an increasing sequence
of integers.  The matrix representation of the system
\eqref{qde:i}--\eqref{qde:z}  depends on the choice of a flat
coordinate system. Choosing a different flat coordinate system
transforms \eqref{qde:i}--\eqref{qde:z} via a constant gauge
transformation. Therefore without lost of generality we may assume
that  $S_0=1$. 

It is convenient to introduce the following notation. Let
$\operatorname{spec}(\delta)$ be the set of eigenvalues of the
operator
\ben
\operatorname{ad}_\delta:\mathfrak{gl}(H) \to \mathfrak{gl}(H),\quad
X\mapsto [\delta,X].
\een
Let us denote by $\mathfrak{gl}_a(H)$ the eigensubspace of
$\operatorname{ad}_\delta$ with eigenvalue $a$. Then we have a direct
sum decomposition of vector spaces
\ben
\mathfrak{gl}(H)=\bigoplus_{a\in \operatorname{spec}(\delta)} \mathfrak{gl}_a(H).
\een
Let us denote by $X_{[a]}$ the projection of $X$ on
$\mathfrak{gl}_a(H)$. The matrices $S$, $\delta$, and $\nu$ are
identified with elements of $\mathfrak{gl}(H)$ via the basis
$\{\phi_i\}_{i=1}^N\subset H$ that we fixed above. 

Substituting the fundamental series $\Phi(t,z)$ in \eqref{qde:z} and
comparing the coefficients in front of powers of $z$ we get that 
\ben
\theta = \delta+\nu_{[0]},\quad 
kS_k +[\theta,S_k] =E\bullet S_{k-1} +\sum_{\ell=1}^k S_{k-\ell}
\nu_{[-\ell]},\quad k>0.
\een
In particular, a weak Levelt solution can be constructed by setting
$(S_k)_{[-k]}=0$ and solving recursively for $\nu_{[-k]}$ and
$(S_k)_{[a]}$ for all $k>0$ and all $a\neq -k$.  This is a standard
procedure, so we skip the details. 

\begin{proposition}\label{prop:sympl-cali}
There exists a weak Levelt solution such that 
\ben
S(t,-z)^T\ S(t,z)=1,
\een
where ${}^T$ is transposition with respect to the Frobenius pairing on
$H=T_{u^\circ}U$. 
\end{proposition}
Proposition \ref{prop:sympl-cali} is known if $\theta$ is
diagonalizable (see \cite{Du2}). In fact, the polynomiality
of the primitive form might be sufficient to prove that $\theta$ is
daigonalizable. However, at this point this is unknown. Let us modify
the argument from \cite{Du2} in order to cover the case of $\theta$
non-diagonalizable.  
\begin{lemma}
The eigenvalues of $\theta$ are rational numbers.
\end{lemma}
\proof
The eigenvalues of $e^{2\pi\sqrt{-1}\theta}$ coincide with the eigenvalues of the
monodromy matrix of the fundamental solution $\Phi(t,z)$, because
$\theta=\delta+\nu_{[0]}$.  Therefore,
it is enough to prove that the monodromy matrix of 
$\Phi(t,z)$ has eigenvalues that are roots of 1.  
On the other hand the system \eqref{qde:i}--\eqref{qde:z} can be
solved in terms of oscillatory integrals. Therefore, there exists a
linear isomorphism 
\ben
\Pi: H_1(X_{u^\circ},\operatorname{Re}(F/z)\ll 0,\CC)\to H
\een
such that 
\ben
J_{\Gamma}(t,z)=\Phi(t,z)\Pi(\Gamma) ,\quad \forall \Gamma.
\een
Therefore, the eigenvalues of the monodromy matrix are the same as the
eigenvalues of the monodromy operator 
$$
M: H_1(X_{u^\circ},\operatorname{Re}(F/z)\ll 0,\CC)\to
H_1(X_{u^\circ},\operatorname{Re}(F/z)\ll 0,\CC) 
$$ 
corresponding to the parallel transport  around $z=\infty$ with respect
to the Gauss--Manin connection. Note that the subspace
$$
\operatorname{Im}(M-1) \subset
H_1(X_{u^\circ},\operatorname{Re}(F/z)\ll 0,\CC)
$$ 
admits a basis consisting of
compact cycles and cycles supported in a neighborhood of a puncture
$\infty_i$. If $\Gamma$ is a compact cycle then the corresponding
integral is analytic near $z=\infty$, so $(M-1)\Gamma=0$. If $\Gamma$
is supported near the puncture $\infty_i$, then we may assume that
$X_{u^\circ}=\CC$ and $F(p)=p^{m_i}$ and an easy computation yields
that the integral depends analytically on $z^{-1/m_i}$, so
$(M^{m_i}-1) \Gamma=0$. We get that $M$ satisfies the following equation:
\ben
(M-1)^2\prod_{i=1}^d (M^{m_i}-1) = 0.
\een
The eigenvalues of $M$ must be also solutions of the above equation,
so they are roots of unity as claimed.
\qed

Let $G_\delta$ be the subgroup of $\operatorname{GL}(H)$ consisting of
linear operators $C$ such that 
\ben
C=1+\sum_{\ell=1}^\infty C_{[-\ell]}.
\een
The group $G_\delta$ acts on the set of weak Levelt solutions
\ben
\Phi(t,z)\mapsto \Phi(t,z) C=:\widetilde{S}(t,z) z^\delta z^{\widetilde{\nu}},
\een
where 
\ben
\widetilde{S}(t,z)= {S}(t,z) z^{\delta} C z^{-\delta},\quad
\widetilde{\nu} = C^{-1} \nu C.
\een
Note that $\nu$ belongs to the Lie algebra $\mathfrak{g}_\delta$ of
$G_\delta$ consisting of matrices $x$ such that $1+x\in G_\delta$. In
particular, $e^{2\pi\sqrt{-1}\delta}$ commutes with $\nu$. 
Proposition \ref{prop:sympl-cali} is a corollary of the following
lemma.
\begin{lemma}
Let $\Phi(t,z)=S(t,z)z^\delta z^\nu$ be a weak Levelt solution. There
exists a unique $C\in G_\delta$ such that $C_{[\ell]}^T = (-1)^\ell
C_{[\ell]}$ for all $\ell\leq 0$ and 
$\widetilde{S}(t,z):= {S}(t,z) z^{\delta} C z^{-\delta}$ satisfies the
symplectic condition $\widetilde{S}(t,-z)^T \widetilde{S}(t,z)=1$. 
\end{lemma}
\proof
Let us first point out that the projection $X\mapsto X_{[a]}$ commutes
with transposition, i.e., $(X_{[a]})^T = (X^T)_{[a]}$ for all $X\in
\operatorname{GL}(H)$ and $a\in \operatorname{spec}(\delta)$. This
follows from the skew-symmetry of $\theta=\delta+\nu_{[0]}$. Namely,
using that there exists an integer $m$ for which
$e^{2\pi\sqrt{-1}\theta m}$ is a unipotent operator (recall that
$\theta$ has rational eigenvalues) and that
$\theta+\theta^T=0$, we get $\nu_{[0]}+(\nu_{[0]})^T=0$ and
$\delta+\delta^T = 0$. The latter implies $[\delta,X]^T=[\delta,X^T],$
so ${}^T:\mathfrak{gl}_a(H)\to \mathfrak{gl}_a(H)$ is a linear
isomorphism for all $a\in \operatorname{spec}(\delta)$.  Our claim
follows easily. In the rest of the proof we put $X_{[a]}^T:=
(X_{[a]})^T = (X^T)_{[a]}$.  

Let us fix $z$ to be a positive real number and define 
\ben
A:=\Phi(t,-z)^T \Phi(t,z),
\een
where we define $(-z)^x:=e^{x\log(-z)}$ with $\log(-z):=\log z
+\pi\sqrt{-1}.$ The differential equations
\eqref{qde:i}--\eqref{qde:z} imply that $A$ is a constant matrix
independent of $t$ and $z$. We will also make use of the following two
properties of $A$ 
\beq\label{A:eqn-1}
A^T = A e^{2\pi\sqrt{-1} \delta} e^{2\pi\sqrt{-1} \nu}.
\eeq
and 
\beq\label{A:eqn-2}
\nu^T = -A\, \nu\, A^{-1}.
\eeq
The first one is proved by 
transposing the identity defining $A$ and analytically
continuing in $z$ from $z$ to $-z$ along an arc in the upper
half-plane. To prove the second one, first  we pick an integer $m$
such that $e^{2\pi\sqrt{-1}m \delta}=1$.  Let us analytically continue
the identity defining $A$ along a loop that goes $m$ times around
$z=\infty$. We get 
\ben
e^{2\pi\sqrt{-1} m \nu^T} A e^{2\pi\sqrt{-1} m \nu}= A\quad
\Rightarrow 
\quad
e^{2\pi\sqrt{-1} m \nu^T} = e^{2\pi\sqrt{-1} m (-A\nu A^{-1})}.
\een
It remains only to use that we can take the logarithm, because $\nu^T$
and $-A\nu A^{-1}$ are nilpotent matrices. 

Let us try to find $C\in G_\delta$, such that $\widetilde{S}(t,-z)^T
\widetilde{S}(t,z)=1.$ Using \eqref{A:eqn-2} we get 
\ben
S(t,-z)^T S(t,z) = (-z)^{\delta} Ae^{\pi\sqrt{-1} \nu} z^{-\delta}.
\een 
Comparing the coefficients in front of the powers of $z$, we get that 
\ben
B:=e^{\pi\sqrt{-1}\delta} A e^{\pi\sqrt{-1} \nu}
\een
belongs to $G_\delta$. Moreover, equation \eqref{A:eqn-1} implies that 
$B_{[\ell]}^T = (-1)^\ell B_{[\ell]}$ for all $\ell\leq 0$. 
After a straightforward computation we get the following equation for
$C$:
\ben
(e^{\pi\sqrt{-1}\delta} C^T e^{-\pi\sqrt{-1}\delta})\, B\, C=1.
\een
The projection of this equation onto $\mathfrak{gl}_{-\ell}(H)$
$(\ell>0)$ yields
\ben
B_{[-\ell]} + C_{[-\ell]} + (-1)^\ell C_{[-\ell]}^T
+\sum_{\substack{i+k+j=\ell\\ 0\leq i,k,j<\ell }}
(-1)^i C_{[-i]}^T B_{[-k]} C_{[-j]} = 0.
\een
The above equations determine a unique solution satisfying the
additional constraint $C_{[-\ell]}^T=(-1)^\ell C_{[-\ell]}$.
\qed

The operator series $S(t,z)$ satisfying the conditions of Proposition
\ref{prop:sympl-cali} are called {\em calibrations}, i.e., these are
operator series of the form $1+S_1(t)z^{-1}+\cdots$ such that
$S(t,-z)^TS(t,z)=1$ and the gauge transformation of the Dubrovin's
connection takes the form
\ben
S(t,z)^{-1} \, \nabla \, S(t,z) = d -\Big(\theta z^{-1} +
\sum_{\ell=1}^\infty \nu_{[-\ell]} z^{-\ell -1}\Big) dz
\een
for some $\nu\in G_\delta$ where $\delta:=\theta_{ss}$ is the
semi-simple part of $\theta$ in the Jordan--Chevalley decomposition of
$\theta$.

\subsection{The total descendant potential }

Let us fix a calibration $S(t,z)$ and view the coefficients $S_k$ as matrix-valued
functions on $U$. Recall that when we introduced the
flat coordinate system on $U$ we had the freedom to choose
$t_a^\circ=t_a(u^\circ)$ as we wish. Let us define
$t_a^\circ:=(S_1(u^\circ)1,\phi^a)$, $1\leq a\leq N$. Using the differential equations
\eqref{qde:i}--\eqref{qde:z} it is easy to check that 
\ben
t_a = (S_1(u)1,\phi^a),\quad 1\leq a\leq N.
\een 
Using the flat coordinates we embed our Frobenius manifold $U$ in $\CC^N$ as an open
neighborhood of $t^\circ=(t_1^\circ,\dots,t_N^\circ)$. We will identify points $u\in U$ with their
coordinates $t=(t_1,\dots,t_N)\in \CC^N$.  

Our goal is to define the following set of {\em
  descendant correlators}  
\ben
\langle \phi_{a_1}\psi^{k_1},\dots,\phi_{a_n}\psi^{k_n}\rangle_{g,n},
\een
where $g,n\geq 0$, $k_1,\dots,k_n\geq 0$, and $1\leq a_1,\dots,a_n\leq
N$ are arbitrary. By definition the descendant correlators are
analytic functions on the Frobenius manifold $U$. If $2g-2+n>0$, then
the above correlator is defined in terms of the calibration $S(t,z)$
and the ancestor correlators as follows 
\ben
\langle [S(t,\overline{\psi})
\phi_{a_1}\overline{\psi}^{k_1}]_+,\dots, 
[S(t,\overline{\psi})
\phi_{a_n}\overline{\psi}^{k_n}]_+
\rangle_{g,n},
\een
$[\ ]_+$ denotes the operation truncating the terms containing
negative powers of $\overline{\psi}$ and the ancestor correlator is evaluated at
a point $u\in U$ with coordinates $t=(t_1,\dots,t_N)$. 

If $2g-2+n\leq 0$, then there are 4 cases. If $(g,n)=(0,2)$, then 
\ben
\langle \phi_a\psi^k,\phi_b\psi^\ell\rangle_{0,2} := (W_{k\ell}\phi_b,\phi_a),
\een
where the linear operators $W_{k\ell}$ are defined by 
\ben
\sum_{k,\ell=0}^\infty W_{k\ell}z^{-k}w^{-\ell} = 
\frac{S(t,z)^T\, S(t,w)-1}{z^{-1}+w^{-1}}.
\een
If $(g,n)=(0,1)$, then 
\ben
\langle \phi_a\psi^k\rangle_{0,1} := 
\langle \phi_a\psi^{k+1},1\rangle_{0,2}= 
(S_{k+2}(t)\phi_a,1).
\een
If $(g,n)=(0,0)$, then 
\ben
\langle\ \rangle_{0,0}:= 
-\frac{1}{2}\langle\psi-t \rangle_{0,1} =
\frac{1}{2} ((S_2(t) S_1(t)-S_3(t))1,1),
\een
where we identified $\CC^N\cong H$ via $t\mapsto \sum_a t_a \phi_a$. 
Note that $F^{(0)}(t):=\langle\ \rangle_{0,0}$ is the primary
potential of the Frobenius structure. 
Finally if $(g,n)=(1,0)$, then 
\ben
\langle\ \rangle_{1,0}:= \frac{1}{2} \int_0^t \sum_{i=1}^N R_1^{ii}(u)
du_i -\frac{1}{24} \sum_{i=1}^N \log c_i(u,0),
\een
where $R_1(u)$ is the coefficient in front of $z$ in the matrix
$R(u,z)$ defined via the stationary phase asymptotic (see Section \ref{sec:spa}),
$R_1^{ij}(u)$ is the $(i,j)$ entry of $R_1(u)$, and
$c_i(u,z)=K(\omega,[\omega_i])$ (see Section \ref{sec:spa}). The
function $F^{(1)}(t):= \langle\ \rangle_{1,0}$ is also known as the
genus-1 primary potential. 

The {\em total descendant potential} is by definition the following
generating series for the descendant correlators
\ben
\mathcal{D}_t(\hbar,\mathbf{t}) =\exp \Big(\sum_{g,n=0}^\infty 
\frac{\hbar^{g-1}}{n!}
\langle\mathbf{t}(\psi),\dots,\mathbf{t}(\psi)\rangle_{g,n}\Big),
\een
where $\mathbf{t}(\psi)=\sum_{k=0}^\infty \sum_{a=1}^N t_{k,a} \phi_a
\psi^k$, the correlators are expanded multilinearly in the formal
variables $t_{k,a}$, and all descendant correlators are evaluated at
the point $t\in U$. 

\begin{remark}
The total descendant potential has the following translation
symmetry. Put $\mathcal{D}^\circ := \mathcal{D}_{t^\circ}$, then 
\ben
\mathcal{D}_t(\hbar,\mathbf{t}) = \mathcal{D}^\circ(\hbar,\mathbf{t}+t-t^\circ)
\een
for all $t$ sufficiently close to $t^\circ$.  
\end{remark}

\begin{remark}
There is an elegant way to write the relation between descendants and
ancestors using Givental's quantization formalism (see \cite{G2}). The
above definition although a bit cumbersome is more convenient for our
purposes. 
\end{remark}

\subsection{Topological recursion and descendants} 

Let $\omega$ be the primitive form corresponding to a primary
differential. We are going to express the
descendant correlators in terms of the forms $\omega_{g,n}$. 
Let us fix a reference point $(u^\circ,z^\circ)\in U\times \CC^*$,
where $\CC^*:=\CC\setminus{\{0\}},$ such that $z^\circ$ is a positive
real number such that the fundamental solution
$\Phi(u,z)=S(u,z)z^\delta z^\nu$ is convergent for all $(u,z)$
sufficiently close to $(u^\circ,z^\circ)$. There exists a unique linear isomorphism 
\ben
\Pi: H_1(X_{u^\circ},\operatorname{Re}(F/z^\circ)\ll 0;\CC)\to H
\een
such that 
\ben
J_\Gamma(u,z) = S(u,z) z^\delta z^\nu \ \Pi(\Gamma)
\een
for all $(u,z)\in U\times \CC^*$ sufficiently close to
$(u^\circ,z^\circ)$. Here the value of $z^x:=e^{x\,\log z}$ for $x=\delta,\nu$
is chosen such that $\log z^\circ = \ln z^\circ$.  

Let $\omega_{g,n}(q_1,\dots,q_n)$ be the correlator forms defined by
the Eynard--Orantin recursion. Note that if $\Gamma$ is a small loop
in $X_u$ around one of the ramification points $p_j(u)$, then 
$$
\oint_{q_i\in \Gamma} e^{F(q_i)/z} \omega_{g,n}(q_1,\dots,q_n) = 0. 
$$
Indeed, this is easy to check when $(g,n)=(0,2)$ or $n=1$, while for
the remaining cases we can argue by induction on $(g,n)$ using that
the correlator forms are symmetric. 

For a given set of semi-infinite cycles
$$
\Gamma_a\in H_1(X_{u^\circ},\operatorname{Re}(F/z^\circ)\ll 0;\CC),
\quad
1\leq a\leq n
$$ 
we define the following integrals
\beq\label{omega:lt}
\prod_{a=1}^n (-2\pi z_a)^{-1/2}\ 
\int_{q_1\in \Gamma_1}\cdots \int_{q_n\in \Gamma_n}
e^{\frac{F(q_1)}{z_1}+\cdots + \frac{F(q_n)}{z_n} }
\omega_{g,n}(q_1,\dots,q_n),
\eeq
where each cycle $\Gamma_i$ is represented by a path avoiding the
ramification points (where the integrand might have poles). According
to our previous remark, the integral is independent of the choice of
representative paths. 
\begin{proposition}\label{EO-desc}
If $2g-2+n>0$, then the integral \eqref{omega:lt} coincides with the
following descendant correlator
\ben
\Big\langle
\frac{z_1^\delta z_1^\nu \, \Pi(\Gamma_1)}{\psi-z_1},\dots,
\frac{z_n^\delta z_n^\nu \, \Pi(\Gamma_n)}{\psi-z_n}
\Big\rangle_{g,n}.
\een
\end{proposition}
\proof
Let us fix $u\in U$ and $z_1,\dots, z_n$. It is enough to prove the
formula when each cycle 
$\Gamma_a=\varprojlim \gamma^{(i_a)}_{\lambda}$ is the limit as
$\lambda\to \infty$ of the Lefshetz thimbles corresponding to a
critical value $u_{i_a}$, where $\lambda$ varies along a path from
$u_{i_a}$ to $\infty$ such that $\operatorname{Re}(\lambda/z_i)<0$
(see Section \ref{sec:periods}). For such cycles the integral
\eqref{omega:lt} takes the form
\ben
\prod_{a=1}^n (-2\pi z_a)^{-1/2}\ 
\int_{u_{i_1}}^\infty  \cdots \int_{u_{i_n}}^\infty 
e^{\frac{\lambda_1}{z_1}+\cdots + \frac{\lambda_n}{z_n} }
\omega^{\gamma^{(i_1)},\dots, \gamma^{(i_n)}}_{g,n}
(\lambda_1,\dots,\lambda_n).
\een
Recalling Theorem \ref{t3} we can write the above integral as
an $n$-pointed genus-$g$ ancestor correlator in which the insertion 
on the $a$th place is
\ben
\sum_{k=0}^\infty 
(-2\pi z_a)^{-1/2}\ \int_{u_{i_a}}^\infty e^{\lambda_a/z_a}
I^{(k+1)}_{\gamma^{(i_a)}}(u,\lambda_a) (-\overline{\psi})^k=
-J_{\Gamma_a}(u,z_a) \ \sum_{k=0}^\infty  
\overline{\psi}^k z_a^{-k-1}.
\een
The above expression can be written as 
\ben
-\sum_{k=0}^\infty 
[S(u,\overline{\psi})\overline{\psi}^k]_+ z_a^{-k-1}
\, 
z_a^\delta z_a^\nu \Pi(\Gamma_a). 
\een
To complete the proof it remains only to recall the definition of the
descendant correlators (in the stable range).
\qed

\subsection{Generalization}

Let $x:\Sigma\to \PP^1$ be a branched covering with ramification
profile the same as in the EO-recursion. Let us define the twisted de
Rham cohomology group $H_{\rm twdR}(\Sigma^n,x) $ as the following quotient
\ben
(\Omega^{1,\infty})^{\boxtimes n}
(\Sigma^n)[z_1^{\pm 1},\dots,z_n^{\pm 1}]/\Big(\sum_{i=1}^n z \  d_i + dx_i\wedge\Big)
(\Omega^{0,\infty})^{\boxtimes n} (\Sigma^n) [z_1^{\pm 1},\dots,z_n^{\pm 1}],
\een
where $d_i$ is the de Rham differential on the $i$th copy of $\Sigma$
in the direct product $\Sigma^n:=\Sigma\times\cdots \Sigma$, the
restriction of $dx_i$ on the $j$th slot of the direct product
$\Sigma^n$ is $dx$ if $j=i$ and $0$ if $j\neq i$. Using our results in
Section \ref{sec:per-iso} (for the case $n=1$) we get that $H_{\rm twdR}(\Sigma^n, x)$
is a free $\CC[z_1^{\pm 1},\dots,z_n^{\pm 1}]$-module of rank $n^N$ and a basis is given
by 
\ben
[\omega_{i_1,\dots,i_n}]:=[\omega_{i_1}]\boxtimes \cdots \boxtimes
[\omega_{i_n}], \quad 1\leq i_1,\dots,i_n\leq N, 
\een 
where $\omega_i$ is the good basis constructed in \ref{subsec:gb}.  

We would like to define a recursion that defines classes
\ben
\omega_{g,n}\in  H_{\rm twdR}(\Sigma^n,x),
\een
that are symmetric with respect to simultaneous permutations of 
$(q_1,\dots,q_n)\in \Sigma^n$ and $(z_1,\dots,z_n)$. 
The recursion depends on the choice of a differential $\omega\in
\Omega^{1,\infty}(\Sigma)[z]$  and a symmetric bi-differential 
$W \in (\Omega^{1,\infty})^{\boxtimes 2}[z,w]$ satisfying the
following conditions. The differential $\omega$ should be a
holomorphic volume form in a neighborhood of each finite ramification
point $p_i$, i.e., if we represent $\omega$ by $\sum_{n=0}^{n_0}
\omega^{(n)}(-z)^n$, then $\omega^{(0)}(p_i)\neq 0$ for all finite ramification
points $p_i$. The symmetric bi-differential is required to have the
form
\ben
W = B+\sum_{m,n\geq 0 } \sum_{i,j=1}^N V_{m,n}^{ij}
[\omega_i]\boxtimes[\omega_j] (-z_1)^m (-z_2)^n,
\een
where $B$ is the fundamental bi-differential (in the usual
EO-recursion) and $V:=\sum_{m,n} V_{m,n} (-z_1)^m(-z_2)^n$ is a polynomial
whose coefficients $V_{m,n}$ are square matrices of size $N$,
satisfying the symmetry condition $V_{mn}^{ij} = V_{nm}^{ji}$, where
for a matrix $A$ we denoted by $A^{ij}$ the entry in position
$(i,j)$.  

In order to define the recursion, we need to introduce the operator
\ben
\partial_x : \Omega^{1,\infty}(V)\to
\Omega^{1,\infty}(V),\quad \theta\mapsto  d\Big(\frac{\theta}{dx}\Big)
\een
defined for all open subsets $V\subset \Sigma$ that do not
contain finite ramification points. If $\theta$ has finite order poles
at some ramification point, then so does $\partial_x \theta$. We will
need also to work with an operator inverse to $\partial_x$
\ben
\partial_x^{-1} \theta (p) = dx(p)\wedge d_p^{-1} \theta(p) .
\een
The choice of $d_p^{-1} \theta(p)$ is not unique so the inverse
operation yields multi-valued analytic forms. However, note that if
$\theta\in \Omega^{1,\infty}(V)$, then the ambiguity in the definition
of $\partial_x^{-n}\theta\in \Omega^{1,\infty}(V)$  is up to $dg$,
where $g$ is a polynomial in $x$ of degree  at most $n$.  

Let us assume that 
\ben
\theta\in \Omega^{1,\infty}(\Sigma\setminus{\{p_1,\dots,p_N\}})
\een
has a finite order pole at every ramification point $p_i$. Then for
all $n\gg0$ the 1-form $\partial_x^{-n}\theta$ is analytic in a
neighborhood of all ramification points, so using the excision
principle (see Proposition \ref{per-iso}, Part b)) we can define a
twisted de Rham cohomology class $[\partial_x^{-n}\theta]\in
H_{\rm twdR}(\Sigma,x)$. It is easy to check that $[\partial_x^{-n}\theta]
$ is independent of the choice of a branch of $\partial_x^{-n}\theta$.
Using the relation
\ben
[\partial_x^{-n}\theta]= (-z)^{-1} [\partial_x^{-n-1}\theta] ,\quad
n\gg 0.
\een
we get that we have the following map 
\ben
\Omega^{1,\infty}(\Sigma\setminus{\{p_1,\dots,p_N\}})\to H_{\rm
  twdR}(\Sigma,x),\quad \theta\mapsto [\theta]:= (-z)^{-n}[\partial_x^{-n}\theta].
\een

We would like to generalize the EO-recursion as follows. For initial
condition put 
\ben
\omega_{0,2}(q_1,q_2) = W(q_1,q_2).
\een
The recursion is defined by the same formula as before except that we
replace the recursion kernel 
\ben
\frac{\int_p^{\tau_i(p)} B(q_0,p')}
{
dx(p)
\int_p^{\tau_i(p)} \phi(p')
}
\een
by 
\ben
\frac{\int_p^{\tau_i(p)} \Big(B(q_0,p')+\sum_{m,n\geq 0}\sum_{i,j=1}^N
V_{mn}^{ij} [\omega_i(q_0)](-z_0)^m (\partial_x^{-n}\omega_j)(p') \Big)}
{
dx(p)
\int_p^{\tau_i(p)}  \Big(\sum_{n=0}^{n_0} (\partial_x^{-n}\omega^{(n)})(p')\Big)
}.
\een

\subsection{Topological recursion for polynomial primitive forms}
The generalization of the EO-recursion proposed above allows us to
extend the statements of Theorem \ref{t3} and \ref{EO-desc}
to the case of polynomial primitive forms. 

To begin with let us assume that $\omega\in \H(U)$ is a polynomial
primitive form represented by $\sum_{n=0}^{n_0}
\omega^{(n)}(-z)^n$. Let us recall also the matrix $R(u,z)$ defining the
corresponding ancestor invariants (see Proposition \ref{prop:R}).  In
order to define a corresponding generalized  
topological recursion we have to specify a differential and a
bi-differential satisfying the conditions described in the previous
section. We choose the differential to be the primitive form. While
the bi-differential 
$$
W(q_1,q_2)=B(q_1,q_2)+\sum_{m,n=0}^\infty \sum_{i,j=1}^N V_{mn}^{ij}
\omega_i(q_1)\omega_j(q_2) (-z_1)^m (-z_2)^n
$$ 
is chosen in such a way that 
\ben
\frac{R(u,z_1)^T R(u,z_2)-1}{z_1+z_2} =\sum_{m,n=0}^\infty W_{2m,2n}
\, (2m-1)!!\, (2n-1)!!\, (-z_1)^m (-z_2)^n,
\een
where $W_{m,n}$ is the matrix whose $(i,j)$-entry $W_{mn}^{ij}$ is
defined as the coefficient in front of
$t_i(q_1)^mt_j(q_2)^ndt_i(q_1)dt_j(q_2)$ in  the Laurent series expansion of 
\ben
B(q_1,q_2)+\sum_{m,n=0}^\infty \sum_{i,j=1}^N V_{mn}^{ij}
(\partial_x ^{-m}\omega_i)(q_1)\, (\partial_x^{-n} \omega_j) (q_2)
\een
at the point $(p_i,p_j).$ 

Using the Taylor series expansion at $q=p_i$
\ben
\partial^{-m}_x\omega_{i'}(q) = dt_i\Big( -\frac{t_i^2m}{(2m-1)!!} + 
\sum_{k=0}^\infty B_{k0}^{ii'} \frac{t_i^{k+2(m+1)}}{(k+1)(k+3)\cdots
  (k+2m+1)}
\Big) 
\een
and Lemma \ref{V=B} we get that the matrices $V_{m,n}$ are uniquely
determined from the identity
\ben
\frac{R_{\omega}(u,z_1)R_\omega(u,z_2)^T-1}{z_1+z_2} =
\sum_{m,n=0}^\infty V_{m,n} (-z_1)^m (-z_2)^n.
\een
Note that since $R_\omega$ is polynomial in $z$ only finitely many
$V_{m,n}\neq 0$. 

Let 
$$
\omega_{g,n}=\sum_{\kappa=(k_1,\dots k_n)} \omega_{g,n;\kappa}
(-z_1)^{k_1}\cdots (-z_n)^{k_n},\quad 
\omega_{g,n;\kappa}\in (\Omega^{1,\infty})^{\boxtimes n}(\Sigma^n),
$$ 
be the forms defined by the generalized  topological recursion. The
multivalued correlator forms
$\omega_{g,n}^{\beta_1,\dots,\beta_n}(\lambda_1,\dots,\lambda_n)$ are
define by the same formulas as before, except that we identify
$\omega_{g,n}$ with 
\ben
\sum_{\kappa=(k_1,\dots k_n)} 
(\partial_{x_1})^{-k_1}\cdots (\partial_{x_n})^{-k_n}\omega_{g,n;\kappa}
\een
The arguments used in the proof of Theorem \ref{t3} and
\ref{EO-desc} can be repeated, so the conclusions of both propositions
hold.

\bibliographystyle{amsalpha}

\begin{thebibliography}{FKRW}

\bibitem{ACG}
E.~Arbarello, M.~Cornalba, and P.~ Griffiths.
\textit{Geometry of algebraic curves II.}
Grundlehren volume 268, Springer--Verlag, Berlin Heidelberg, 2011.

\bibitem{AGV}
V.I.~Arnol'd, S.M.~Gusein-Zade, and A.N.~Varchenko.
\textit{Singularities of differentiable maps. Vol. II. Monodromy and asymptotics of integrals}. 
Monographs in Mathematics, 83. Birkh\"auser Boston, Inc., Boston, MA,
1988.

\bibitem{BM}
B.~Bakalov, T.~Milanov.
\textit{$\mathcal{W}$-constraints for the total descendant potential of a simple singularity.}
Compositio Math. 149 (2013), no. 5, 840--888.

\bibitem{BE}{V.~Bouchard and B.~Eynard.}
\emph{Think globally, compute locally.}
J. of High Energy Phys. (2013), no. 2, Article 143.


\bibitem{DSa1}{A. Douai and C. Sabbah.}
\emph{Gauss-Manin systems, Brieskorn lattices and Frobenius
  structures (I).}
Ann. Inst. Fourier, vol. 53, no. 4(2003): 1055--1116.

\bibitem{DSa2}{A. Douai and C. Sabbah.}
\emph{Gauss-Manin systems, Brieskorn lattices and Frobenius
  structures (II). } 
Frobenius manifolds (Quantum cohomology and singularities), Aspects of Mathematics, vol. E36, Vieweg(2004): 1--18


\bibitem{Du}
B.~Dubrovin.
\textit{Geometry of 2D topological field theories}. 
In: ``Integrable systems and quantum groups'' 
(Montecatini Terme, 1993), 120--348, Lecture Notes
in Math., 1620, Springer, Berlin, 1996.

\bibitem{BOSS}{P.~Dunin--Barkowski, N.~Orantin, S.~Shadrin, and  L.~Spitz.}
{\em Identification of the Givental formula with the spectral curve topological recursion procedure.} 
Comm. in Math. Phys. 328 (2014), no. 2, 669--700.


\bibitem{Du2}{B.~Dubrovin.}
\emph{
Painlev\'e transcendents in two dimensional topological field theory.}
arXiv: 9803.107

\bibitem{DNOPS1}{P.~Dunin--Barkowski, P.~Norbury, N.~Orantin,
    A.~Popolitov, and S.~Shadrin.}
{\em Dubrovin's superpotential as a global spectral curve.} 
arXiv: 1509.06954.

\bibitem{DNOPS2}{P.~Dunin--Barkowski, P.~Norbury, N.~Orantin,
    A.~Popolitov, and S.~Shadrin.}
{\em Primary invariants of Hurwitz Frobenius manifolds.} 
arXiv: 1605.07644.

\bibitem{Ey}{B. Eynard.}
\emph{Invariants of spectral curves and intersection theory of moduli
  spaces of complex curves.} 
Comm. Numb. Theory and Phys., vol. 8, no. 3(2014): 541--588.

\bibitem{EO}{B. Eynard and N. Orantin.}
\emph{Invariants of algebraic curves and topological expansion.} 
Comm. Numb. Theory and Phys., vol. 1, no. 2(2007): 347--452.

\bibitem{FLZ}{B. Fang, C.-C. M. Liu, Z. Zong.}
\emph{ The Eynard-Orantin
recursion and equivariant mirror symmetry for the projective line.}
 Geometry \& Topology, vol. 21, no. 4 (2017): 2049–2092.

\bibitem{Fu}{W. Fulton.}
\textit{Hurwitz schemes and the irreducibility of moduli of algebraic
  curves.}
Ann. of Math., vol. 90, no. 3(1969), 542--575

\bibitem{G1}
A.~Givental. 
\textit{Semisimple Frobenius structures at higher genus}. 
Internat. Math. Res. Notices, vol. 23(2001), 1265--1286.

\bibitem{G2}
A.~Givental. 
\textit{Gromov--Witten invariants and quantization of quadratic Hamiltonians}. 
Mosc. Math. J. vol. 1(2001), 551--568. 

\bibitem{GR}
H.~Grauert and R. Remmert.
\emph{Coherent analytic sheaves.}
Springer--Verlag, Berlin Heidelberg, 1984.

\bibitem{He}{C.~Hertling.}
\emph{Frobenius Manifolds and Moduli Spaces for Singularities.}
Cambridge Tracts in Mathematics, 151. Cambridge University Press,
Cambridge, 2002. x+270 pp. 

\bibitem{Iv}{B.~Iversen}
\emph{Cohomology of sheaves.} 
Springer--Verlag, Berlin Heidelberg, 1986.

\bibitem{LZ}{Si-Qi Liu and Y. Zhang.}
\emph{Uniqueness theorem of $\mathcal{W}$-constraints for simple
  singularities.} 
Lett. Math. Phys., Vol. 103, no. 12(2013): 1329–1345.

\bibitem{ML}{T. Milanov and D. Lewanski.}
\textit{$\mathcal{W}$-algebra constraints and topological recursion for
$A_N$-singularity.}
International Journal of Math., Vol. 27, no. 13(2016) 1650110 (21
pages). 

\bibitem{Mi1}{T.~Milanov.}
\emph{
The Eynard--Orantin recursion for the total ancestor potential.}
Duke Math. J. 163 (2014), no. 9, 1795--1824.

\bibitem{Mi2}{T.~Milanov}
\emph{
The Eynard--Orantin recursion for simple singularities.}
Commun. Number Theory Phys., Vol. 9, no. 4(2015): 707--739.

\bibitem{Ra}{H. Rauch.}
\emph{Weierstrass points, branch points, and moduli of Riemann
  surfaces.}
Comm. on Pure and Appl. Math., vol. 12(1959), 543--560.

\bibitem{Sa}{K.~Saito.}
{\em Primitive forms for a universal unfolding of a function with an isolated critical point.} 
J. Fac. Sci. Univ. Tokyo Sect. IA Math. 28 (1981), no. 3, 775-792 (1982).

\bibitem{SaT}{K.~Saito and A.~Takahashi.}
{\em From primitive forms to Frobenius manifolds.}
From Hodge theory to integrability and TQFT tt*-geometry, 31-48,
Proc. Sympos. Pure Math., 78, Amer. Math. Soc., Providence, RI, 2008. 

\bibitem{MSa}{M. Saito.}
\emph{On the structure of Brieskorn lattice.} Ann. Inst. Fourier,
vol. 39, no. 1(1989): 27--72.


\bibitem{Shr}{V.~Shramchenko.}
\emph{Deformations of Frobenius structures on Hurwitz spaces.}
Internat. Math. Res. Notices, vol. 6(2005), 339--387.

\bibitem{Te}
C.~Teleman,
\emph{The structure of 2D semi-simple field theories.} 
Invent. Math., 188 (2012), no. 3, 525--588.

\end{thebibliography}

\end{document}